\documentclass[makeidx,hyperref,reqno]{amsart}
\usepackage{mathtools}
\usepackage{array}
\usepackage{todonotes}
\usepackage{xspace}
\usepackage{MnSymbol}
\usepackage{enumitem}
\usepackage{xfrac}
\usepackage{url}
\usepackage{tikz}
\usetikzlibrary{decorations,matrix,arrows}
\tikzset{closed/.style={decoration={markings,mark=at position 0.5 with {\node[transform shape,xscale=.8,yscale=.4]{/};} },postaction=decorate},open/.style={decoration={markings, mark=at position 0.5 with {\node[transform shape,scale=.7]{$\circ$};}},postaction=decorate}}
\definecolor{Light-Blue}	{rgb}{0.04,0.3,1}
\definecolor{Dark-Blue}		{rgb}{0.00,0.00,0.50}
\definecolor{Dark-Green}	{rgb}{0.00,0.50,0.00}
\definecolor{Dark-Gray}		{rgb}{0.5,0.5,0.5}
\definecolor{Dark-Red}		{rgb}{0.75,0,0}
\usepackage[bookmarks=true,bookmarksopen=true,colorlinks=true,breaklinks=true,backref=page,hyperindex=true,hyperfootnotes=true,frenchlinks=true]{hyperref}
\hypersetup{linkcolor=Dark-Blue,urlcolor=blue,citecolor=Dark-Green,anchorcolor=red,filecolor=cyan}
\renewcommand*{\backrefalt}[4]{{\tiny \ifcase #1 \ \or \textcolor{red}{$\uparrow$}c.p.~#2 \else \textcolor{red}{$\uparrow$}c.pp.~#2 \fi}}
\setlist[description]{leftmargin=\parindent,labelindent=\parindent}
\newcommand{\mathematician}[2]{\def#1{{#2}\xspace}}
\mathematician{\bondarko}			{Bondarko}
\mathematician{\borisov}			{Borisov}
\mathematician{\bousfield}			{Bousfield}
\mathematician{\calabi}				{Calabi}
\mathematician{\chow}				{Chow}
\mathematician{\dejong}				{De\ Jong}
\mathematician{\deligne}			{Deligne}
\mathematician{\eilenberg}			{Eilenberg}
\mathematician{\euler}				{Euler}
\mathematician{\fano}				{Fano}
\mathematician{\fermat}				{Fermat}
\mathematician{\fraenkel}			{Fraenkel}
\mathematician{\frobenius}			{Frobenius}
\mathematician{\gabber}				{Gabber}
\mathematician{\gabriel}			{Gabriel}
\mathematician{\galois}				{Galois}
\mathematician{\gersten}			{Gersten}
\mathematician{\gillet}				{Gillet}
\mathematician{\grothendieck}		{Grothendieck}
\mathematician{\guletskii}			{Guletski\u \i}
\mathematician{\hasse}				{Hasse}
\mathematician{\hironaka}			{Hironaka}
\mathematician{\hirzebruch}			{Hirzebruch}
\mathematician{\hodge}				{Hodge}
\mathematician{\ivorra}				{Ivorra}
\mathematician{\kan}				{Kan}
\mathematician{\kapranov}			{Kapranov}
\mathematician{\kelly}				{Kelly}
\mathematician{\kleiman}			{Kleiman}
\mathematician{\krull}				{Krull}
\mathematician{\larsen}				{Larsen}
\mathematician{\lunts}				{Lunts}
\mathematician{\maclane}			{Mac Lane}
\mathematician{\morel}				{Morel}
\mathematician{\nagata}				{Nagata}
\mathematician{\nisnevich}			{Nisnevich}
\mathematician{\noetherian}			{Noetherian}
\mathematician{\poincare}			{Poincar\'e}
\mathematician{\quillen}			{Quillen}
\mathematician{\roendigs}			{R\"ondigs}
\mathematician{\serre}				{Serre}
\mathematician{\soule}				{Soul\'e}
\mathematician{\tate}				{Tate}
\mathematician{\voevodsky}			{Voevodsky}
\mathematician{\waldhausen}			{Waldhausen}
\mathematician{\weil}				{Weil}
\mathematician{\yau}				{Yau}
\mathematician{\yoneda}				{Yoneda}
\mathematician{\zakharevich}		{Zakharevich}
\mathematician{\zermelo}			{Zermelo}
\newcommand {\eg}		{\textit{e.g.}\xspace}				

\newcommand {\Etale}	{\'{E}tale\xspace}
\newcommand {\cf}		{\textit{cf.}\xspace}				
\newcommand {\ie}		{\textit{i.e.}\xspace}				
\let\iff\undefined
\newcommand {\iff}		{if and only if\xspace}
\newcommand {\resp} 	{resp\xspace}
\newcommand {\two}		{two-out-of-three property\xspace}
\makeatletter
\newcommand\@dotsep{4.5}
\def\@tocline#1#2#3#4#5#6#7{\relax
	\ifnum #1>\c@tocdepth 
	\else
	\par \addpenalty\@secpenalty\addvspace{#2}%
	\begingroup \hyphenpenalty\@M
	\@ifempty{#4}{%
		\@tempdima\csname r@tocindent\number#1\endcsname\relax
	}{%
		\@tempdima#4\relax
	}%
	\parindent\z@ \leftskip#3\relax \advance\leftskip\@tempdima\relax
	\rightskip\@pnumwidth plus4em \parfillskip-\@pnumwidth
	#5\leavevmode\hskip-\@tempdima #6\nobreak\relax
	\leaders\hbox{$\m@th
		\mkern 4 mu\hbox{.}\mkern 4
		mu$}
	\hfil\hbox to\@pnumwidth{\@tocpagenum{#7}}\par
	\nobreak
	\endgroup
	\fi
}
\def\l@chapter{\@tocline{0}{8pt plus1pt}{0pt}{}{}}
\def\l@section{\@tocline{1}{0pt}{0pt}{5pc}{}}
\def\l@subsection{\@tocline{2}{0pt}{20pt}{5pc}{}}
\def\paragraph{\@startsection{paragraph}{4}{\z@}{2.25ex \@plus1ex \@minus.2ex}{-2em}{\normalfont\normalsize\itshape}}
\def\subparagraph{\@startsection{subparagraph}{5}{\z@}{1.75ex \@plus1ex \@minus .2ex}{-1em}	{\normalfont\normalsize\itshape}}
\@addtoreset{subsection}{section}
\@addtoreset{subsubsection}{subsection}
\@addtoreset{paragraph}{subsubsection}
\@addtoreset{subparagraph}{paragraph}
\newcommand{\prebibitem}[3][]%
{%
	\ifthenelse{\equal{#3}{Yes}}%
	{%
		\ifthenelse{\equal{#1}{}}%
		{%
			\@skiphyperreftrue\H@item\@skiphyperreffalse%
			\Hy@raisedlink{%
				\hyper@anchorstart{cite.#2\@extra@b@citeb}\relax\hyper@anchorend%
			}%
			\if@filesw%
			\begingroup%
			\let\protect\noexpand%
			\immediate\write\@auxout{%
				\string\prebibcite{#2}{\the\value{\@listctr}}%
			}%
			\endgroup%
			\fi%
			\ignorespaces%
		}%
		{%
			\@skiphyperreftrue%
			\H@item[%
			\ifx\Hy@raisedlink\@empty%
			\hyper@anchorstart{cite.#2\@extra@b@citeb}%
			\@BIBLABEL{#1}%
			\hyper@anchorend%
			\else%
			\Hy@raisedlink{%
				\hyper@anchorstart{cite.#2\@extra@b@citeb}\hyper@anchorend%
			}%
			\@BIBLABEL{#1}%
			\fi%
			\hfill%
			]%
			\@skiphyperreffalse%
			\if@filesw%
			\begingroup%
			\let\protect\noexpand%
			\immediate\write\@auxout{%
				\string\prebibcite{#2}{#1}%
			}%
			\endgroup%
			\fi%
			\ignorespaces%
		}%
	}%
	{\bibitem[#1]{#2}}%
}%
\def\MR[#1 (#2)]{\ifthenelse{\equal{#2}{}}{MR\href{http://www.ams.org/mathscinet-getitem?mr=#1}{#1}}{MR\href{http://www.ams.org/mathscinet-getitem?mr=#1}{#1 (#2)}}}
\def\doi[#1]{DOI:\href{http://dx.doi.org/#1}{#1}}
\def\MathOverflow[#1]{MathOverflow.net:\href{http://mathoverflow.net/#1}{#1}}
\def\MO[#1]{MO:\href{http://mathoverflow.net/questions/#1}{#1}}
\def\MathStackExchange[#1]{Math.StackExchange.com:\href{http://math.stackexchange.com/#1}{#1}}
\def\MSE[#1]{MSE:\href{http://math.stackexchange.com/questions/#1}{#1}}
\def\HSM[#1]{HSM:\href{http://hsm.stackexchange.com/questions/#1}{#1}}
\def\documenta[#1]{Documenta:\href{https://www.math.uni-bielefeld.de/documenta/#1}{#1}}
\def\chicago[#1]{UCP:\href{http://press.uchicago.edu/ucp/books/book/chicago/#1}{#1}}
\def\numdam[#1]{Numdam:\href{http://www.numdam.org/item?id=#1}{#1}}
\def\EuDML[#1]{{Eu}DML:\href{https://eudml.org/#1}{#1}}
\def\karchive[#1]{K-theory archive:\href{http://www.math.uiuc.edu/K-theory/#1}{#1}}
\def\MMJ[#1]{MMJ:\href{http://www.ams.org/distribution/mmj/#1}{#1}}
\def\ICM[#1]{ICM:\href{http://www.mathunion.org/ICM/ICM#1}{#1}}
\def\arxiv[#1]{arXiv:\href{http://arxiv.org/abs/#1}{#1}}
\def\wiki[#1]{wikipedia:\href{https://en.wikipedia.org/wiki/#1}{#1}}
\def\hal[#1]{hal-\href{https://hal.archives-ouvertes.fr/hal-#1}{#1}}
\def\euclid[#1]{Euclid:\href{http://projecteuclid.org/euclid.#1}{#1}}
\def\PUP[#1]{PUP title:\href{http://www.pupress.princeton.edu/titles/#1.html}{#1}}
\def\TAC[#1]{TAC:\href{http://www.tac.mta.ca/tac/#1}{#1}}
\def\ICTP[#1]{ICTP:\href{http://users.ictp.it/~pub_off/lectures/lns#1}{#1}}
\def\link[#1]{\href{#1}{#1}}
\def\claymath[#1]{Claymath:\href{http://www.claymath.org/library/#1}{#1}}
\def\bookstore[#1]{AMS bookstore:\href{http://bookstore.ams.org/#1}{#1}}
\def\youtube[#1]{YouTube:\href{https://www.youtube.com/watch?v=#1}{#1}}
\def\SP[#1]{SP:\href{http://stacks.math.columbia.edu/tag/#1}{#1}}
\DeclareSymbolFont	{stmry}			{U}{stmry}{m}{n}		 
\DeclareSymbolFont	{cmrOML}		{OML}{cmr}{m}{it}
\DeclareSymbolFont	{rsfs}			{U}{rsfs}{m}{n}
\DeclareSymbolFontAlphabet			{\mathscr}{rsfs}
\DeclareSymbolFont	{Operators} 	{OT1}{txtt}{m}{n}
\SetSymbolFont		{Operators} 	{bold}{OT1}{txtt}{b}{n}
\DeclareSymbolFont	{categories}	{OT1}{txss}{m}{n}
\SetSymbolFont		{categories} 	{bold}{OT1}{txss}{b}{n}
\DeclareSymbolFont	{ordinary}		{OT1}{lmr}{m}{n}
\SetSymbolFont		{ordinary}		{bold}{OT1}{lmr}{bx}{n}
\newcommand{\newOperator}{\@ifstar{\@declmathOp m}{\@declmathOp o}}
\long\def\@declmathOp#1#2#3{\@ifdefinable{#2}{\DeclareRobustCommand{#2}{\qOpnname\newmcodes@#1{#3}}}}
\DeclareRobustCommand{\qOpnname}[3]{\mathop{#1\kern\z@\Op@font#3}\csname n#2limits@\endcsname}%
\newcommand{\newBinaryOp}{\@ifstar{\@declmathbin m}{\@declmathbin o}}
\long\def\@declmathbin#1#2#3{\@ifdefinable{#2}{\DeclareRobustCommand{#2}{\qbinname\newmcodes@#1{#3}}}}
\DeclareRobustCommand{\qbinname}[3]{\mathbin{#1\kern\z@\binary@font#3}}
\newcommand{\newRelation}{\@ifstar{\@declmathrel m}{\@declmathrel o}}
\long\def\@declmathrel#1#2#3{\@ifdefinable{#2}{\DeclareRobustCommand{#2}{\qrelname\newmcodes@#1{#3}}}}
\DeclareRobustCommand{\qrelname}[3]{\mathrel{#1\kern\z@\relation@font#3}}
\newcommand{\newSymbol}{\@ifstar{\@declmathord m}{\@declmathord o}}
\long\def\@declmathord#1#2#3{\@ifdefinable{#2}{\DeclareRobustCommand{#2}{\qordname\newmcodes@#1{#3}}}}
\DeclareRobustCommand{\qordname}[3]{{#1\kern\z@#3}}
\newcommand{\newIndex}{\@ifstar{\@declmathind m}{\@declmathind o}}
\long\def\@declmathind#1#2#3{\@ifdefinable{#2}{\DeclareRobustCommand{#2}{\qindname\newmcodes@#1{#3}}}}
\DeclareRobustCommand{\qindname}[3]{{#1\kern\z@\ordinary@font#3}}
\newcommand{\newCategory}{\@ifstar{\@declmathcat m}{\@declmathcat o}}
\long\def\@declmathcat#1#2#3{\@ifdefinable{#2}{\DeclareRobustCommand{#2}{\qcatname\newmcodes@#1{#3}}}}
\DeclareRobustCommand{\qcatname}[3]{\mathord{#1\kern\z@\categories@font#3}}
\@onlypreamble\newOperator
\@onlypreamble\newBinaryOp
\@onlypreamble\newRelation
\@onlypreamble\newSymbol
\@onlypreamble\newIndex
\@onlypreamble\newCategory
\def\Op@font		{\mathgroup\symOperators}
\def\binary@font	{\mathgroup\symOperators}
\def\relation@font	{\mathgroup\symOperators}
\def\ordinary@font	{\mathgroup\symoperators}
\def\categories@font{\mathgroup\symcategories}
\let\newSet\newOperator
\makeatother
\makeindex
\newcommand{\indx}[1]{\index{#1}{\em #1}}
\let\mathdsr\undefined	\newcommand{\mathdsr}[1]	{\text{\usefont{U}{dsrom}{m}{n}#1}}
\let\mathsf\undefined	\newcommand{\mathsf}[1]		{\text{\usefont{OT1}{txss}{m}{n}#1}}
\let\mathbb\undefined	\newcommand{\mathbb}[1]		{\text{\usefont{U}{msb}{m}{n}#1}}
\let\mathfrak\undefined	\newcommand{\mathfrak}[1]	{\text{\usefont{U}{euf}{m}{n}#1}}
\let\mathrm\undefined	\newcommand{\mathrm}[1]		{\text{\usefont{OT1}{lmr}{m}{n}#1}}
\let\mathcal\undefined	\newcommand{\mathcal}[1]	{\text{\usefont{OMS}{cmsy}{m}{n}#1}}
\let\mathbf\undefined	\newcommand{\mathbf}[1]		{\text{\usefont{OT1}{cmr}{bx}{n}#1}}
\let\mathbftt\undefined	\newcommand{\mathbftt}[1]	{\text{\usefont{OT1}{txtt}{b}{n}#1}}
\newcommand \bbS {\mathbb S}
\newcommand \bfc {\mathbf c}

\newcommand \bfw {\mathbf w}

\newcommand \caK {\mathcal K}
 
\newcommand \dsA {\mathdsr A}
\newcommand \dsF {\mathdsr F}
\newcommand \dsN {\mathdsr N}
\newcommand \dsP {\mathdsr P}
\newcommand \dsQ {\mathdsr Q}
\newcommand \dsZ {\mathdsr Z}
\let\frm\undefined
\newcommand \frm {\mathfrak m}
\newcommand \frU {\mathfrak U}
\newcommand \rmC {\mathrm C}

\newcommand \scA {\mathscr A} 
 
\newcommand \scC {\mathscr C}
\newcommand \scD {\mathscr D}
\newcommand \scE {\mathscr E} 
 
\newcommand \scI {\mathscr I}
\newcommand \scJ {\mathscr J}

\newcommand \scP {\mathscr P}

\newcommand \scS {\mathscr S} 

\newcommand \scU {\mathscr U}
\newcommand \scV {\mathscr V}
\newcommand \scX {\mathscr X}
\newcommand \scY {\mathscr Y}
\let\1\undefined
\let\AA\undefined
\let\deg\undefined
\let\dim\undefined
\let\Oldint\int
\let\int\undefined
\let\lim\undefined
\let\OldWedge\wedge
\let\SS\undefined	
\let\wedge\undefined
\newBinaryOp	\wedge	{\OldWedge}
\newCategory 	\CHM 	{CHM}										
\newCategory	\CRing	{CRing}										
\newCategory	\CAT	{CAT}										
\newCategory	\Comp	{Comp}
\newCategory	\deloop	{\mathbftt B}								
\newCategory	\El		{El}										
\newCategory	\fdNoe	{Noe^{\!\!\!\!\!^\fd}}						
\newCategory	\FSet	{FSet}										
\newCategory	\Funct	{Fun}										
\newCategory	\Mor	{Mor}
\newCategory	\PropSch{Prop}										
\newCategory	\PSh	{PSh}										
\newCategory	\Ring	{Ring}										
\newCategory	\Set	{Set}										
\newCategory	\sftSch {Sch^{\!\!\!\!^\ft}}						
\newCategory	\Shv 	{Shv}										
\newCategory	\Site	{Site}										
\newCategory	\Spt	{Spt}										
\newCategory	\Var 	{Var}										
\newCategory	\Wald	{Wald}										
\newcommand 	\ds		{\displaystyle}
\newIndex 		\eff	{eff}										
\newIndex 		\equi	{equi}										
\newIndex 		\open 	{open}										
\newIndex 		\pr 	{pr}										
\newIndex 		\prop 	{prop}										
\newIndex 		\red 	{red}
\newIndex 		\shriek	{_!}										
\newIndex		\can	{\ensuremath {can}}							
\newIndex		\cdf	{\ensuremath {cdf}}							
\newIndex		\cdh	{\ensuremath {cdh}}							
\newIndex		\cdp	{\ensuremath {cdp}}							
\newIndex		\fd		{fd}										
\newIndex		\Fr		{Fr}
\newIndex		\ft		{ft}										
\newIndex		\htop	{\ensuremath {h}}							
\newIndex		\id		{id}
\newIndex		\inv	{^{-1}}
\newIndex		\ldh	{\ensuremath {\ell{dh}}}					
\newIndex		\Nis 	{\ensuremath {Nis}}							
\newIndex		\op		{op}
\newIndex		\pre	{pre}
\newIndex		\ps		{c}
\newIndex		\rh		{\ensuremath {rh}}							%
\newIndex		\uh		{uh}				
\newOperator 	\Cone 	{Cone}
\newOperator 	\core 	{\heartsuit}
\newOperator 	\Cyl 	{Cyl}
\newOperator 	\lop	{\Omega}
\newOperator 	\Spec	{Spec}
\newOperator 	\sus	{\Sigma}
\newOperator 	\Sym 	{Sym}										
\newOperator 	\trdeg 	{tr.deg}									
\newOperator	\ash	{a}											
\newOperator	\asp	{s}											
\newOperator	\Bl		{Bl}										
\newOperator	\coker	{coker}
\newOperator	\deg	{deg}
\newOperator	\dim	{dim}
\newOperator	\Graph	{\reflectbox{$\Gamma$}}						
\newOperator	\im		{im}
\newOperator	\imDelta{\Delta}									
\newOperator	\imGraph{\Gamma}									
\newOperator	\Kc		{\caK}										
\newOperator	\Kg		{K}
\newOperator	\Lan	{Lan}
\newOperator	\Ner 	{N}											
\newOperator	\yon	{h}
\newOperator*	\colim	{colim}
\newOperator*	\lim	{lim}
\newOperator*	\int	{\Oldint}
\newRelation 	\cim	{\nhooklongrightarrow}
\newRelation 	\incl	{\hookrightarrow}
\newRelation 	\iso	{\xrightarrow{\simeq}}
\newRelation 	\oim	{\lhook\joinrel\mathrel{\mkern9mu\circ\mkern-19mu}\longrightarrow}
\newRelation 	\onto	{\twoheadrightarrow}
\newRelation	\nhooklongrightarrow	{\lhook\joinrel\mathrel{\mkern3.5mu\notarrow\mkern-3.5mu}\longrightarrow}
\newSet			\Cg		{\rmC}										
\newSet			\CH		{\mathrm{CH}}								
\newSet			\Cov	{Cov}										
\newSet			\Gal	{Gal}
\newSet			\Hom	{Hom}										
\newSet			\iHom	{\underline{\Hom}}							
\newSet			\Ob		{Ob}
\newSymbol 		\1 		{\mathdsr 1} 								
\newSymbol 		\PP		{\dsP}										
\newSymbol 		\QQ		{\dsQ}										
\newSymbol 		\rf 	{\kappa} 									
\newSymbol		\0 		{\mathbf 0}									
\newSymbol		\AA		{\dsA}
\newSymbol		\FF		{\dsF}										
\newSymbol		\fk		{k}											
\newSymbol		\NN		{\dsN}										
\newSymbol		\SS 	{\bbS}										
\newSymbol		\ZZ		{\dsZ}										
\newcommand{\Theorems}[5]{\ifthenelse{\equal{#3}{}}{\ifthenelse{\equal{#1}{}}{\begin{#4}{\em{#2}}\end{#4}}{\begin{#4}[#1]{\em{#2}}\end{#4}}}{\ifthenelse{\equal{#1}{}}{\begin{#4}\label{#5:#3}{\em{#2}}\end{#4}}{\begin{#4}[#1]\label{#5:#3}{\em{#2}}\end{#4}}}}

\newcommand{\Th}[3]			{\Theorems{#1}{#2}{#3}{theorem}{Th}}
\newcommand{\Conj}[3]		{\Theorems{#1}{#2}{#3}{conjecture}{Conj}}
\newcommand{\Cor}[3]		{\Theorems{#1}{#2}{#3}{corollary}{Cor}}
\newcommand{\Def}[3]		{\Theorems{#1}{#2}{#3}{definition}{Def}}
\newcommand{\Eg}[3]			{\Theorems{#1}{#2}{#3}{example}{Ex}}
\newcommand{\Lem}[3]		{\Theorems{#1}{#2}{#3}{lemma}{Lem}}
\newcommand{\Prop}[3]		{\Theorems{#1}{#2}{#3}{proposition}{Prop}}
\newcommand{\Rem}[3]		{\Theorems{#1}{#2}{#3}{remark}{Rem}}
\newtheorem*{Utheorem}{Theorem}
\newtheorem*{Uconjecture}{Conjecture}
\DeclareMathSymbol		{\lhook} 		{\mathrel}{cmrOML}{'54}
\DeclareMathSymbol		{\notarrow}		\mathrel{stmry}{"58}
\newcommand {\abs}[1]		{{\left|#1\right|}}
\newcommand {\stim}[1]		{\overline{#1}}
\newcommand {\andd}			{\qquad\text{ and }\qquad}
\title{Motivic Measures through Waldhausen K-Theories}
\author{Alameddin, Anwar} 
\email{\texttt{anwar.alameddin@liv.ac.uk}}
\date{\today}
\keywords{Motivic measures, Waldhausen $\Kg$-theory, Grothendieck group, \cdp-topology}
\subjclass[2010]{14E18, 18F30, 19A99, 14C35, 18F10, 14F20}
\begin{document}
\maketitle
\begin{abstract}In this paper we introduce the notion of a \cdp-functor to a Waldhausen category. We show that such functors admit extensions that satisfy the excision property, to which we associate \euler-\poincare characteristics that send the class of a proper scheme to the class of its image. As an application, we show that the \yoneda embedding gives rise to a monoidal proper-fibred Waldhausen category over Noetherian schemes of finite Krull dimensions, with canonical \cdp-functors to its fibres. 
\end{abstract}

\tableofcontents

\section{Introduction}
In this paper, we set out to investigate the connection between motivic measures and Waldhausen $\Kg$-theories by means of the completely decomposed proper (\cdp) Grothendieck topology.

\medskip

Motivic measures are connected to fundamental questions in algebraic geometry. For instance, the motivic measure of counting rational points over a finite field gives rise to the \hasse-\weil zeta function through applying it to symmetric powers, as it was first shown by Kapranov in \cite{Kap:00:ECSDTESKMG}. Also, Larsen-Lunts motivic measure, which takes value in the monoid ring of stable birational classes of algebraic varieties over a field, has important applications in birational algebraic geometry, see \cite{LL:03:MMSBG} and \cite{GS:14:FVLRPCH}. Other important questions are tackled through the universal motivic measure, which takes value in the Grothendieck ring of varieties, see \cite{NS:11:GRV} and \cite{DL:04:SRGSOAG}. This ring is not fully understood; for instance, only recently the class of the affine line was shown to be a zero divisor for a field of characteristic zero, see \cite{Bor:15:CALZDGR}.

\medskip

More generally, for a category with a set of distinguished sequences (\eg subtraction sequences, exact sequences, cofibre sequences, distinguished triangles), its Grothendieck group is the group generated by isomorphism classes of objects module splitting the sequences. It can be though of as a decategorification of the category, with respect to the considered sequences. For a category with an exact structure, \quillen introduced a higher {algebraic $\Kg$-theory}, that extends the \grothendieck group, see \cite{Qui:73:HAKT}. {That was {generalised} by \waldhausen in \cite{Wal:85:AKTS}, who defined what is now called a \waldhausen structure, to which he associated an algebraic $\Kg$-theory spectrum whose path components group coincides with its Grothendieck group. 
	
\medskip

Several motivic measures arise from (co)homology theories with proper\footnote{They are usually called cohomology theories with compact support.} support, \ie they satisfy the excision property. For instance, the \hodge measure and the \hodge characteristic arise from the mixed \hodge structure on singular cohomology with rational coefficients and {proper support} over {the complex numbers}, see \cite{Sri:02:HC}; whereas the $\ell-$adic motivic measure arises from the $\ell$-adic cohomology with {proper support}. The latter, also gives rise to the classical measure of counting rational points over a finite field through the trace formula, see \cite{Mus:13:ZFAG}.

	\medskip

	Some motivic measures arise from weak monoidal functors
	$({\sftSch/S}^\prop_\open,\times,\id_{\!_S})\to (\scC,\wedge,\1)$ that satisfy the excision property, where $\scC$ is a symmetric monoidal Waldhausen category and ${\sftSch/S}^\prop_\open$ is the category whose objects are $S$-schemes and whose morphisms are finite compositions of proper morphisms and formal inverses of open immersions, for a base scheme $S$. In fact, every weak monoidal functor $G:({\sftSch/S}^\prop_\open,\times,\id_{\!_S})\to (\scC,\wedge,\1)$ that satisfies the excision property induces a motivic measure $\mu_G:\Kg_0(\sftSch/S)\to \Kg_0(\scC)$, which sends the class of an $S$-scheme $x$ to the class of $G(x)$. 
	
	\medskip

	The aforementioned cohomology theories arise from plain cohomology theories (that do not satisfy the excision property), and both versions coincide for proper schemes, over the base. Then, it becomes natural to ask for a scheme $S$,
\begin{itemize}
	\item[] \em when does a weak monoidal functor $F:\PropSch/S\to \scC$, from the category of proper $S$-schemes to a {symmetric monoidal \waldhausen category}, admit a properly supported extension? I.e., when does $F$ give rise to a weak monoidal functor $F^\ps:{\sftSch/S}^\prop_\open\to \scC$ that satisfies the excision property and restricts to $F$ on the category $\PropSch/S$?
\end{itemize}

\medskip

On the one hand, \gillet-\soule motivic measure, that takes value in the Grothendieck group of $K_0$-motives with rational coefficients is based on sending \cdp-squares\footnote{A \cdp-square is a generalisation of a blow up square, see Definition \ref{Def:cd:topology}.} to homotopy Cartesian squares, see \cite{GS:96:DMK} and \cite{GS:09:MWCAV}. Also, the aforementioned cohomology theories send \cdp-squares of proper $S$-schemes to (homotopy) pushout squares. On the other hand, when $S=\Spec \fk$, for a field $\fk$ of characteristic zero, there exists a motivic measure to a \waldhausen $\Kg_0$-group of simplicially stable motivic spectra, introduced in \cite{Ron:16:GRVAKTS}, which relies on a presentation of the \grothendieck group of $\fk$-varieties, in which the generators are classes of smooth projective $\fk$-varieties and the relations are induced by blow up squares, recalled in Theorem \ref{Th:Bittner:Blow:up:Relations}. That led us to distinguish functors from proper $S$-schemes to a Waldhausen category that satisfy the properties:
\begin{enumerate}	[label=\textup{(\textbf{PS\arabic*})},ref=\textup{(\textbf{PS\arabic*})}]
	\item\label{itm:PS:Closed} closed immersions of proper $S$-schemes are mapped to cofibrations;
	\item\label{itm:PS:Initial} the empty $S$-scheme is mapped to a zero object; and
	\item\label{itm:PS:Proper} \cdp-squares of proper $S$-schemes are mapped to pushout squares. 
\end{enumerate}
We call a functor that satisfies these properties a \indx{\cdp-functor}, and we use \nagata's {Compactification} Theorem to show that such functors admit properly supported extensions. Sending \cdp-squares of proper $S$-schemes to pushout squares accounts for the extensions being independent of the choice of compactifications. The below theorem is our main result in that direction.
	\begin{Utheorem}[{\ref{Th:Induced:Measure}}] 
		Let $S$ be a \noetherian scheme of finite Krull dimension, and assume that $F:(\PropSch/S,\times,\id_{\!_S})\to (\scC,\wedge,\1)$ is a weak monoidal \cdp-functor to a {symmetric monoidal \waldhausen category}. Then, there exists a functor 
		\[
		F^\ps:({\sftSch/S}^\prop_\open,\times,\id_{\!_S}) \to (\scC,\wedge,\1),
		\]
		and there exists a natural isomorphism $\varphi:F\stackrel{\sim}{\Rightarrow} F^\ps_{|_{\PropSch/S}}$, such that 
		\begin{itemize}
			\item $F^\ps$ satisfies the excision property, \ie for every closed immersion $i:v\cim x$ of $S$-schemes with complementary open immersion $j:u\oim x$, the sequence 
			\[
			F^\ps(v)\stackrel{i_\shriek}{\longrightarrow} F^\ps(x) \stackrel{j^!}{\longrightarrow} F^\ps {(u)}
			\]
			is a cofibre sequence in $\scC$, where $i_\shriek\coloneqq F^\ps(i)$ and $j^!\coloneqq F^\ps(j^\op)$; and
			\item $F^\ps$ is weak monoidal, \ie $F^\ps$ is lax monoidal whose coherence morphisms
			\[
			\phi^\ps_{x,y}:F^\ps(x)\wedge F^\ps(y) \to F^\ps(x \times y)\andd \phi^\ps_{\!_S}:\1\to F(\id_{\!_S})
			\]
			are weak equivalences in $\scC$, for every $S$-schemes $x$ and $y$.
		\end{itemize}
		Therefore, there exists a motivic measure 
		\[
		\mu_{\!_F}:\Kg_0(\sftSch/S)\to \Kg_0(\scC),
		\]
		that sends the class of a proper $S$-scheme $p$ to the class of $F(p)$.
	\end{Utheorem}

In general, a Grothendieck group is thought of as a shadow of a richer structure, a $\Kg$-theory spectrum, that encodes deeper information about the category one started with. The category of $S$-schemes does not admit a \waldhausen structure, due to the lack of enough cokernels. Recently, \zakharevich introduced, in \cite{Zak:17:KTA}, the notion of an assembler, and used it to define a spectrum whose path components group coincides with the Grothendieck group of $\fk$-schemes, for a field $\fk$. Then, Campbell defined a variation of a \waldhausen structure, called a \indx{semi-\waldhausen structure}, on the category $\fk$-schemes, in which closed immersions play the role of cofibrations, resulting in an $E_\infty$-ring spectrum with the same property, see \cite{Cam:17:KTSV}. 
	
	\medskip
	For a \noetherian scheme $S$ of finite Krull dimension, we propose a $\Kg$-theory commutative ring spectrum $\Kg\big( \scC^\omega(S)\big)$ as an extension of the (modified) Grothendieck ring of $S$-schemes. The spectrum $\Kg\big( \scC^\omega(S)\big)$ is obtained through the \cdp-cosheafification of the pointed \yoneda embedding on proper $S$-schemes, and it comes with a canonical \cdp-functor $\underline{\yon}:\PropSch/S \to \scC^\omega(S)$ to its presenting Waldhausen category.

	\begin{Uconjecture}[{\ref{Conj:Universal:Measure:Modified}}]
	The canonical motivic measure $\mu_{\!_\uh}:\Kg_0(\sftSch/S)\to \Kg_0^\uh(\sftSch/S)$ that takes value in the modified Grothendieck ring of $S$-schemes, is isomorphic to the motivic measure $\mu_{\!_\cdp}:\Kg_0(\sftSch/S)\to \Kg_0\big(\scC^\omega(S)\big)$, given in \eqref{K:Schemes}.
	\end{Uconjecture}
	In fact, the spectrum $\Kg\big(\scC^\omega(S)\big)$ arises from a fibre of a monoidal proper-fibred Waldhausen category over Noetherian schemes, as in \S.\ref{Fibred:Noetherian}.
	
	\medskip
	
	Investigating the geometric information encoded in the higher $\Kg$-groups of the spectrum $\Kg\big(\scC^\omega(S)\big)$ and validating Conjecture \ref{Conj:Universal:Measure:Modified} will be considered in a future work. 
	
	\subsection{Outline}Sections \ref{Ch:MM} and \ref{Ch:KT} recall basics of motivic measures and Waldhausen $\Kg$-theory, respectively.
	
	\medskip
	
	Section \ref{Sec:Compactification} begins with an account on compactifications, needed to extend \cdp-functors. {Afterwards,} we prove the existence of properly supported extensions for \cdp-functors. Then, we provide a brief outline on how to compactify functors that do not satisfy the properties \ref{itm:PS:Initial}-\ref{itm:PS:Proper}. That is applied to the \yoneda embedding, in \S.\ref{Compactified:Yoneda}, to obtain a monoidal proper-fibred Waldhausen category over Noetherian schemes, with canonical \cdp-functors to its fibres. 
	
	\medskip
	The paper assumes the reader's familiarity with basics of category theory, as in \cite{ML:98:CWM}. In the Appendix \ref{Cat:GroSite}, we briefly recall some of the needed results on Grothendieck sites. Also, we study representable sheaves in some non-canonical topologies, which are utilised in section \ref{Compactified:Yoneda}.

\subsection{Conventions and Notations}\label{Notations}
Throughout this paper, all \indx{schemes} are assumed to be Noetherian of finite Krull dimensions and separated over the ring of integers\footnote{One may alternatively assume that relative schemes are separated over the base.}, and $\fdNoe$ denotes the category of such schemes and their morphisms. In particular, all morphisms in $\fdNoe$ are separated. For a scheme $S$ in $\fdNoe$,
\begin{itemize}
\item an \indx{$S$-scheme} refers to a finite type morphism of schemes to $S$;
\item $\sftSch/S$ denotes the category of $S$-schemes and their morphisms;
\item $\PropSch/S$ denotes the full subcategory in $\sftSch/S$ of proper $S$-schemes; and
\item $\Var/S$ denotes the full subcategory in $\sftSch/S$ of reduced $S$-schemes, where an $S$-scheme is said to be \indx{reduced} when its underlying scheme is reduced. An object in $\Var/S$ is called an \indx{$S$-variety}.
\end{itemize}
We use small Latin letters to denote $S$-schemes, and the corresponding capital letters to denote their underlying schemes.

\medskip

For a set $\scC$ of objects in $\sftSch/S$, and for sets $\scP$ and $\scI$ of morphisms in $\sftSch/S$ that are closed under compositions and contain the isomorphisms of $S$-schemes, we denote the subcategory in $\sftSch/S$ whose objects belong to $\scC$ and whose morphisms belong to $\scP$ by $\scC^\scP$, whereas the category $(\scC^\scI)^\op$ is denoted by $\scC_\scI$. Also, we denote the category whose objects belong to $\scC$ and whose morphisms are finite compositions of morphisms in $\scP$ and formal inverses of morphisms in $\scI$ by $\scC^\scP_\scI$.

\medskip

Following \maclane's proposal in \cite{ML:69:OUFCT}, we assume and fix a model for \zermelo-\fraenkel set theory with the axiom of choice (\textbf{ZFC}), in which {we} assume and fix an uncountable {\index{Universe!Grothendieck --}\em \grothendieck universe} {$\frU$}, see \cite[Expos\'e I.\S.0]{SGA4}. Then, a {\index{Set}\em set} refers to an object in the assumed \textbf{ZFC} model, a {\index{Set!small -}\em small set} refers to an element in $\frU$, and a {\em class} refers to a subset in $\frU$. We denote by
\begin{itemize}
\item $\Set$ the category of small sets and maps between them; and
\item $\CAT$ the category of large categories and functors between them.
\end{itemize}
Reader interested in motivations for the adopted foundations may consult \cite{ML:69:OUFCT} and \cite{Shu:08:STCT}.

\begin{description}
	\item[Acknowledgements] I am thankful to Vladimir \guletskii for his directions during my PhD study and for introducing me to the subjects that resulted in this paper. Also, I would like to thank Oliver E. Anderson for useful discussions about schemes and the \htop-topology. I am thankful for the Department of Mathematical Sciences at the University of Liverpool for the financial support during my PhD study, and for Qian Wang for the financial support while working on this paper.
\end{description}

\section{Motivic Measures}\label{Ch:MM}
Throughout this section, we fix a scheme $S$ in $\fdNoe$. A \indx{generalised \euler-\poincare characteristic} over $S$ with values in a group $(G,+)$ is a map $\chi:\Ob(\sftSch/S)\to G$ that is invariant under isomorphisms and respects the \indx{scissors relations}, \ie $\chi(x)=\chi(z)+\chi(u)$ if there exists a closed immersion $z\cim x$ of $S$-schemes with complementary open immersion $u\oim x$, see \cite[p.73]{Mus:13:ZFAG} and \cite[p.5]{DL:01:GASAV}. A \indx{motivic measure} over $S$ with values in a ring $(R,+,\cdot)$ is a generalised \euler-\poincare characteristic $\mu:\Ob(\sftSch/S)\to R$ 
with value in $(R,+)$ that respects the Cartesian product of $S$-schemes, \ie for $S$-schemes $x$ and $y$, one has $\mu(x\times y)=\mu(x)\cdot \mu(y)$. In particular, when $\mu$ is surjective and $(R,+,1)$ is unitary, one has $\mu(\id_{\!_S})=1_{\!_R}$.

\medskip

The scissors relations imply that $\mu(\emptyset_{\!_S})=0_G$, and $\mu(z)=\mu(x)$, if there exits a surjective closed immersion $z\cim x$ of $S$-schemes. Therefore, one may equivalently define motivic measures over $S$ to be maps from $\Ob(\Var/S)$ to rings that are invariant under isomorphisms and respect the scissors relations and the Cartesian product in $\Var/S$.

\Eg{}{Let $\FF_q$ be a finite field with $q$ elements. Then, there exists a motivic measure $\mu_\#:\Ob(\sftSch/\FF_q)\to \ZZ$, given for an $\FF_q$-scheme $X$ by the cardinality of the set of $\FF_{q}$-points in $X$, see \cite[Prop.2.1, Rem.2.2 and Rem.2.3]{Mus:13:ZFAG}. It is called the \indx{motivic measure of counting rational points} over $\FF_q$.
}{Counting:Rational:Points}

\Eg{}{For a field $\fk$, H. \gillet and C. \soule used \dejong's alterations of singularities {\cite{dJ:97:FCA}} to define a \indx{covariant weight complex} functor $W:\Var^\prop/\fk\to \Kc^b(\CHM^\eff_\QQ(\fk))$, where $\Kc^b(\CHM^\eff_\QQ(\fk))$ is the bounded homotopy category of complexes in effective \chow motive with rational coefficient over $\fk$. In fact, that was obtained in {greater} generality. Then, using the $\Kg$-theory of coherent sheaves, {they} showed {in} \cite[Th.5.9 and Cor.5.13]{GS:09:MWCAV} that the covariant weight complex functor induces a motivic measure to the Grothendieck group of effective \chow motive with rational coefficients over $\fk$, which sends the class of a smooth projective $\fk$-variety to the class of its effective \chow motive. For a field of characteristic zero, that was obtained using \hironaka's resolution of singularities and \gersten complexes in \cite{GS:96:DMK}.
}{}

The \indx{\grothendieck group of $S$-schemes}, denoted by $\Kg_0(\sftSch/S)$, is the abelian group generated by isomorphism classes of $S$-schemes module the scissors relations, \ie module the relations
\[
\left\lbrace
\begin{tabular}{l|l}
[x]=[z]+[u]&{\tiny \text{there exists a closed immersion }$z\cim x$\text{ of }$S$-\text{schemes, with complementary open immersion }$u\oim x$}
\end{tabular}
\right\rbrace,
\]
where $[y]$ denotes the isomorphism class of an $S$-scheme $y$. We abuse notations, and let $[y]$ also denote image of the isomorphism class of $y$ in the \grothendieck group. The Cartesian product of $S$-schemes defines a commutative ring structure on the group $\Kg_0(\sftSch/S)$, whose multiplication $\cdot$ is given for $S$-schemes $x$ and $y$ by $[x]\cdot [y]\coloneqq [x\times y]$. The resulting ring is called \indx{\grothendieck ring of $S$-schemes}, in which one has $0=[\emptyset_{\!_S}]$, $1=[\id_{\!_S}]$, and $[z]=[x]$ if there exits a surjective closed immersion $z\cim x$ of $S$-schemes.

\Eg{}{
The canonical map $[-]:\Ob(\sftSch/S)\to \Kg_0(\sftSch/S)$, that sends an $S$-scheme to its class, is an initial universal motivic measure over $S$. Hence, one might abuse notation and call {any} ring homomorphism from $\Kg_0(\sftSch/S)$ {a} motivic measure over $S$.
}{}
{The} \grothendieck ring of $S$-schemes was first introduced by \grothendieck in a letter to \serre in 1964. Yet, it was {not} until 2002, when it was shown to contain zero divisors over a field of characteristic zero, see \cite{Poo:02:GRVND}. Also, in 2014, \borisov constructed two smooth \calabi-\yau varieties over the complex numbers, and showed that a multiple of {the} difference {between} their classes annihilates the class of the affine line, see \cite[Th.2.12]{Bor:15:CALZDGR}. That, in particular, answers negatively the cut-and-past question of \larsen and \lunts, as in \cite[Question 1.2]{LL:03:MMSBG}.

\paragraph{\texorpdfstring{\grothendieck}{Grothendieck} Ring of Varieties in Characteristic Zero}
For a field $\fk$ of characteristic zero, the ring $\Kg_0(\sftSch/\fk)$ admits alternative presentations with a better-behaved set of generators.

\Th{}{
	Let $\fk$ be a field of characteristic zero. Then, the group $\Kg_0(\sftSch/\fk)$ is isomorphic to the abelian group generated by isomorphism classes of smooth connected projective (\resp. proper) $\fk$-varieties modulo the relations
	\begin{itemize}
		\item $[\emptyset]=0$; and 
		\item $[\Bl_Y X] -[E]=[X]-[Y]$, for every smooth connected projective (\resp. proper) $\fk$-variety $X$ and a closed smooth subvariety $Y\cim X$, where $\Bl_Y X$ is a blow-up of $X$ along $Y$ with an exceptional divisor $E$.
	\end{itemize}
}{Bittner:Blow:up:Relations}
\begin{proof}
	See \cite[Th.3.1]{Bit:04:UECVCZ}.
\end{proof}

\paragraph{The Modified \texorpdfstring{\grothendieck}{Grothendieck} Ring of Varieties}
In motivic integration, it is convenient to consider a modified version of the Grothendieck ring of $S$-schemes, see \cite{NS:11:GRV}. Also, see \cite{Har:16:QMQGRV}.

\medskip

Let $I_{\!_S}^\uh$ be the ideal in the Grothendieck ring of $S$-schemes generated by the set
\[
\left\lbrace
\begin{tabular}{l|l}
[x]-[y]&\text{there exists a universal homeomorphism }$x\to y$\text{ of }$S$-\text{schemes}
\end{tabular}
\right\rbrace.
\]
The \indx{modified Grothendieck ring of $S$-schemes}, denoted by $\Kg_0^\uh(\sftSch/S)$, is the quotient ring $ \Kg_0(\sftSch/S)/I_{\!_S}^\uh$. We denote the quotient homomorphism $\Kg_0(\sftSch/S)\to \Kg_0^\uh(\sftSch/S)$
by $\mu_{\!_\uh}$.

\Prop{}{Let $f:x\to y$ be a universal homeomorphism of $\QQ$-schemes. Then, $f$ is a piecewise isomorphism. Thus, for a $Q$-scheme $S$, the motivic measure $\mu_{\!_\uh}$ is an isomorphism.
}{}
\begin{proof}
	See \cite[Prop.3.10 and Cor.3.11]{NS:11:GRV}.
\end{proof}
The motivic measure of counting rational points over a finite field, the $\ell$-adic motivic measure over a Noetherian base, and the Euler characteristic over a field factorise through the motivic measure $\mu_\uh$, see \cite[Prop.4.13]{NS:11:GRV}. In general, it is not know if $\mu_\uh$ is an isomorphism.

\section{\texorpdfstring{\waldhausen}{Waldhausen} \texorpdfstring{$\Kg$}{K}-Theory}\label{Ch:KT}
For a category with a set of distinguished sequences (\eg exact sequences, cofibre sequences, distinguished triangles), its Grothendieck group is the abelian group generated by isomorphism classes of objects module splitting the sequences. It can be though of as a decategorification of the category, with respect to the considered sequences. For a category with an exact structure, \quillen introduced an {algebraic $\Kg$-theory}, that extends the \grothendieck group, see \cite{Qui:73:HAKT}. {That was {generalised} by \waldhausen in \cite{Wal:85:AKTS}, who defined what is now called a \waldhausen structure, to which he associated an algebraic $\Kg$-theory spectrum whose path components group coincides with its Grothendieck group.

\medskip

Let $\scC$ be a pointed category, \ie a category with a zero object $\0$. A \indx{Waldhausen structure} on $\scC$ is a pair $(\bfc\scC,\bfw\scC)$ of subcategories in $\scC$ that contain all isomorphisms in $\scC$ such that
\begin{enumerate}[label=\textup{\textbf{W\arabic*}},ref=\textup{\textbf{W\arabic*}}]
	\item\label{itm:W1}the initial morphism $\0\to U$ belongs to $\bfc\scC$, for every object $U\in \scC$; 
	\item\label{itm:W2}all pushouts along morphisms in $\bfc\scC$ exist in $\scC$, and $\bfc\scC$ is closed under pushouts; and 
	\item\label{itm:W3} for every solid commutative diagram 
\[
\begin{tikzpicture}[descr/.style={fill=white},>=angle 90,scale=2,text height=1.5ex, text depth=.25ex,row sep=4em, column sep=4em]
\node (F110) at (-1,.5) {$U$};
\node (F120) at (.5,.5) {$V$};
\node (F210) at (-1,-1) {$W$};
\node (F220) at (.5,-1) {$Z$};
\node at (.35,-.85) {$\lrcorner$};
\node (F11) at (-.5,-0.25) {$U'$};
\node (F12) at (1,-0.25) {$V'$};
\node (F21) at (-.5,-1.75) {$W'$};
\node (F22) at (1,-1.75) {$Z'$};
\node at (.85,-1.6) {$\lrcorner$};
\draw[dashed,->,font=\scriptsize]
(F120) edge node[descr,pos=.5,sloped] {} (F220)
(F12) edge (F22)
;
\draw[dashed,>->,font=\scriptsize]
(F210) edge node[descr,pos=.28,sloped] {} (F220)
(F21) edge (F22)
;
\draw[->,font=\scriptsize]
(F11) edge (F21)
(F110) edge (F210)
(F110) edge (F11)
(F210) edge (F21)
(F120) edge (F12)
;

\draw[->,dotted,font=\scriptsize]
(F220) edge (F22)
;
\draw[>->,font=\scriptsize]
(F11) edge node[descr]{$i'$}(F12)
(F110) edge node[descr]{$i$}(F120)
;
\end{tikzpicture}
\]
	in $\scC$, in which $i,i'$ belong to $\bfc\scC$ and the diagonal solid morphisms belong to $\bfw\scC$, the induced morphism $Z\to Z'$ belongs to $\bfw\scC$, for $Z:=W\ds\coprod_U V$ and $Z':=W'\ds\coprod_{U'} V'$.
\end{enumerate} 
For a \waldhausen structure $(\bfc\scC,\bfw\scC)$ on $\scC$, the triple $(\scC,\bfc\scC,\bfw\scC)$ is called a \indx{Waldhausen category}. and a morphism in $\bfc\scC$ (\resp. $\bfw\scC$) is called a \indx{cofibrations} (\resp. \indx{weak equivalences}) in the \waldhausen structure. A cofibrations in $\bfc\scC$ is denoted by a feathered arrow $\rightarrowtail$. Moreover, if $\bfw\scC$ also satisfies the \two then the \waldhausen category is said to be \indx{saturated}.

\medskip

When no confusion arise, we refer to the Waldhausen category $(\scC,\bfc\scC,\bfw\scC)$ by $\scC$.

\Eg{}{Let $(\scC,\tau)$ be an essentially small site. Then, the category of (finitely presented objects of) pointed $\tau$-sheaves of sets $\Shv_\tau(\scC)$ admits a canonical Waldhausen structure, in which cofibrations coincide with monomorphisms and weak equivalence coincide with isomorphisms.
}{Waldhausen:Sheaves}

\Eg{}{Let $(\scA,\scE)$ be a \quillen exact category, let $\bfc\scA$ be the set of admissible monomorphisms in $\scE$, and let $\bfw\scA$ be the set of isomorphisms in $\scA$. Then, $\big(\scA,\bfc\scA,\bfw\scA\big)$ is a \waldhausen category, see \cite[\S.1.2.9]{TT:90:HAKTSDC}.}{Waldhausen:Exact}

\waldhausen associated to an essentially small \waldhausen category $\scC$ a spectrum $\Kg(\scC)$, through what is known as the \waldhausen $S$-construction, see \cite[\S.1.3]{Wal:85:AKTS}. This {spectrum} $\Kg(\scC)$ is called the \indx{Waldhausen $\Kg$-theory} of $\scC$. For an integer $n\geq 0$, the \indx{Waldhausen $\Kg$-group} $\Kg_n(\scC)$ is the homotopy group $\Kg_n(\scC)\coloneqq \pi_n \Kg(\scC)$. In particular, the \waldhausen $\Kg_0(\scC)$ group admits a simple expression. That is, \indx{$\Kg_0(\scC)$} is isomorphic to the abelian group generated by isomorphism classes of objects in $\scC$ modulo the relations 
\begin{enumerate}
	\item $[U]=[V]$ if there exists a weak equivalence $U\to V$; and
	\item $[V]=[U]+[V/U]$ for every cofibre sequence $U\rightarrowtail V\twoheadrightarrow V/U$;
\end{enumerate}
where a \indx{cofibre sequence} in $\scC$ is a cokernel sequence of a cofibration in $\scC$, and $[X]$ denotes the isomorphism class of an object $X$ in $\scC$, see \cite[\S.1.5.6]{TT:90:HAKTSDC}.

\Eg{\eilenberg Swindle}{Let $\scC$ be an essentially small \waldhausen category that is closed under countable coproducts. Then, the spectrum $\Kg(\scC)$ is connected, \ie the group $\Kg_0(\scC)$ vanishes.}{Eilenberg:Swindle}
One is usually interested in \waldhausen categories that satisfy some finiteness conditions, which do not admit \eilenberg Swindle. Since finitely presented objects are closed under finite colimits, given a \waldhausen category $(\scC,\bfc\scC,\bfw\scC)$, the full subcategory of its finitely presented objects $\scC^c$ admits a \waldhausen structure $(\bfc\scC^c,\bfw\scC^c)$, given by the restriction of the structure $(\bfc\scC,\bfw\scC)$ to $\scC^c$, \cf \eqref{Waldhausen:Omega}.

\Eg{}{Let $\scC$ be the category of spectra of pointed simplicial sets. Then, $\Kg_0(\scC)=\0$ and $\Kg_0(\scC^c)=\ZZ$, see \cite[Prop.5.5.1]{Bon:10:WStS}.}{}

\Th{}{Let $(\scA, \scE)$ be an essentially small \quillen exact category. Then, there exists a natural homotopy equivalence between \quillen's $\Kg$-theory spectrum of $(\scA, \scE)$ and \waldhausen's $\Kg$-theory spectrum of the associated \waldhausen category, as {in Example} \ref{Ex:Waldhausen:Exact}.
}{Waldhausen:Quillen}
\begin{proof}
	See \cite[\S.1.9]{Wal:85:AKTS}.
\end{proof}

For \waldhausen categories $\scC$ and $\scD$, a functor $F:\scC\to \scD$ is said to be \indx{exact} with respect to the \waldhausen structures, if it preserves cofibrations, weak equivalences, and pushouts along cofibrations. A natural transformation $\alpha:F\Rightarrow G:\scC\to \scD$ between exact functors is said to be a \indx{weak equivalence} if the components of $\alpha$ belong to $\bfw\scD$, see \cite[p.330]{Wal:85:AKTS}. Then, an exact functor $F:\scC\to \scD$ between essentially small Waldhausen categories induces a map of spectra
$\Kg(F):\Kg(\scC)\to \Kg(\scD)$, and a weak equivalence $\alpha:F\Rightarrow G$ between exact functors induces a homotopy $\Kg(F)\Rightarrow \Kg(G)$, see \cite[Prop.1.3.1]{Wal:85:AKTS}.

\subsubsection{Symmetric Monoidal \texorpdfstring{\waldhausen}{Waldhausen} Categories}
A symmetric monoidal \waldhausen category admits a symmetric monoidal structure and a \waldhausen structure that are compatible, which endues its $\Kg$-theory spectrum with a homotopy commutative monoid structure.

\Def{}{A \indx{symmetric monoidal \waldhausen category} $(\scC,\scS)$ is a pair of a \waldhausen category $\scC$ and a symmetric monoidal structure $\scS=(\wedge,\1_{\wedge},\psi,\alpha,\lambda,\rho)$ on $\scC$ such that
	\begin{enumerate}
		\item\label{MW:Tensor} the endofunctors $X\wedge -$ and $-\wedge X$ are exact, for every $X\in \scC$; and 
		\item\label{MW:Pushout} for cofibrations $i:U\to V$ and $i':U'\to V'$ in $\scC$, the pushout product
		\[
		i \square i':U\wedge V'\coprod_{U\wedge U'}V\wedge U'\rightarrow V\wedge V'
		\]
		is a cofibration in $\scC$.
	\end{enumerate}
}{Waldhausen:Model:Structure}

\Eg{}{Let $\FSet_\bullet$ be the Waldhausen category of pointed finite sets whose cofibrations (\resp. weak equivalences) are monomorphisms (\resp. isomorphisms). Recall that Barratt-Priddy-\quillen Theorem implies $\Kg(\FSet_\bullet)\cong \SS$, where $\SS=(S^0,S^1,S^2,\ldots)$ is the sphere spectrum, see \cite[Th.8.9.3]{Rog:10:LNAKT}. The smash product endues the category $\FSet_\bullet$ with a symmetric monoidal structure, with a unite $\ast_+=(\ast\coprod \ast, \ast)$, making it into a symmetric monoidal \waldhausen category. For a pointed finite set $(X,x)$, one has $[(X,x)]=\abs{X\setminus \{x\}}\cdot [\ast_+]\in \Kg_0(\FSet_\bullet)$.
}{Waldhausen:Finite:Sets}

For an essentially small symmetric monoidal \waldhausen category $(\scC,\scS)$, there exists an exact functor of Waldhausen categories
\begin{equation}\label{Unit:Waldhausen:Upsilon}
\upsilon_\scC:\FSet_\bullet \to \scC,
\end{equation}
for which $\upsilon_\scC\big((X,x)\big)\cong \coprod_{X\setminus \{x\}} 1_{\wedge}$, for every pointed finite set $(X,x)$. Hence, there exists a map of spectra $\Kg(\upsilon_\scC):\SS\cong \Kg(\FSet_\bullet)\to \Kg(\scC)$. On the other hand, the monoidal product defines a paring $\Kg(\wedge):\Kg(\scC)\wedge \Kg(\scC)\to \Kg(\scC)$, see \cite[p.342]{Wal:85:AKTS}. Then, the coherence natural isomorphisms of the monoidal structure induce a homotopy commutative monoid structure on the \waldhausen $\Kg$-theory spectrum $\Kg(\scC)$, \ie they make $\Kg(\scC)$ into a ring spectrum, see \cite[Cor.2.8]{BM:11:DKDIAKTS}. In particular, $\Kg_0(\scC)$ is a ring and its ring characteristic $\ZZ\to \Kg_0(\scC)$ is given by the ring homomorphism $\Kg_0(\upsilon_\scC)$.

\Eg{}{Recall Example \ref{Ex:Waldhausen:Sheaves}, for an essentially small site $(\scC,\tau)$, the Waldhausen category of pointed $\tau$-sheaves of sets on $\scC$ is a symmetric monoidal \waldhausen category, whose monoidal product is given by the smash product. The unit of the symmetric monoidal structure is given by $\ast_+=(\ast \coprod \ast,\ast)$, whereas the smash product of pointed $\tau$-sheaves $(\scX,x)$ and $(\scY,y)$ is given by the pushout of the span
\[
\begin{tikzpicture}[descr/.style={fill=white}]
\node (X0) at (0,0) {$\scX\coprod \scY$};
\node (XBar0) at (4,0) {$\scX\times \scY$};
\node (X1) at (0,-1) {$\ast$};
\path[->,font=\scriptsize]
(X0) edge (X1)
(X0) edge node[above]{${{ (\id_\scX\times y)}\coprod{ (x\times \id_\scY)}}$} (XBar0)
;
\end{tikzpicture}
\]
in $\Shv_\tau(\scC)$, with the canonical base-point.
}{Waldhausen:Sheaves:Monoidal}

An exact functor $F:\scC\to \scD$ between symmetric monoidal \waldhausen categories is said to be \indx{weak monoidal} if it is lax monoidal such that the coherence morphism $F(X) \wedge_\scD F(Y)\to F(X \wedge_\scC Y)$ belongs to $\bfw\scD$ for every $X,Y\in \scC$, and so is the coherence morphism $\1_\scD\to F(\1_\scC)$. For a weak monoidal exact functor $F:\scC\to \scD$ between essentially small symmetric monoidal \waldhausen categories, the map of spectra $\Kg(F):\Kg(\scC)\to \Kg(\scD)$ is a morphism of ring spectra.

\section{Properly Supported Extensions}\label{Sec:Compactification}
In this section\footnote{This section does not require restricting Noetherian schemes to those of finite Krull dimensions.}, we show that a weak monoidal functor $F:\PropSch/S\to \scC$, to a symmetric monoidal \waldhausen category, admits a properly supported extension $F^\ps:{\sftSch/S}^\prop_\open\to \scC$, when it satisfies the properties \ref{itm:PS:Closed}-\ref{itm:PS:Proper}. In which case, it defines a motivic measure to $\Kg_0(\scC)$, given by sending the class of a proper $S$-scheme to the class of its image along $F$, see Theorem \ref{Th:Induced:Measure}.

\medskip

We begin by defining compactifications of $S$-schemes, and we show the category of compactifications to be cofiltered\footnote{Our notion of a morphism of compactifications differs from that usually used in the literature.}, as in Corollary \ref{Cor:Compactification:Cofiltered}, which is the main ingredient used to define properly supported functors on morphisms. Then, in \S.\ref{Extension:PS}, we define the extension $F^\ps$, and study its properties leading to the construction of the desired motivic measure in \S.\ref{Motivic:Measure:PS}. Finally, we describe how one may proceed when $F$ is not weak monoidal or does not satisfy the properties \ref{itm:PS:Initial}-\ref{itm:PS:Proper}.

\subsection{Compactifications}\label{Compactifications}Our notions of compactifications of $S$-schemes and their morphisms differ slightly from the notions usually used in the literature, \eg as in \cite[\S.2.2.8]{CD:13:TCMM}. These notions were essentially motivated by the argument in \cite{GS:96:DMK}, {{Corollary} \ref{Cor:Compactification:Cofiltered}, and {Definition} \ref{Def:Compact:Support:Functor:Proper}}. Then, after becoming aware of \indx{subtraction sequences} in \cite{Cam:17:KTSV}, they were altered to the current form, which both strengthens the result and simplifies the proofs.

\Def{}{
	Let $x$ be an $S$-scheme. A \indx{compactification} of $x$ is a closed immersion $i:z\cim p$ of proper $S$-schemes with complementary open immersion $j_{\!_i}:x\oim p$. Let $i:z\cim p$ and $l:w\cim q$ be a pair of compactifications of $S$-schemes $x$ and $y$, respectively. {A \indx{morphism of compactifications} $(f,g):i\to l$ is a solid commutative square
		\begin{equation}\label{Eq:Compactifications:Morphism}
		\begin{tikzpicture}[descr/.style={fill=white}]
		\node (fp) at (0,0) {$w\times_{q}p$};
		\node (X0) at (-1,1) {$z$};
		\node (XBar0) at (2,0) {$p$};
		\node (X1) at (0,-2) {$w$};
		\node (XBar1) at (2,-2) {$q$};
		\node at (.25,-.4) {$\ulcorner$};
		\path[->,font=\scriptsize]
		(XBar0) edge node[right]{$g$} (XBar1);
		\path[right hook->,font=\scriptsize]
		(X1) edge node{$\diagup$}node[below]{$l$} (XBar1);
		\path[right hook->,bend left,font=\scriptsize]
		(X0) edge node{$\diagup$}node[above]{$i$} (XBar0);
		\path[->,bend right,font=\scriptsize]
		(X0) edge node[left]{$f$} (X1);
		\path[dashed, ->,font=\scriptsize]
		(fp) edge (X1);
		\path[dashed, right hook->,font=\scriptsize]
		(fp) edge node{$\diagup$} (XBar0);
		\path[dotted, right hook ->,font=\scriptsize]
		(X0) edge node{$\diagup$} (fp);
		\end{tikzpicture}
		\end{equation}
		in $\sftSch/S$, for which the unique morphism $z\cim w\times_q p$ of $S$-schemes, induced by the universal property of pullbacks, is surjective.} 
}{Compactification}For a morphism of compactifications $(f,g):i\to l$, the morphism $f$ is uniquely determined by $g$, when it exists, as  $l$ is a monomorphism of $S$-schemes. Hence, when no confusion arise, we may denote this morphism of compactifications by $g:i\to l$.

\medskip

A compactification $i:z\cim p$ of an $S$-scheme $x$ induces a complementary open immersion $x\oim p$, which is unique up to isomorphisms, and we denote by $j_{\!_i}$. Since {open complements are closed under pullbacks} and both $z$ and $w\times_q p$ in \eqref{Eq:Compactifications:Morphism} have the same open complement in $p$, the morphism of compactifications $g:i\to l$ induces a Cartesian square
\begin{equation}\label{Eq:Compactifications:Morphism:Open}
\begin{tikzpicture}[descr/.style={fill=white}]
\node (X0) at (0,0) {$x$};
\node (XBar0) at (2,0) {$p$};
\node (X1) at (0,-2) {$y$};
\node (XBar1) at (2,-2) {$q$};
\node at (.25,-.4) {$\ulcorner$};
\path[->,font=\scriptsize]
(XBar0) edge node[right]{$g$} (XBar1)
(X0) edge node[left]{$g_{|_x}$} (X1);
\path[right hook->,font=\scriptsize]
(X0) edge node[descr]{$\!\!\!\circ\!\!\!$}node[above]{$j_{\!_i}$} (XBar0)
(X1) edge node[descr]{$\!\!\!\circ\!\!\!$}node[below]{$j_{\!_l}$} (XBar1);
\end{tikzpicture}
\end{equation}
in $\sftSch/S$. One may alternatively define the morphism of compactifications $g:i\to l$ to be the solid outer square in \eqref{Eq:Compactifications:Morphism} that induces the Cartesian square \eqref{Eq:Compactifications:Morphism:Open}.

\Rem{}{Although, the existence of the Cartesian square \eqref{Eq:Compactifications:Morphism:Open} does not imply the existence of a morphism of compactifications $(f,g):i\to l$, it defines a morphism of compactifications $i_{\!_\red}\to l$, where $i_{\!_\red}$ is the composition of $i$ with the surjective closed immersion $z_{\!_\red} \cim z$. One may be tempted to define a morphism of compactifications as a Cartesian square, without invoking the additional surjective closed immersion. However, our need to induce a morphism of compactifications from the Cartesian square \eqref{Eq:Compactifications:Morphism:Open}, to prove Proposition \ref{Prop:Compactification:Proper:Morphism} and Proposition \ref{Prop:Compactification:Refinement}, is the reason for the adopted notion of a morphism of compactifications. 
}{Cartesian:Square:Morphism:Compactifications}

Compactifications of $S$-schemes and their morphisms form a category, with the evident composition and identity maps, which we denote by $\Comp_{\!_S}$. For an $S$-scheme $x$, let $\Comp_{\!_S}(x)$ denote the subcategory in $\Comp_{\!_S}$ whose objects are compactifications of $x$ and whose morphisms are morphisms of compactifications that restrict to isomorphisms on $x$. That is, a morphism $g:i'\to i$ of compactifications of $x$ belongs to $\Comp_{\!_S}(x)$ \iff $\id_x$ is a base change in $\sftSch/S$ of $g$ along $j_{\!_i}$. The restriction imposed on morphisms in $\Comp_{\!_S}(x)$ is needed for Corollary \ref{Cor:Compactification:Cofiltered}, and for the cofibres in Remark \ref{Rem:Cofiltered:Extension} to be independent of the choice of compactifications.

\medskip

Before we proceed, we need to recall the following technical result, which we utilise on multiple occasions.

\Lem{}{Let $i:v\cim x$ be a closed immersion and $j:x\oim q$ be an open immersion of $S$-schemes, and let $i':p\cim q$ be the scheme-theoretic image of the immersion $j\circ i$. Then, the unique morphism $j':v\to p$ of $S$-schemes for which $j\circ i=i'\circ j'$ is an open immersion. Moreover, the square
	\[
	\begin{tikzpicture}[descr/.style={fill=white}]
	\node (Z) at (2,0) {$p$};
	\node (Y) at (0,-2) {$x$};
	\node (Q) at (2,-2) {$q$};
	\node (X) at (0,0) {$v$};
	\node at (.25,-.4) {$\ulcorner$};
	\path[right hook->,font=\scriptsize]
	(Z) edge node{$\diagup$} node[right]{$i'$}(Q)
	(X) edge node{$\diagup$} node[left]{$i$}(Y)
	(Y) edge node[descr]{$\!\!\!\circ\!\!\!$}node[below]{$j$} (Q)
	(X) edge node[descr]{$\!\!\!\circ\!\!\!$}node[above]{$j'$} (Z)
	;
	\end{tikzpicture}
	\]
	is Cartesian in $\sftSch/S$.
}{Open:Closed:Closed:Open}
\begin{proof}Consider the commutative diagram 
	\begin{equation}\label{Open:Closed:Closed:Open:Proof}
\begin{tikzpicture}[descr/.style={fill=white}]
\node (Z) at (2,0) {$p$};
\node (fp) at (0,0) {$x\times_q p$};
\node (Y) at (0,-2) {$x$};
\node (Q) at (2,-2) {$q$};
\node (X) at (-1,1) {$v$};
\node at (.25,-.4) {$\ulcorner$};
\path[right hook->,font=\scriptsize]
(Z) edge node{$\diagup$} node[right]{$i'$}(Q)
(Y) edge node[descr]{$\!\!\!\circ\!\!\!$}node[below]{$j$} (Q)
(fp) edge node{$\diagup$} node[right]{$\underline{i}'$}(Y)
(fp) edge node[descr]{$\!\!\!\circ\!\!\!$}node[below]{$\underline{j}$} (Z)
;
\path[->,dotted,font=\scriptsize]
(X) edge node[descr]{$\!\!\!l\!\!\!$}(fp)
;
\path[bend right,right hook->,font=\scriptsize]
(X) edge node{$\diagup$} node[left]{$i$}(Y);
\path[bend left,->,font=\scriptsize]
(X) edge node [above]{$j'$} (Z)
;
\end{tikzpicture}
\end{equation}
in $\sftSch/S$, where $l$ is the unique morphism $v\to x\times_q p$ of $S$-schemes that makes the diagram commute. Since $j\circ i$ is an immersion, so is $j'$. All the underlying schemes of the $S$-schemes in the diagram \eqref{Open:Closed:Closed:Open:Proof} are Noetherian; hence the immersion $j'=\underline{j}\circ l$ is quasi-compact, see \cite[Tags 01OX and 01T6]{stacks-project}, and it factorises in $\sftSch/S$ as an open immersion followed by a closed immersion, see \cite[Tag 01QV]{stacks-project}. Then, $j'$ is an open immersion, as $i'$ is the scheme-theoretic image of $j\circ i=i'\circ j'$. Hence, $l$ is also an open immersion. On the other hand, since $i$ and $\underline{i}'$ are closed immersions, so is $l$, which is also surjective because $i'$ the scheme-theoretic image of $j\circ i$. Therefore, $l$ is a surjective open immersion, and hence an isomorphism.
\end{proof}

\paragraph{The Category of Compactifications}
{Since} the notions of compactifications and their morphisms used here differ from those in the literature, {we} {need to prove} that the category $\Comp_{\!_S}(x)$, and certain subcategories of which, are {cofiltered}, for every $S$-scheme $x$. This is the main tool used to extend a functor $F:\PropSch/S\to \scC$ that satisfies the properties \ref{itm:PS:Closed}-\ref{itm:PS:Proper} to a functor $F^\ps:{\sftSch/S}^\prop_\open\to \scC$ that satisfies the excision property.

\medskip
Recall that a nonempty category $\scJ$ is \indx{cofiltered} if 
\begin{itemize}
	\item for every $X_0,X_1\in \scJ$ there exists a \indx{span} $X_0\leftarrow X \rightarrow X_1$ in $\scJ$; and
	\item for every parallel morphisms $f_0,f_1:X_0\rightrightarrows X_1$ in $\scJ$, there exists a \indx{refining morphism} $f:X\to X_0$ in $\scJ$ for which $f_0 \circ f=f_1\circ f$.
\end{itemize}
For every $S$-scheme $x$, we start by showing the category $\Comp_{\!_S}(x)$ to be nonempty, then Proposition \ref{Prop:Compactification:Proper:Morphism} provides the existence of the desired spans, and Proposition \ref{Prop:Compactification:Refinement} gives the refining morphisms.

\Rem{}{
	Due to \nagata's Compactification Theorem, as in \cite{Nag:62:IAVCV} and \cite{Nag:63:GIPAVCV}, every $S$-scheme $x$ admits an open immersion $j:x\oim p$ into a proper $S$-scheme $p$. Let $i_{\!_j}:z\cim p$ be a complementary closed immersion to $j$, endued with the reduced induced structure. Then, $i_{\!_j}$ is a compactification of $x$, and hence $\Comp_{\!_S}(x)\neq\emptyset$. In particular, when $p$ is a proper $S$-scheme, the category $\Comp_{\!_S}(p)$ has an initial object, namely $\emptyset_p:\emptyset_{\!_S}\cim p$.
}{Nagata}

\Prop{}{Assume that $f:x\to y$ is a proper morphism of $S$-schemes, and let $i:z\cim p$ and $l:w\cim q$ be compactifications of $x$ and $y$, respectively. Then, there exists a compactification $i':z'\cim p'$ of $x$ and morphisms of compactifications ${h'}:i'\to i$ and ${g'}:i'\to l$, such that $\id_x$ (\resp. $f$) is a base change in $\sftSch/S$ of ${h'}$ (\resp. ${g'}$) along $j_{\!_i}$ (\resp. $j_{\!_l}$), where $j_{\!_i}:x\oim p$ and $j_{\!_l}:y\oim q$ are complementary open immersions of $i$ and $l$, respectively.
}{Compactification:Proper:Morphism}
\begin{proof}
	In line with the arguments in \cite[\S.2.3,p.141]{GS:96:DMK} and {\cite[Tags 0ATU and 0A9Z]{stacks-project}}, consider the solid {commutative} diagram 
	\begin{equation}\label{Compactification:Product}
	\begin{tikzpicture}[descr/.style={fill=white}]
	\matrix(m)[matrix of math nodes,row sep=2em, column sep=2em,text height=1.5ex, text depth=.25ex]
	{x&&&&p\\
		&\imGraph_f&&\stim{\imGraph_f}\\
		&x\times y&&p\times q\\
		y&&&&q\\
	};
	\path[dotted,->,font=\scriptsize]
	(m-3-2) edge node[descr,pos=.195,sloped] {} (m-1-1)
	(m-3-2) edge (m-4-1);
	\path[dotted,->,font=\scriptsize]
	(m-2-4) edge node[left]{${{h'}}$} (m-1-5)
	(m-2-4) edge node[right]{${{g'}}$} node[descr,pos=.18,sloped] {} (m-4-5)
	;
	\path[->,font=\scriptsize]
	(m-2-2) edge node[right]{$h$}(m-1-1)
	(m-2-2) edge node[left]{$g$} (m-4-1)
	(m-3-4) edge (m-1-5)
	(m-3-4) edge (m-4-5)
	(m-1-1) edge node[left]{$f$} (m-4-1)
	;
	\path[right hook->,font=\scriptsize]
	(m-1-1) edge node[descr]{$\!\!\!\circ\!\!\!$}node[above]{$j_{\!_i}$} (m-1-5)
	(m-4-1) edge node[descr]{$\!\!\!\circ\!\!\!$}node[below]{$j_{\!_l}$} (m-4-5)
	(m-3-2) edge node[descr]{$\!\!\!\circ\!\!\!$}node[below]{$j_{\!_i}\times j_{\!_l}$} (m-3-4);
	\path[right hook->,font=\scriptsize]
	(m-2-2) edge node{$\diagup$}node[right]{${c}$}(m-3-2);
	\path[right hook->,dashed,font=\scriptsize]
	(m-2-2) edge node[descr]{$\!\!\!\circ\!\!\!$}node[above]{$j$} (m-2-4);
	\path[right hook->,dashed,font=\scriptsize]
	(m-2-4) edge node{$\diagup$}node[left]{$c'$} (m-3-4);
	\end{tikzpicture}
	\end{equation}
	in ${\sftSch/S}$, that is induced by the existence of Cartesian products in ${\sftSch/S}$ and the definition of the graph $\imGraph_f$ of $f$. Since open immersions are closed under pullbacks and compositions, the morphism $j_{\!_i}\times j_{\!_l}$ is an open immersion. {Let $h$ (\resp. $g$) be the composition of the Cartesian product projection $x\times y\onto x$ (\resp. $x\times y\onto y$) with the closed immersion $c:\imGraph_f\cim x\times y$, let} $c':\stim{\imGraph_f}\cim p\times q$ be the scheme-{theoretic} image of $(j_{\!_i}\times j_{\!_l})\circ c$, and let ${h'}$ (\resp. ${g'}$) be the composition of the Cartesian product projection $p\times q\onto p$ (\resp. $p\times q\onto q$) with the closed immersion $c':\stim{\imGraph_f}\cim p\times q$.	

	\medskip
	
	There exists an open immersion $j$ for which $(j_{\!_i}\times j_{\!_l})\circ c=c' \circ j$, by Lemma \ref{Lem:Open:Closed:Closed:Open}. The composition $h:\imGraph_f\cim x\times y\onto x$ is an isomorphism, see \cite[p.134]{EGA1}. Thus, there exists a compactification $i':z'\cim \stim{\imGraph_f}$ of $x$, where $i'$ is a complementary closed immersion to $j\circ h^\inv:x\oim \stim{\imGraph_f}$, endued with the reduced induced structure.

	\medskip
	
	Consider the commutative diagram 
	\[\label{Compactification:Product:2}
	\begin{tikzpicture}[descr/.style={fill=white}]
	\node (Gamma) at (-1,1) {$x$};
	\node (GammaBar) at (2,0) {$\stim{\imGraph_f}$};
	\node (Y) at (0,-2) {$y$};
	\node (YBar) at (2,-2) {$q,$};
	\node (GammaY) at (0,0) {$\stim{\imGraph_f}\times_{q} y $};
	\node at (.25,-.45) {$\ulcorner$};
	\path[->,bend right,font=\scriptsize]
	(Gamma) edge node[left]{$f$} (Y);
	\path[->,font=\scriptsize]
	(GammaBar) edge node[right]{${g'}$}(YBar)
	(GammaY) edge node[right]{$\underline{{g}}'$} (Y);
	\path[right hook->,bend left,font=\scriptsize]
	(Gamma) edge node[descr]{$\!\!\!\circ\!\!\!$}node[above]{$j \circ h^\inv$} (GammaBar);
	\path[right hook->,font=\scriptsize]
	(Y) edge node[descr]{$\!\!\!\circ\!\!\!$}node[below]{$j_{\!_l}$} (YBar)
	(GammaY) edge node[descr]{$\!\!\!\circ\!\!\!$}node[below]{$\underline{j_{\!_l}}$} (GammaBar);
	\path[right hook->,dashed,font=\scriptsize]
	(Gamma) edge node{${j'}$} (GammaY);
	\end{tikzpicture}
	\]
	induced by the universal property of pullbacks in $\sftSch/S$. Since $j \circ h^\inv$ is an open immersion, so is $j'$. {The morphism} ${g'}$ is proper, and so is $\underline{{g}}'$, by \cite[Prop.5.4.2]{EGA2}. Since $f$ is also proper, \cite[Cor.5.4.3(i)]{EGA2} implies that the immersion ${j}'$ is proper, and hence a closed immersion, by \cite[Cor.18.12.6]{EGA4.4}. Since $c'$ is the scheme-{theoretic} image of $(j_{\!_i}\times j_{\!_l})\circ c\circ h^\inv$ and $j'$ is a closed immersion, $j'$ is also surjective. Thus, $j'$ is an isomorphism for being a surjective open immersion. Therefore, ${g'}$ defines a morphism of compactifications ${g'}:i'\to l$, as in Remark \ref{Rem:Cartesian:Square:Morphism:Compactifications}, because the underlying scheme of $z'$ is reduced. Moreover, $f$ is a base change in $\sftSch/S$ of ${g'}$ along $j_{\!_l}$.

	Similarly, one sees that there exists a morphism of compactifications ${h'}:i'\to i$ such that $\id_x$ is a base change in $\sftSch/S$ of ${h'}$ along $j_{\!_i}$.
\end{proof}

\Prop{}{Let $x$ and $y$ be $S$-schemes, let $i:z\cim p$ and $l:w\cim q$ be compactifications of $x$ and $y$, respectively, and suppose that $(f_0,g_0),(f_1,g_1):i\to l$ are parallel morphisms of compactifications. Then, there exist an $S$-scheme $x'$, a compactification $i':z'\cim p'$ of $x'$, and a morphism of compactifications $(f,g):i'\to i$ for which $	(f_0,g_0) \circ (f,g)=(f_1,g_1) \circ (f,g)$. Moreover, when ${g_0}_{|_x}={g_1}_{|_x}$, {the $S$-scheme} $x'$ can be chosen to be $x$, and the morphism $g$ can be chosen such that $\id_x$ is a base change in $\sftSch/S$ of $g$ along $j_{\!_i}$.
}{Compactification:Refinement}
\begin{proof} Let $j_{\!_i}:x\oim p$ and $j_{\!_l}:y\oim q$ be complementary open immersions of $i$ and $l$, respectively, and let ${g_k}_{|_x}:x\to y$ be a base change in $\sftSch/S$ of $g_k$ along $j_{\!_l}$, for $k=0,1$. Consider the solid diagram
	\begin{equation}\label{Compactification:Diagram:3}
	\begin{tikzpicture}[descr/.style={fill=white}]
	\node (XBar0) at (-1.5,1.5) {${p}$};
	\node (Gamma1) at (3,0) {${\imGraph_{g_1}}$};
	\node (Gamma0) at (0,-3) {${\imGraph_{g_0}}$};
	\node (XProd) at (3,-3) {$p\times q$};
	\node (Gamma) at (0,0) {$p'$};
	\node at (.25,-.4) {$\ulcorner$};
	\node (YBar0) at (4.5,-4.5) {$q$};
	\node (X0) at (-6.5,-0.5) {${x}$};
	\node (G1) at (-2,-2) {${\imGraph_{{g_1}_{|_x}}}$};
	\node (G0) at (-5,-5) {${\imGraph_{{g_0}_{|_x}}}$};
	\node (XP) at (-2,-5) {${x}\times {y}$};
	\node (G) at (-5,-2) {${x}'$};
	\node at (-4.75,-2.4) {$\ulcorner$};
	\node (Y0) at (-0.5,-6.5) {${y}$};
	\path[right hook->,dotted,font=\scriptsize]
	(Gamma) edge node{$\diagup$}node[right]{$\overline{i}'_1$} node[descr,pos=.39,sloped] {} (Gamma0)
	(Gamma) edge node{$\diagup$}node[above]{$\overline{i}'_0$} (Gamma1);
	\path[->,bend right,font=\scriptsize]
	(XBar0) edge node[below]{$\overline{h}_0$} node[descr,pos=.43,sloped] {} node[descr,pos=.7,sloped] {}(Gamma0)
	(Gamma0) edge node[descr,pos=.175,sloped] {} (YBar0)
	(G0) edge (Y0)
	;
	\path[right hook->,dashed,font=\scriptsize]
	(G0) edge node[descr]{$\!\!\!\!\circ\!\!\!\!$}node[below]{$j_0$}node[descr,pos=.61,sloped] {}node[descr,pos=.89,sloped] {}(Gamma0)
	(G1) edge node[descr]{$\!\!\!\!\circ\!\!\!\!$}node[right]{$j_1$}(Gamma1);
	
	\path[right hook->, dotted,font=\scriptsize]
	(G) edge node[descr]{$\!\!\!\!\circ\!\!\!\!$}node[above]{$j'$}node[descr,pos=.415,sloped] {}(Gamma);	
	\path[right hook->,font=\scriptsize]
	(XP) edge node[descr]{$\!\!\!\!\circ\!\!\!\!$}node[below]{$j_{\!_i}\times j_{\!_l}$}   node[descr,pos=.27,sloped] {} (XProd)
	(X0) edge node[descr]{$\!\!\!\!\circ\!\!\!\!$}node[above]{$j_{\!_i}$}(XBar0)
	(Y0) edge node[descr]{$\!\!\!\!\circ\!\!\!\!$}node[below]{$j_{\!_l}$}(YBar0)
	;

	\path[->,bend left,font=\scriptsize]
	(XBar0) edge node[above]{$\overline{h}_1$} (Gamma1)
	(Gamma1) edge (YBar0)
	(G1) edge (Y0);

	\path[->,font=\scriptsize]
	(XProd) edge node[above]{$\pi_{q}$} (YBar0);
	
	\path[right hook->,font=\scriptsize]	
	(Gamma1) edge node{$\diagup$}node[right]{$\overline{i}_1$}(XProd)
	(Gamma0) edge node{$\diagup$}node[above]{$\overline{i}_0$} (XProd);

	\path[->,font=\scriptsize]
	(XP) edge node[below]{$\pi_{y}$} (Y0);
	
	\path[->,bend right,font=\scriptsize]
	(X0) edge node[left]{${h}_0$} (G0);
	
	\path[->,bend left,font=\scriptsize]
	(X0) edge node[above]{${h}_1$} (G1);
	
	\path[right hook->,font=\scriptsize]
	(G1) edge node{$\diagup$}node[left]{$i_1$}(XP)
	(G0) edge node{$\diagup$}node[below]{$i_0$} (XP);
	
	\path[right hook->,dotted,font=\scriptsize]
	(G) edge node{$\diagup$}node[left]{$i'_1$} (G0)
	(G) edge node{$\diagup$}node[below]{$i'_0$} (G1);
	\end{tikzpicture}
	\end{equation}
	of $S$-schemes, which is induced by the definition of the graphs $\imGraph_{{g_k}_{|_x}}$ and $\imGraph_{g_k}$ of ${g_k}_{|_x}$ and ${g_k}$, respectively, for $k=0,1$. In the solid diagram \eqref{Compactification:Diagram:3}, the side subdiagrams are commutative, but the front and back faces are not necessarily commutative. The morphisms $h_k$ and $\overline{h}_k$ are the unique morphisms that factorise $(\id_{x},{g_k}_{|_x})$ and $(\id_p,g_k)$ in $\sftSch/S$ as $ {i}_k \circ {h}_k$ and $ \overline{i}_k \circ \overline{h}_k$, respectively, for $k=0,1$. Whereas, the morphisms {$\pi_x$, $\pi_y$, $\pi_p$ and $\pi_{q}$ are the Cartesian products projections}.

	\medskip
	
	The proof is based on basic constructions on {this} solid diagram, and follows through commutative subdiagrams chase. 
	
	\medskip

	The morphisms $h_k$ and $\overline{h}_k$ are isomorphisms with inverses $\pi_x \circ i_k$ and $\pi_p \circ \overline{i}_k$, respectively, for $k=0,1$. In particular, {we} have an open immersion $j_k:\imGraph_{{g_k}_{|_x}}\oim \imGraph_{g_k}$, given by $j_k=\overline{h}_k \circ j_{\!_i} \circ h_k^\inv $, for $k=0,1$. Then, the horizontal square containing $j_0$ and the vertical square containing $j_1$ are commutative, \ie $\overline{i}_k \circ j_k=(j_{\!_i}\times j_{\!_l}) \circ i_k$, for $k=0,1$.
	
	\medskip
	
	Let ${x}'$ be the fibre product $\imGraph_{{g_0}_{|_x}}\times_{x\times y}\imGraph_{{g_1}_{|_x}}$, with the fibre product projections $i'_0$ and $i'_1$, and let $p'$ be the fibre product $\imGraph_{g_0}\times_{p\times q}\imGraph_{g_1} ${, with the fibre product projections $\overline{i}'_0$ and $\overline{i}'_1$}. Then, there exists a unique morphism $x'\to p'$ of $S$-schemes, induced by the universal property of fibre products, making the squares containing it commute, which we denote by $j'$.
	
	\medskip
	
	In fact, a diagram chase shows the two squares that contain the morphism $j'$ to be Cartesian in $\sftSch/S$, as $j_{\!_i}\times j_{\!_l}$, $\overline{i}'_{k}$ and $i'_{k}$ are monomorphisms of $S$-schemes, for $k=0,1$. Then, in particular, $j'$ is an open immersion, and there exists a compactification $i':z'\cim p'$ of $x'$, where $i'$ is a complementary closed immersion to $j'$, endued with the reduced induced structure. 
	
	\medskip
	
	To establish the desired morphism of compactifications, notice that 
	\[
	h_1^\inv \circ {i}'_{0}=(\pi_x \circ i_1 \circ h_1) \circ (h_1^\inv \circ {i}'_{0})=\pi_x \circ i_1 \circ {i}'_{0}
	=\pi_x \circ i_0 \circ {i}'_{1}=(\pi_x \circ i_0 \circ h_0) \circ (h_0^\inv \circ {i}'_{1})=h_0^\inv \circ {i}'_{1},
	\]
	and similarly, $\overline{h}_1^\inv \circ \overline{i}'_{0}=\overline{h}_1^\inv \circ \overline{i}'_{1}$. Let $g\coloneqq \overline{h}_1^\inv \circ \overline{i}'_{0}=\overline{h}_0^\inv \circ \overline{i}'_{1}$, then the morphism ${h}_1^\inv \circ {i}'_{0}={h}_0^\inv\circ  {i}'_{1}$ is a base change in $\sftSch/S$ of $g$ along $j_{\!_i}$, which we denote by $g_{|_{x'}}$. Since the square
	\[
	\begin{tikzpicture}[descr/.style={fill=white}]
	\node (X0) at (0,0) {$x'$};
	\node (XBar0) at (2,0) {$p'$};
	\node (X1) at (0,-2) {$x$};
	\node (XBar1) at (2,-2) {$p$};
	\node at (.25,-.4) {$\ulcorner$};
	\path[right hook->,font=\scriptsize]
	(XBar0) edge node{$\diagup$}node[right]{${g}$} (XBar1)
	(X0) edge node{$\diagup$}node[left]{$g_{|_{x'}}$} (X1);
	\path[right hook->,font=\scriptsize]
	(X0) edge node[descr]{$\!\!\!\circ\!\!\!$}node[above]{$j'$} (XBar0)
	(X1) edge node[descr]{$\!\!\!\circ\!\!\!$}node[below]{$j_{\!_i}$} (XBar1);
	\end{tikzpicture}
	\]
	is Cartesian in $\sftSch/S$ and the underlying scheme of $z'$ is reduced, there exists a morphism of compactifications $(f,g):i'\to i$, for the unique morphism $f:z'\to z$ that factorises $g\circ i'$ in $\sftSch/S$ as $i \circ f$, see Remark \ref{Rem:Cartesian:Square:Morphism:Compactifications}. Then, one has
	\[
	g_0 \circ g=(\pi_{q} \circ \overline{i}_0 \circ \overline{h}_0) \circ (\overline{h}_0^\inv \circ \overline{i}'_1)=
	\pi_{q} \circ \overline{i}_0 \circ \overline{i}'_1=
	\pi_{q} \circ \overline{i}_1 \circ\overline{i}'_0=
	( \pi_{q} \circ \overline{i}_1 \circ \overline{h}_1) \circ (\overline{h}_1^\inv \circ \overline{i}'_0)=
	g_1 \circ g,
	\]
	and hence $l \circ f_0 \circ f=g_0 \circ g \circ i'=g_1 \circ g \circ i'=l \circ f_1 \circ f$. Since $l$ is a monomorphism of $S$-schemes, one has $(f_0,g_0)\circ (f,g)=(f_1,g_1) \circ (f,g)$.

	\medskip
	
	Moreover, when ${g_0}_{|_x}={g_1}_{|_x}$, the universal property of pullbacks implies the existence of a morphism $x\to x'$ in $\sftSch/S$ that factorises the isomorphism $h_0=h_1$. Since $i'_0$ is a {closed} immersion, such a morphism $x \to x'$ is an isomorphism. Pullbacks are determined up to isomorphisms; thus, we may choose $x'$ to be $x$, in which case $\id_x$ is a base change in $\sftSch/S$ of $g$ along $j_{\!_i}$.
\end{proof}

\Cor{}{Let $x$ be an $S$-scheme. Then, the category $\Comp_{\!_S}(x)$ is cofiltered.
}{Compactification:Cofiltered} 
\begin{proof}
	Since $S$ is a Noetherian scheme\footnote{In the light of \cite{Con:07:DNNC:Published}, one may generalise most statements in this section to a quasi-compact quasi-separated base scheme $S$.}, \nagata's Compactification Theorem implies that $\Comp_{\!_S}(x)$ is nonempty, as seen in Remark \ref{Rem:Nagata}. Then, the {statement} of the corollary is a direct result of {Proposition} \ref{Prop:Compactification:Proper:Morphism}, for $f=\id_x$, and {Proposition} \ref{Prop:Compactification:Refinement}.
\end{proof}

Let $f:x\to y$ be a morphism of $S$-schemes, and let $l:w\cim q$ be a compactification of $y$. We denote {by $\Comp_{\!_S}(f,l)$} the full subcategory in $\Comp_{\!_S}(x)$ that satisfies the property
\begin{itemize}
	\item[] a compactification $i$ of $x$ belongs to $\Comp_{\!_S}(f,l)$ \iff it admits a morphism of compactifications $g:i\to l$ such that $f$ is a base change in $\sftSch/S$ of $g$ along $j_{\!_l}$.
\end{itemize}
Also, let $\Comp_{\!_S}(f)$ denote the full subcategory in $\Comp_{\!_S}(x)$ of compactifications of $x$ that belong to $\Comp_{\!_S}(f,l)$ for some compactification $l$ of $y$.

\Cor{}{Assume that $f:x\to y$ is a proper morphism of $S$-schemes, and let $l:w\cim q$ be a compactification of $y$. Then, the category $\Comp_{\!_S}(f,l)$ is co-cofinal in $\Comp_{\!_S}(x)$, and so is $\Comp_{\!_S}(f)$. Moreover, the categories $\Comp_{\!_S}(f,l)$ and $\Comp_{\!_S}(f)$ are cofiltered.
}{Compactification:Cofiltered:fj}
\begin{proof}
	A direct consequence of {Proposition} \ref{Prop:Compactification:Proper:Morphism}, {Proposition} \ref{Prop:Compactification:Refinement}, and \cite[Ch.0.\S.3.2-3]{Tam:94:IEC}.
\end{proof}

\subsection{Extensions of \texorpdfstring{\cdp-}{cdp-}Functors}\label{Extension:PS}
Throughout this subsection, let $F:\PropSch/S\to \scC$ be a functor to a {\waldhausen category} that satisfies the properties \ref{itm:PS:Closed}-\ref{itm:PS:Proper}. The functor $F$ extends to a functor $F^\ps:{\sftSch/S}^\prop_\open\to \scC$ that satisfies the excision property, as in Proposition \ref{Prop:Excision}. Moreover, for a symmetric monoidal \waldhausen category $\scC$, if $F$ is weak monoidal, then so is $F^\ps$, as in Proposition \ref{Prop:Monoidal}. The main statement here is Theorem \ref{Th:Induced:Measure}.

\Def{}{A functor $\PropSch/S\to \scC$ to a {\waldhausen category} that satisfies the properties \ref{itm:PS:Closed}-\ref{itm:PS:Proper} is called a \indx{\cdp-functor}.
}{}
This terminology is motivated by Definition \ref{Def:cd:topology} and Proposition \ref{Prop:cd:Cosheafification}.

\Rem{}{The properties \ref{itm:PS:Closed}-\ref{itm:PS:Proper} imply the following property.
	\begin{enumerate}[label=\textup{(\textbf{PS\arabic*})},ref=\textup{(\textbf{PS\arabic*})}]
		\setcounter{enumi}{3}
		\item\label{itm:PS:Surjective} $F$ maps every surjective closed immersion of proper $S$-schemes to an {isomorphism}. 
		\setcounter{enumi}{0}
	\end{enumerate}	
	That is, for a surjective closed immersion $i:z\cim p$ of proper $S$-schemes, the square
	\[
	\begin{tikzpicture}[descr/.style={fill=white}]
	\node (X0) at (0,0) {$\emptyset_{\!_S}$};
	\node (XBar0) at (2,0) {$\emptyset_{\!_S}$};
	\node (X1) at (0,-2) {$z$};
	\node (XBar1) at (2,-2) {$p$};
	\node at (.25,-.4) {$\ulcorner$};
	\path[right hook->,font=\scriptsize]
	(XBar0) edge node{$\diagup$} (XBar1)
	(X0) edge node{$\diagup$} (X1)
	(X0) edge node{$\diagup$} (XBar0)
	(X1) edge node{$\diagup$}node[below]{$i$} (XBar1);
	\end{tikzpicture}
	\]
	is a \cdp-square of proper $S$-schemes, which is mapped by $F$ to a square of cofibrations in $\scC$, by \ref{itm:PS:Closed}. Then, \ref{itm:PS:Proper} and \ref{itm:PS:Initial} imply that $F(i)$ is the composite isomorphism $F(z)\cong F(z)/F(\emptyset_{\!_S})\cong F(p)/F(\emptyset_{\!_S})\cong F(p)$.
}{Eq:Surjective}
In the sequel, we denote $F(f)$ by $f_\ast$, for every morphism $f$ of proper $S$-schemes. Also, we adopt the same notation for other (contravariant) functors, when no confusion arises.

\medskip

Assume that $x$ and $y$ are $S$-schemes, let $i:z\cim p$ and $l:w\cim q$ be compactifications of $x$ and $y$, respectively, and let $(f,g):i\to l$ be a morphism of compactifications, as in \eqref{Eq:Compactifications:Morphism}. The morphism $(f,g)$ is mapped to the solid commutative square
\[\label{Eq:F:(f,g)}
\begin{tikzpicture}[descr/.style={fill=white}]
\node (Z0) at (0,0) {$F(z)$};
\node (XBar0) at (2,0) {$F(p)$};
\node (C0) at (4,0) {$C_{\!\!_F}(i)$};
\node (Z1) at (0,-2) {$F(w)$};
\node (XBar1) at (2,-2) {$F(q)$};
\node (C1) at (4,-2) {$C_{\!\!_F}(l)$};
\path[->,font=\scriptsize]
(XBar0) edge node[right]{$g_\ast$} (XBar1)
(Z0) edge node[left]{${f}_\ast$} (Z1);
\path[>->,font=\scriptsize]
(Z0) edge node[above]{${i}_\ast$} (XBar0)
(Z1) edge node[below]{${l}_\ast$} (XBar1);
\path[->>,dashed, font=\scriptsize]
(XBar0) edge node[above]{${ \epsilon_{\!_i}}$} (C0)
(XBar1) edge node[below]{${ \epsilon_{\!_l}}$} (C1);
\path[->,dotted, font=\scriptsize]
(C0) edge node[right]{$\big({f_\ast,g_\ast}\big)$} (C1);
\end{tikzpicture}
\]
in $\scC$. Sine $F$ satisfies \ref{itm:PS:Closed}, both $i_\ast$ and $l_\ast$ are cofibrations in $\scC$. Let $C_{\!\!_F}(i)$ and $C_{\!\!_F}(l)$ be the cofibres of $i_\ast$ and $l_\ast$, respectively. Since the left solid square commutes, there exists a unique morphism $C_{\!\!_F}(i)\to C_{\!\!_F}(l)$ in $\scC$ that makes the whole diagram commute, which we denote by $({f_\ast,g_\ast})$. 

\medskip

That defines a functor $\label{Eq:Compact:Support:Cofibre}
C_{\!\!_F}:\Comp_{\!_S}\to \scC$, given on an object $i\in \Comp_{\!_S}$ by $C_{\!\!_{F}}(i)$ and on a morphism $(f,g)$ in $\Comp_{\!_S}$ by $({f_\ast,g_\ast})$.

\Rem{}{For every $S$-scheme $x$, let $C_{\!\!_{F,x}}$ be the composition of the functor $C_{\!\!_{F}}$ with the inclusion functor $\Comp_{\!_S}(x)\incl \Comp_{\!_S}$. For a morphism of compactifications $(f,g):i'\to i$ in $\Comp_{\!_S}(x)$, \ie a commutative diagram
	\begin{equation}\label{Comp:x:Comp}
	\begin{tikzpicture}[descr/.style={fill=white}]
	\node (fp) at (0,0) {$z\times_{p}p'$};
	\node (X0) at (-1,1) {$z'$};
	\node (XBar0) at (2,0) {$p'$};
	\node (X1) at (0,-2) {$z$};
	\node (XBar1) at (2,-2) {$p,$};
	\node at (.25,-.4) {$\ulcorner$};
	\path[->,font=\scriptsize]
	(XBar0) edge node[right]{$g$} (XBar1);
	
	\path[right hook->,font=\scriptsize]
	(X1) edge node{$\diagup$}node[below]{$i$} (XBar1);
	\path[right hook->,bend left,font=\scriptsize]
	(X0) edge node{$\diagup$}node[above]{$i'$} (XBar0);
	\path[->,bend right,font=\scriptsize]
	(X0) edge node[left]{$f$} (X1);
	
	\path[dashed, ->,font=\scriptsize]
	(fp) edge node[right]{$\underline{g}$} (X1);
	\path[dashed, right hook->,font=\scriptsize]
	(fp) edge node{$\diagup$} node[below]{$\underline{i}$} (XBar0);
	\path[dotted, right hook ->,font=\scriptsize]
	(X0) edge node{$\diagup$} node[right]{$c$}(fp);
	\end{tikzpicture}
	\end{equation}
	in $\sftSch/S$, in which $c$ is a surjective closed immersion and $\id_x$ is a base change in $\sftSch/S$ of $g$ along $j_{\!_i}$. The pullback square in \eqref{Comp:x:Comp} is a \cdp-square of proper $S$-schemes, as $g$ is a proper morphism, $i$ is a closed immersion, and $\id_x$ is a base change in $\sftSch/S$ of $g$ along $j_{\!_i}$. Since $F$ satisfies \ref{itm:PS:Closed} and \ref{itm:PS:Proper}, the morphism $(\underline{g}_\ast,g_\ast)$ is an isomorphism in $\scC$. Also, $c_\ast$ is an isomorphism in $\scC$, as $F$ satisfies \ref{itm:PS:Surjective}. Hence, $({f}_\ast,g_\ast)$ is an isomorphism. Therefore, $C_{\!\!_{F,x}}$ is a diagram of isomorphisms, and hence $\lim C_{\!\!_{F,x}}$ exists in $\scC$. For a compactification $i$ of $x$, we denote the limit projection $\lim C_{\!\!_{F,x}} \to C_{\!\!_F}(i)$ by $\iota_{\!_i}$.
}{Cofiltered:Extension}

Since $\lim C_{\!\!_{F,x}}$ satisfies the excision property for proper $S$-schemes, we will define the extension $F^\ps$ on an $S$-scheme $x$ by $F^\ps(x)\coloneqq \lim C_{\!\!_{F,x}}$. Then, in \S.\ref{Sec:Proper:Pushforwards} and \S.\ref{Sec:Pull:Back}, we define $F^\ps$ on proper morphisms and formal inverses of open immersions, respectively.

\subsubsection{Proper Pushforwards}\label{Sec:Proper:Pushforwards} We need to assign for every proper morphism $f:x\to y$ of $S$-schemes a unique morphism $\lim C_{\!\!_{F,x}} \to \lim C_{\!\!_{F,y}}$, that is independent of the choice of compactifications and morphisms between them. We show below that the canonical choices of such morphisms coincide{, as in Corollary \ref{Cor:Compact:Support:Functor:Proper:Morphisms}}.

\Lem{}{Assume that $f:x\to y$ is a proper morphism of $S$-schemes, and let $l:w\cim q$ be a compactification of $y$. Then, the morphism $
	g_\ast \circ \iota_{\!_i}:\lim C_{\!\!_{F,x}}\to C_{\!\!_F}(l)$ is independent of the choice of the compactification $i$ in $\Comp_{\!_S}(f,l)$ and of the morphism of compactifications $g:i\to l$ for which $f$ is a base change in $\sftSch/S$ along $j_{\!_l}$. We denote this morphism by $\varrho_l^f$.
}{Waldhausen:Compact:Support:Limit:Morphisms}
\begin{proof}
	Since $f$ is proper, the category $\Comp_{\!_S}(f,l)$ is nonempty, by Proposition \ref{Prop:Compactification:Proper:Morphism}. Suppose that $i_k:z_k\cim p_k$ is a compactification in $\Comp_{\!_S}(f,l)$, and let $g_k:i_k\to l$ be a morphism of compactifications such that {$f$ is a base change in $\sftSch/S$ of $g_k$ along $j_{\!_l}$,} for $k=0,1$. Since the category $\Comp_{\!_S}(f,l)$ is cofiltered, there exists a compactification $i:z\cim p$ in $\Comp_{\!_S}(f,l)$ and a morphism of compactifications $g'_k:i\to i_k$ such that {$\id_x$ is a base change in $\sftSch/S$ of $g'_k$ along $j_{\!_{i_k}}$,} for $k=0,1$. {Proposition} \ref{Prop:Compactification:Refinement} implies that $i$, $g'_0$, and $g'_1$ can be chosen such that $g_0 \circ g'_0=g_1 \circ g'_1$. Thus, ${g_0}_\ast \circ \iota_{\!_{i_0}}={g_0}_\ast \circ {g'_0}_\ast \circ \iota_{\!_i}={g_1}_\ast \circ {g'_1}_\ast \circ \iota_{\!_i}={g_1}_\ast \circ \iota_{\!_{i_1}}$.

	\medskip
	
	On the other hand, suppose that $i:z\cim p$ is a compactification in $ \Comp_{\!_S}(f,l)$ and let $g_0,g_1:i\to l$ be parallel morphisms of compactifications such that $f$ is a base change in $\sftSch/S$ of $g_k$ along $j_{\!_l}$, for $k=0,1$. By {Proposition} \ref{Prop:Compactification:Refinement}, there exists a refining compactification $i':z'\cim p'$ of $x$ and a morphism of compactifications $g:i'\to i$ in $\Comp_{\!_S}(f,l)$ such that 
	$g_0\circ g=g_1\circ g$. Thus, ${g_0}_\ast \circ \iota_{\!_i}={g_0}_\ast \circ g_\ast \circ \iota_{\!_{i'}}={g_1}_\ast \circ g_\ast \circ \iota_{\!_{i'}}={g_1}_\ast \circ \iota_{\!_i}$.
\end{proof}

\Lem{}{Assume that $f:x\to y$ is a proper morphism of $S$-schemes, let $l_k:w_k\cim q_k$ be a compactification of $y$, for $k=0,1$, and let $g:l_0\to l_1$ be a morphism of compactifications in $\Comp_{\!_S}(y)$. Then, $\varrho_{l_1}^f={g}_\ast \circ \varrho_{l_0}^f$.
}{}
\begin{proof}
	Since $f$ is proper, the category $\Comp_{\!_S}(f,l_k)$ is nonempty, for $k=0,1$. Let $i_k:z_k\cim p_k $ be a compactification in $\Comp_{\!_S}(f,l_k)$, and let $g_k:i_k\to l_k$ be a morphism of compactifications such that $f$ is a base change in $\sftSch/S$ of $g_k$ along $j_{\!_{l_k}}$, for $k=0,1$. Since $\Comp_{\!_S}(x)$ is cofiltered, there exists a compactification $i:z\cim p$ of $x$ and a morphism of compactifications $g'_k: i\to i_k$ in $\Comp_{\!_S}(x)$, for $k=0,1$, by {Corollary} \ref{Cor:Compactification:Cofiltered}. Then, {Lemma} \ref{Lem:Waldhausen:Compact:Support:Limit:Morphisms} implies that
	\[
	\varrho_{l_1}^f=
	{g_1}_\ast\circ \iota_{\!_{i_1}}=
	{g_1}_\ast \circ {g'_1}_\ast \circ \iota_{\!_i}=
	{g}_\ast \circ {g_0}_\ast \circ {g'_0}_\ast \circ \iota_{\!_i}=		
	{g}_\ast \circ{g_0}_\ast \circ \iota_{\!_{i_0}}=
	{g}_\ast \circ \varrho_{l_0}^f,
	\]
	as $g\circ g_0\circ g'_0 $ and $ g_1\circ g'_1$ are parallel morphisms of compactifications and $f$ is a base change in $\sftSch/S$ of both $g\circ g_0\circ g'_0$ and $g_1\circ g'_1$ along $j_{\!_l}$, 
\end{proof}

\Cor{}{Assume that $f:x\to y$ is a proper morphism of $S$-schemes. Then, there exists a unique morphism $f_\shriek:\lim C_{\!\!_{F,x}}\to \lim C_{\!\!_{F,y}}$ in $\scC$ for which
	\begin{equation}\label{PS:Proper:Pushforward}
	\iota_{\!_l} \circ f_\shriek=\varrho_l^f ,
	\end{equation}
	for every compactification $l$ of $y$.
}{Compact:Support:Functor:Proper:Morphisms}
The uniqueness of the morphism $f_\shriek$, for a proper morphism $f$ of $S$-schemes, implies the functoriality of  pushforwards along proper morphisms. That is, for proper morphisms $f:x\to y$ and $g:y\to z$ of $S$-schemes, one has
$(g \circ f)_\shriek=g_\shriek \circ f_\shriek$ and $(\id_x)_\shriek=\id_{\!_{\lim C_{\!\!_{F,x}}}}$.
\Def{}{The \indx{properly supported counterpart} of $F$ is a functor $F^\ps:{\sftSch/S}^\prop \to \scC$	that sends an $S$-schemes $x$ to $\lim C_{\!\!_{F,x}}$, as in Remark \ref{Rem:Cofiltered:Extension}, and sends a proper morphism $f:x\to y$ of $S$-schemes to the unique morphism $f_\shriek:\lim C_{\!\!_{F,x}}\to \lim C_{\!\!_{F,y}}$ asserted by {Corollary} \ref{Cor:Compact:Support:Functor:Proper:Morphisms}.
}{Compact:Support:Functor:Proper}
When $\scC$ has cofiltered limits, the functor $F^\ps$ may be defined on proper morphisms similarly, even when $F$ is not a \cdp-functor. However, such a functor does not necessarily satisfy the excision property.

\Eg{}{Assume that $p$ is a proper $S$-scheme. Since the category $\Comp_{\!_S}(p)$ admits an initial object, namely $\emptyset_p:\emptyset_{\!_S}\cim p$, and $F$ satisfies \ref{itm:PS:Initial}, one has $F^\ps(p)\cong C_{\!\!_F}(\emptyset_p)\cong F(p)$.
}{}

The functor $F^\ps$ satisfies generalisations of the properties \ref{itm:PS:Closed}-\ref{itm:PS:Surjective} to the category ${\sftSch/S}^\prop$, as in the following proposition.
\Prop{}{The functor $F^\ps:{\sftSch/S}^\prop\to \scC$
	\begin{enumerate}	[label=\textup{(\textbf{PS\arabic*$'$})},ref=\textup{(\textbf{PS\arabic*$'$})}]
		\item\label{itm:PS:C} maps closed immersions of $S$-schemes to cofibrations in $\scC$;
		\item\label{itm:PS:P} maps \cdp-squares of $S$-schemes to {pushout squares} in $\scC$; 
		\item\label{itm:PS:I} maps the empty $S$-scheme to the zero object in $\scC$; and 
		\item\label{itm:PS:D} maps surjective closed immersions of $S$-schemes to isomorphisms.
	\end{enumerate}
}{PS:Properties}
\begin{proof}
	The statement \ref{itm:PS:I} is evident; whereas \ref{itm:PS:D} follows from the other statements, as seen in Remark \ref{Rem:Eq:Surjective}.
	\begin{itemize}
		\item[\ref{itm:PS:C}] Assume that $i:v\cim x$ is a closed immersion of $S$-schemes. Let $l:w\cim q$ be a compactification of $x$ with complementary open immersion $j_{\!_l}:x\oim q$, and let $i':p\cim q$ be the scheme-theoretic image of the immersion $j_{\!_l}\circ i$. Then, Lemma \ref{Lem:Open:Closed:Closed:Open} implies the existence of an open immersion $j'_{\!_l}:v\oim p$ for which the solid square
		\begin{equation}\label{Closed:Cofibration}
		\begin{tikzpicture}[descr/.style={fill=white}]
		\node (X0) at (0,0) {$v$};
		\node (W0) at (4,0) {$z$};
		\node (XBar0) at (2,0) {$p$};
		\node (X1) at (0,-2) {$x$};
		\node (W1) at (4,-2) {$w$};
		\node (XBar1) at (2,-2) {$q$};
		\node at (.25,-.4) {$\ulcorner$};
		\node at (3.75,-.4) {$\urcorner$};
		\path[right hook->,font=\scriptsize]
		(XBar0) edge node{$\diagup$} node[right]{$i'$} (XBar1)
		(X0) edge node{$\diagup$}node[left]{$i$} (X1)
		(X0) edge node[descr]{$\!\!\!\circ\!\!\!$}node[above]{$j'_{\!_l}$} (XBar0)
		(X1) edge node[descr]{$\!\!\!\circ\!\!\!$}node[below]{$j_{\!_l}$} (XBar1)
		;
		\path[left hook->,font=\scriptsize]
		(W1) edge node{$\diagup$}node[below]{$l$} (XBar1)
		;
		\path[right hook->,dashed,font=\scriptsize]
		(W0) edge node{$\diagup$}node[right]{$\underline{i}$} (W1)
		;
		\path[left hook->,dashed,font=\scriptsize]
		(W0) edge node{$\diagup$}node[above]{$\underline{l}$} (XBar0)
		;
		\end{tikzpicture}
		\end{equation}
		is Cartesian in $\sftSch/S$. Let $\underline{l}$ be a base change in $\sftSch/S$ of $l$ along $i'$. {Since open complements are closed under pullbacks} and $j_{\!_l}$ is a complementary open immersion to $l$, one finds that $j'_{\!_l}$ is a complementary open immersion to $\underline{l}$. Hence, $\underline{l}:z\cim p$ is a compactification in $\Comp_{\!_S}(i,l)$ and $(\underline{i},i'):\underline{l}\to l$ is a morphism of compactifications such that $i$ is a base change in $\sftSch/S$ of $i'$ along $j_{\!_l}$.

		\medskip
		
		There exists a pushout in {$\PropSch/S$} of the span of closed immersions $\underline{l}$ and $\underline{i}$, which we denote by $\underline{q}$. In fact, since the right square in \eqref{Closed:Cofibration} is Cartesian, {there exists a bicartesian square} 
\begin{equation}\label{Closed:Cofibration:2}
\begin{tikzpicture}[descr/.style={fill=white}]
\node (m-1-1) at (0,0) {$z$};
\node (m-1-2) at (2,0) {$p$};
\node (m-2-1) at (0,-2) {$w$};
\node (m-2-2) at (2,-2) {$\underline{q}$};
\node at (.25,-.4) {$\ulcorner$};
\node at (1.75,-1.6) {$\lrcorner$};
\path[right hook->,font=\scriptsize]
(m-1-1) edge node{$\diagup$}node	[left]	{$\underline{i}$} (m-2-1)
(m-1-2) edge node{$\diagup$}node	[right]	{$\underline{i}'$} (m-2-2)
(m-2-1) edge node{$\diagup$}node	[below]	{$\underline{l}'$} (m-2-2)
(m-1-1) edge node{$\diagup$}node	[above]	{$\underline{l}$} (m-1-2);	
\end{tikzpicture}
\end{equation}
		of closed immersions of proper $S$-schemes, see \cite[Th.3.11]{Sch:05:GSSWCP} and \cite[Tag 0B7M]{stacks-project}. In particular, the square \eqref{Closed:Cofibration:2} is a \cdp-square of proper $S$-schemes, and hence it is mapped by $F$ to a pushout square of cofibrations in $\scC$. Moreover, the unique morphism $k:\underline{q}\to q$, for which $k\circ \underline{i}'=i'$ and $k\circ \underline{l}'=l$, is a closed immersion. Consider the solid diagram \begin{equation}\label{Closed:Cofibration:PushOut}
		\begin{tikzpicture}[descr/.style={fill=white},>=angle 90,scale=1.5,text height=1.5ex, text depth=.25ex,row sep=4em, column sep=4em]
		\node (F11) at (0,0) {$F(z)$};
		\node (F12) at (2,0) {$F(p)$};
		\node (F13) at (4,0) {$C_{\!_F}(\underline{l})$};
		\node (F21) at (0,-2) {$F(w)$};
		\node (F22) at (2,-2) {$F(q)$};
		\node (F23) at (4,-2) {$C_{\!_F}(l)$};
		\node (F) at (1.3,-1.3) {$F(\underline{q})$};
		\node at (1,-1) {$\lrcorner$};
		\draw[>->,font=\scriptsize]
		(F11) edge node[above] {${\underline{l}}_\ast$}(F12)
		(F21) edge node[below] {$l_\ast$} (F22)
		(F12) edge node[left] {${\underline{i}'}_\ast$}(F)
		(F21) edge node[above] {${\underline{l}'}_\ast$}(F)
		(F11) edge node[left] {${\underline{i}}_\ast$} (F21)
		(F12) edge node[right] {$i'_\ast$}(F22)
		(F) edge node[left] {$k_\ast$}(F22)
		;
		\draw[->,dotted,font=\scriptsize]	
		(F13) edge node[right] {$(\underline{i}_\ast,i'_\ast)$}(F23)
		;
		\draw[->>,dashed,font=\scriptsize]
		(F12) edge node[above] {$\epsilon_{\!_{\underline{l}}}$} (F13)
		(F22) edge node[above]{$\epsilon_{\!_{l}}$} (F23)
		;
		\end{tikzpicture}
		\end{equation}		
		of cofibrations in $\scC$. Since $\epsilon_{\!_{\underline{l}}} \circ \underline{l}_\ast=0$, there exists a unique morphism $\gamma:F(\underline{q})\to C_{\!_F}(\underline{l})$ in $\scC$, for which $\gamma\circ \underline{l}'_\ast=0$ and $\gamma\circ \underline{i}'_\ast=\epsilon_{\!_{\underline{l}}}$. Since $\epsilon_{\!_l}$ is an epimorphism in $\scC$, a diagram chase shows that $(\underline{i}_\ast,i'_\ast)$ is a cobase change in $\scC$ of $k_\ast$ along $\gamma$. Recall that $\iota_{\!_l}\circ i_\shriek=(\underline{i}_\ast,i'_\ast)\circ \iota_{\!_{\underline{l}}}$ and that $\iota_{\!_{\underline{l}}}$ is an isomorphisms in $\scC$. Thus, $i_\shriek$ is a cobase change in $\scC$ of $k_\ast$ along $\iota_{\!_{\underline{l}}}^\inv \circ \gamma$. Since $k_\ast$ is a cofibration in $\scC$, as $k$ is a closed immersion of proper $S$-schemes, the morphism $i_\shriek$ is a cofibration in $\scC$.

		\item[\ref{itm:PS:P}] A \cdp-square of $S$-schemes defines a Cartesian cube, in which the ambient proper $S$-schemes fit into a \cdp-square. Then, a diagram chase on the cube implies the statement.
	\end{itemize}
\end{proof}

\Eg{}{Assume that $j_x:x\incl x\bigsqcup y$ is a closed open immersion of $S$-schemes with complementary closed open immersion $j_y:y\incl x\bigsqcup y$. The square 
	\[
	\begin{tikzpicture}[descr/.style={fill=white}]
	\node (X0) at (0,0) {$\emptyset_{\!_S}$};
	\node (XBar0) at (2,0) {$y$};
	\node (X1) at (0,-2) {$x$};
	\node (XBar1) at (2,-2) {$x\bigsqcup y$};
	\node at (.25,-.4) {$\ulcorner$};
	\node at (1.75,-1.6) {$\lrcorner$};
	\path[->,font=\scriptsize]
	(X0) edge node[left]{$\emptyset_x$} (X1)
	(X0) edge node[above]{$\emptyset_y$} (XBar0);
	\path[right hook->,font=\scriptsize]
	(XBar0) edge node[right]{$j_y$} (XBar1)
	(X1) edge node[below]{$j_x$} (XBar1);
	\end{tikzpicture}
	\]
	is a \cdp-square of $S$-schemes. Then, Proposition \ref{Prop:PS:Properties} implies that $F^\ps(x\bigsqcup y)\cong F^\ps(x)\coprod F^\ps(y)$.
}{Coproduct:Preservation}

\subsubsection{A Comparison Morphism}\label{Comparison:Morphism}
Let $G:\sftSch/S\to \scC$ be a functor to a {\waldhausen} category, whose restriction to $\PropSch/S$ is a \cdp-functor. We abuse notation, and use $G$ to also denote its restrictions to ${\sftSch/S}^\prop$ and $\PropSch/S$. We see below that there exists a canonical natural transformation from $G$ to its properly supported counterpart. This natural transformation is particularly useful to define the monoidal coherence morphisms in \S.\ref{MSPS:Weak:Monodial}.

\Lem{}{Let $G:\sftSch/S\to \scC$ be a functor to a {\waldhausen} category, whose restriction to $\PropSch/S$ is a \cdp-functor. Then, there exists a unique natural transformation $\varphi:G\Rightarrow G^\ps$ for which $\iota_{\!_i}\circ \varphi_x=\epsilon_{\!_i}\circ{j_{\!_i}}_\ast$,	for every $S$-scheme $x$ and for every compactification $i$ of $x$.
}{}
\begin{proof}
	Suppose that $x$ is an $S$-scheme, let $i_k:z_k\cim p_k$ be a compactification of $x$ with complementary open immersion $j_{\!_{i_k}}:x\oim p_k$, for $k=0,1$, and let $g:i_0\to i_1$ be a morphism of compactifications in $\Comp_{\!_S}(x)$. Then, $g_\ast \circ \epsilon_{\!_{i_0}}\circ {j_{\!_{i_0}}}_\ast=
	\epsilon_{\!_{i_1}}\circ g_\ast \circ {j_{\!_{i_0}}}_\ast=
	\epsilon_{\!_{i_1}}\circ {(g\circ {j_{\!_{i_0}}})}_\ast=
	\epsilon_{\!_{i_1}}\circ {j_{\!_{i_1}}}_\ast$. Therefore, by the universal property of limits, there exists a unique morphism $\varphi_x:G(x)\to G^\ps(x)$ in $\scC$ for which 
	\begin{equation}\label{Eq:Comparison:Morphism}
	\iota_{\!_i}\circ \varphi_x=\epsilon_{\!_i}\circ{j_{\!_i}}_\ast,
	\end{equation}
	for every compactification $i$ of $x$.

	\medskip
	
	Suppose that $f:x\to y$ is a proper morphism of $S$-schemes and let $l:w\cim q$ be a compactification of $y$, then the category $\Comp_{\!_S}(f,l)$ is nonempty, by {Proposition} \ref{Prop:Compactification:Proper:Morphism}. Assume that $i:z\cim p$ is a compactification of $x$ in $\Comp_{\!_S}(f,l)$, let $g:i\to l$ be a morphism of compactifications such that $f$ is a base change in $\sftSch/S$ of $g$ along $j_{\!_l}$, and consider the diagram
	\begin{equation}\label{Eq:Comparison:Morphism:2}
	\begin{tikzpicture}[descr/.style={fill=white}]
	\node (X) at (0,0) {$G(x)$};
	\node (XBar) at (2,0) {$G^\ps(x)$};
	\node (Ci) at (4,0) {$C_{\!\!_G}(i)$};
	\node (Y) at (0,-2) {$G(y)$};
	\node (YBar) at (2,-2) {$G^\ps(y)$};
	\node (Cj) at (4,-2) {$C_{\!\!_G}(l)$};
	\path[->,font=\scriptsize]
	(X) edge node[above]{$\varphi_x$} (XBar)
	(Y) edge node[below]{$\varphi_y$} (YBar)
	(X) edge node[left]{${f}_\ast$} (Y)
	(XBar) edge node[right]{$f_\shriek$} (YBar);
	\path[->, font=\scriptsize]
	(XBar) edge node[above]{${ \iota_{\!_i}}$} (Ci);
	\path[->, font=\scriptsize]
	(YBar) edge node[below]{${ \iota_{\!_l}}$} (Cj);
	\path[->, font=\scriptsize]
	(Ci) edge node[right]{${ g}_\ast$} (Cj);
	\end{tikzpicture}
	\end{equation}	
	in $\scC$. The right square in \eqref{Eq:Comparison:Morphism:2} is commutative due to {Corollary} \ref{Cor:Compact:Support:Functor:Proper:Morphisms}. Then, one has
	\begin{align*}
		\iota_{\!_l}\circ \varphi_y\circ f_\ast &=
		\epsilon_{\!_l}\circ {j_{\!_l}}_\ast \circ f_\ast=
		\epsilon_{\!_l}\circ{(j_{\!_l}\circ f)}_\ast=
		\epsilon_{\!_l}\circ{(g\circ j_{\!_i})}_\ast=
		\epsilon_{\!_l}\circ g_\ast\circ {j_{\!_i}}_\ast \\&=
		g_\ast	\circ \epsilon_{\!_i}\circ {j_{\!_i}}_\ast=
		g_\ast\circ \iota_{\!_i}\circ \varphi_x=\iota_{\!_l}\circ f_\shriek \circ \varphi_x.
	\end{align*}
	By the universal property of limits, one has $\varphi_y\circ f_\ast=f_\shriek \circ \varphi_x$. Therefore, there exists a natural transformation $\varphi:G\Rightarrow G^\ps$, whose component at an $S$-scheme $x$ is given by the unique morphism $\varphi_x$ in $\scC$ that satisfies \eqref{Eq:Comparison:Morphism}.
	
	\medskip

	Assume that $\varphi':G\Rightarrow G^\ps$ is a natural transformation for which $\iota_{\!_i}\circ \varphi'_x=\epsilon_{\!_i}\circ {j_{\!_i}}_\ast=\iota_{\!_i}\circ \varphi_x$, for every $S$-scheme $x$ and for every compactification $i:z\cim p$ of $x$. Then, the universal property of limits implies that $\varphi'_x=\varphi_x$, and hence $\varphi'=\varphi$.
\end{proof}

\Cor{}{There exists a unique natural isomorphism $\varphi:F \stackrel{\sim}{\Rightarrow} F^\ps_{|_{\PropSch/S}}:{\PropSch/S}\to \scC$, such that $\iota_{\!_i}\circ \varphi_p=\epsilon_{\!_i}\circ {j_{\!_i}}_\ast$, for every proper $S$-scheme $p$ and for every compactification $i$ of $p$.}{Comparison:Morphism}

\subsubsection{Open Pullbacks}\label{Sec:Pull:Back}The functor $F^\ps$ in Definition \ref{Def:Compact:Support:Functor:Proper} can be extended to the category ${\sftSch/S}^\prop_\open$ whose objects are $S$-schemes and whose morphisms are finite compositions of proper morphisms and formal inverses of open immersions of $S$-schemes. We define below pullbacks along open immersions, and we show in \S.\ref{Sec:Base:Change} the comparability between proper pushforwards and open pullbacks.

\medskip

Suppose that $f:x\oim y$ is an open {immersion} of $S$-schemes, let $l:w\cim q$ be a compactification of $y$ with complementary open immersion $j_{\!_l}:y\oim q$, and let $l_{\!_\red}$ be the composition of $l$ with the surjective closed immersion $w_{\!\!_\red} \cim w$. 

\medskip

The open immersion $j_{\!_l}\circ f:x\oim q$ defines a compactification $i_{l}^{f}:z\cim q$, that is a complementary closed immersion to $j_{\!_l}\circ f$, endued with the reduced induced structure. Then, $l_{\!_\red}$ factorises uniquely in $\sftSch/S$ as $i_l^f\circ c$ for a closed immersion ${c}:w_{\!_\red}\cim z$. The resulting commutative square 
\begin{equation}\label{Pull:Back:CLosed:Immersions}
\begin{tikzpicture}[descr/.style={fill=white}]
\node (X0) at (0,0) {$z$};
\node (XBar0) at (2,0) {$q$};
\node (X1) at (0,-2) {$w_{\!_\red}$};
\node (XBar1) at (2,-2) {$q$};
\path[->,font=\scriptsize]
(XBar1) edge node[right]{$\id_{q}$} (XBar0);
\path[right hook->,font=\scriptsize]
(X1) edge node{$\diagup$}node[left]{$c$} (X0)
(X0) edge node{$\diagup$}node[above]{$i_l^f$} (XBar0)
(X1) edge node{$\diagup$}node[below]{$l_{\!_\red}$} (XBar1);
\end{tikzpicture}
\end{equation}
of closed immersions of proper $S$-schemes induces the solid commutative square
\begin{equation}
\begin{tikzpicture}[descr/.style={fill=white}]
\node (Z0) at (0,0) {$F(z)$};
\node (XBar0) at (2,0) {$F(q)$};
\node (C0) at (4,0) {$C_{\!\!_F}(i_l^f)$};
\node (Z1) at (0,-2) {$F(w_{\!_\red})$};
\node (XBar1) at (2,-2) {$F(q)$};
\node (C1) at (4,-2) {$C_{\!\!_F}(l_{\!_\red}).$};

\path[->,font=\scriptsize]
(XBar1) edge node[right]{$\id_{{F(q)}}$} (XBar0);

\path[>->,font=\scriptsize]
(Z1) edge node[left]{${c}_\ast$} (Z0)
(Z0) edge node[above]{${i_l^f}_\ast$} (XBar0)
(Z1) edge node[below]{${l_{\!_\red}}_\ast$} (XBar1);
\path[->>,dashed, font=\scriptsize]
(XBar0) edge node[above]{${ \epsilon_{\!_{i_l^f}}}$} (C0)
(XBar1) edge node[below]{${ \epsilon_{\!_{l_{\!_\red}}}}$} (C1);
\path[->,dotted, font=\scriptsize]
(C1) edge node[right]{${f_{\!_l}^\ast}$} (C0);
\end{tikzpicture}
\end{equation}
of cofibrations in $\scC$. Let $f_{\!_l}^\ast$ denote the unique morphism $C_{\!\!_F}(l_{\!_\red})\to C_{\!\!_F}(i_l^f)$, induced the universal property of cokernels, which makes the diagram {commute}; and denote the morphism $
f_{\!_l}^\ast\circ \iota_{\!_{l_{\!_\red}}}\colon {F^\ps(y) }\to C_{\!\!_F}(i_l^f)$ by $\rho_l^f$. Notice that the square \eqref{Pull:Back:CLosed:Immersions} is not a morphism of compactifications unless $f$ is an isomorphism.

\Lem{}{Assume that $f:x\oim y$ is an open immersion of $S$-schemes, let $l_k:w_k\oim q_k$ be a compactification of $y$, for $k=0,1$, and let $g:l_0\to l_1$ be a morphism of compactifications in $\Comp_{\!_S}(y)$. Then, $g$ induces a morphism of compactifications $g:i_{l_0}^f\to i_{l_1}^f$ in $\Comp_{\!_S}(x)$ for which $\rho_{l_1}^f={g}_\ast \circ \rho_{l_0}^f$.
}{Open:Independence}
\begin{proof}
	The morphism of compactifications $g:l_0\to l_1$ defines a morphism compactifications $g:{l_0}_{\!_\red}\to {l_1}_{\!_\red}$ in $\Comp_{\!_S}(y)$. Let $i_{l_k}^f:z_k\cim q_k$ be a complementary closed immersion to $j_{\!_{l_k}}\!\!\circ f$, endued with the reduced induced structure, and hence a compactification of $x$, and let $c_k:{w_k}_{\!_{\red}}\cim z_k$ be the unique closed immersion for which ${l_k}_{\!_{\red}}=i_{l_k}^f\!\!\circ c_k$, for $k=0,1$. Since $g:{l_0}_{\!_\red}\to {l_1}_{\!_\red}$ is a morphism of compactifications in $\Comp_{\!_S}(y)$, the morphism $\id_y$ is a base change in $\sftSch/S$ of $g$ along $j_{\!_{ {l_1} }}$, and hence $\id_x$ is a base change in $\sftSch/S$ of $g$ along $j_{\!_{i_{l_1}^f}}=j_{\!_{{l_1}}}\!\!\circ f$. Therefore, $g$ induces a morphism of compactifications $g:i_{l_0}^f\to i_{l_1}^f$ in $\Comp_{\!_S}(x)$, see Remark \ref{Rem:Cartesian:Square:Morphism:Compactifications}. Consider the commutative diagram
	\begin{equation}\label{Open:Independence:1}
	\begin{tikzpicture}[descr/.style={fill=white},scale=3]
	\node (Z0) at (.3,0) {$z_0$};
	\node (XBar0) at (1.2,0) {$q_0$};
	\node (W0) at (.3,-.9) {${w_0}_{\!_\red}$};
	\node (YBar0) at (1.2,-.9) {$q_0$};
	\node (Z1) at (.9,.3) {$z_1$};
	\node (XBar1) at (1.8,.3) {$q_1$};
	\node (W1) at (.9,-.6) {${w_1}_{\!_\red}$};
	\node (YBar1) at (1.8,-.6) {$q_1$};
	\path[right hook->,font=\scriptsize]
	(W1) edge node{$\diagup$}node[above]{${{{l_1}}_{\!_\red}}$}  node[descr,pos=.21,sloped] {} (YBar1)
	(W1) edge node{$\diagup$}node[right]{${{c_1}}$}  node[descr,pos=.70,sloped] {} (Z1)
	;
	\path[->,font=\scriptsize]	
	(W0) edge  (W1)
	;
	\path[->,font=\scriptsize]
	(YBar0) edge node[descr]{$\!\!\!\id\!\!\!$} (XBar0)
	(YBar1) edge node[descr]{$\!\!\!\id\!\!\!$} (XBar1)
	;
	\path[->,font=\scriptsize]
	(Z0) edge (Z1)
	(XBar0) edge node[descr]{$\!\!\!g\!\!\!$} (XBar1)
	(YBar0) edge node[descr]{$\!\!\!g\!\!\!$} (YBar1)
	;

	\path[right hook->,font=\scriptsize]
	(W0) edge node{$\diagup$}node[left]{${{c_0}}$} (Z0)
	(Z0) edge node{$\diagup$}node[below]{${i_{l_0}^f}$} (XBar0)
	(W0) edge node{$\diagup$}node[below]{${{{l_0}_{\!_\red}}}$} (YBar0)
	(Z1) edge node{$\diagup$}node[above]{${i_{l_1}^f}$} (XBar1)
	;
	\end{tikzpicture}
	\end{equation}
	of proper $S$-schemes, where the commutativity of the left face is a result of the commutativity of the other faces and having ${i_{l_1}^f}$ a monomorphism of proper $S$-schemes. The diagram \eqref{Open:Independence:1} induces the solid commutative diagram
	\[
	\begin{tikzpicture}[descr/.style={fill=white},scale=3]
	\node (Z0) at (0,0) {$F(z_0)$};
	\node (XBar0) at (1.5,0) {$F(q_0)$};
	\node (Ci0) at (3,0) {$C_{\!\!_F}(i_{l_0}^f)$};
	\node (W0) at (0,-1.2) {$F({w_0}_{\!_\red})$};
	\node (YBar0) at (1.5,-1.2) {$F(q_0)$};
	\node (Cj0) at (3,-1.2) {$C_{\!\!_F}({l_0}_{\!_\red})$};
	\node (Z1) at (1.1,.4) {$F(z_1)$};
	\node (XBar1) at (2.6,.4) {$F(q_1)$};
	\node (Ci1) at (4.1,.4) {$C_{\!\!_F}(i_{l_1}^f)$};
	\node (W1) at (1.1,-.8) {$F({w_1}_{\!_\red})$};
	\node (YBar1) at (2.6,-.8) {$F(q_1)$};
	\node (Cj1) at (4.1,-.8) {$C_{\!\!_F}({l_1}_{\!_\red})$};
	\path[>->,font=\scriptsize]
	(W1) edge node[descr]{$\!\!\! {{{l_1}_{\!_\red}}}_\ast \!\!\!$} node[descr,pos=.13] {} (YBar1)
	(W0) edge node[descr]{$\!\!\!{{c_0}}_\ast\!\!\!$} (Z0)
	(W1) edge node[descr]{$\!\!\!{{c_1}}_\ast\!\!\!$} node[descr,pos=.7,sloped] {} (Z1)
	(Z0) edge node[descr]{$\!\!\! {i_{l_0}^f}_\ast \!\!\!$} (XBar0)
	(W0) edge node[descr]{$\!\!\! {{{l_0}_{\!_\red}}}_\ast \!\!\!$} (YBar0)
	(Z1) edge node[descr]{$\!\!\!{i_{l_1}^f}_\ast\!\!\!$} (XBar1)
	
	;
	\path[->,font=\scriptsize]
	(YBar0) edge node[descr]{$\!\!\!\id\!\!\!$} (XBar0)
	(YBar1) edge node[descr]{$\!\!\!\id\!\!\!$}   node[descr,pos=.70,sloped] {}(XBar1)
	;
	\path[->,font=\scriptsize]
	(Z0) edge (Z1)
	(W0) edge (W1)
	(XBar0) edge node[descr]{$\!\!\!{g}_\ast\!\!\!$} (XBar1)
	(YBar0) edge node[descr]{$\!\!\!{g}_\ast\!\!\!$} (YBar1)
	;
	
	\path[->>,dashed, font=\scriptsize]
	(XBar0) edge node[descr]{$\!\!\!{ \epsilon_{\!_{i_{l_0}^f}}} \!\!\!$} (Ci0)
	(YBar0) edge node[descr]{$\!\!\! { \epsilon_{\!_{{l_0}_{\!_\red}}}}\!\!\!$} (Cj0)
	(XBar1) edge node[descr]{$\!\!\!{ \epsilon_{\!_{i_{l_1}^f}}} \!\!\!$} (Ci1)
	(YBar1) edge node[descr]{$\!\!\! { \epsilon_{\!_{{l_1}_{\!_\red}}}}\!\!\!$}  node[descr,pos=.195,sloped] {} (Cj1)
	;
	\path[->,dotted, font=\scriptsize]
	(Cj0) edge node[descr]{$\!\!\!f_{l_0}^\ast\!\!\!$} (Ci0)
	(Cj1) edge node[descr]{$\!\!\!f_{l_1}^\ast\!\!\!$} (Ci1)
	(Ci0) edge node[descr]{$\!\!\!{g}_\ast\!\!\!$} (Ci1)
	(Cj0) edge node[descr]{$\!\!\!{g}_\ast\!\!\!$} (Cj1)
	;
	\end{tikzpicture}
	\]
	in $\scC$. Since ${ \epsilon_{\!_{{l_0}_{\!_\red}}}}\!\!\!$ is an epimorphism in $\scC$, the universal propriety of cokernels implies that $f_{l_1}^\ast\circ {g}_\ast={g}_\ast \circ f_{l_0}^\ast$, \ie the whole diagram commutes. Then, one has
	\[
	\rho_{l_1}^f=f_{l_1}^\ast\circ \iota_{\!_{{l_1}_{\!_\red}}}=f_{l_1}^\ast\circ {g}_\ast\circ \iota_{\!_{{l_0}_{\!_\red}}}={g}_\ast \circ f_{l_0}^\ast\circ \iota_{\!_{{l_0}_{\!_\red}}}={g}_\ast \circ \rho_{l_0}^f.
	\]
\end{proof}
This proof, in particular, shows that there exists a faithful (not necessarily full) functor $\theta^f:\Comp_{\!_S}(y)\to \Comp_{\!_S}(x)$ that sends a compactification $l$ of $y$ to the compactification $i_l^f$ of $x$, and sends a morphism of compactifications $g:l_0\to l_1$ in $\Comp_{\!_S}(y)$ to the morphism of compactifications $g:i_{l_0}^f\to i_{l_1}^f$ in $\Comp_{\!_S}(x)$. Let $C_{\!\!_{F,x}}^f$ be the composition of the functor $C_{\!\!_{F,x}}$, given in Remark \ref{Rem:Cofiltered:Extension}, with the functor $\theta^f$. The universal properties of limits induces a canonical morphism $\vartheta^f:\lim C_{\!\!_{F,x}}\to \lim C_{\!\!_{F,x}}^f$ in $\scC$. Since $F$ is a \cdp-functor, the morphism $\vartheta^f$ is an isomorphism, which allows one to deduce the following corollary.

\Cor{}{Assume that $f:x\oim y$ is an open immersion of $S$-schemes. Then, there exists a unique morphism $f^!:F^\ps(y)\to F^\ps(x)$ in $\scC$ for which
	\begin{equation}\label{PS:Open:Pullback}
	{\iota_{\!_{i_l^f}}}\circ f^{!}=f_{l}^\ast\circ \iota_{\!_{l_{\!_\red}}},
	\end{equation}
	for every compactification $l$ of $y$.
}{Compact:Support:Functor:Open:Immersions}
The uniqueness of the morphism $f^!$, for an open morphism $f$ of $S$-schemes, implies the functoriality of pullbacks along open immersions. That is, for open immersions $f:x\oim y$ and $g:y\oim z$ of $S$-schemes, one has ${(g\circ f)}^!=f^! \circ g^!$ and ${(\id_x)}^!=\id_{\!_{\lim C_{\!\!_{F,x}}}}$.

\Cor{}{{The functor $F^\ps$, in Definition \ref{Def:Compact:Support:Functor:Proper}, extends to a functor $F^\ps:{\sftSch/S}^\prop_\open \to \scC$ that sends $f^\op$, for an open immersion $f:x\oim y$ of $S$-schemes, to the unique morphism $f^!:\lim C_{\!\!_{F,y}}\to \lim C_{\!\!_{F,x}}$ asserted vy {Corollary} \ref{Cor:Compact:Support:Functor:Open:Immersions}.}
}{Compactly:Supported:Counterpart}

\Rem{}{In contrast to pushforwards along proper morphisms, pullbacks along open immersions do not necessarily exist when $F$ is not a \cdp-functor, even if $\scC$ is has cofiltered limits. That is, when $F$ is not a \cdp-functor, the morphism $\vartheta^f$ is not necessarily an isomorphism, as the functor $\theta^f$ does not have to be co-cofinal. For instance, when $p$ is a proper $S$-scheme and $f:p\oim q$ is an open immersion of $S$-schemes that is not an isomorphism, the initial compactification $\emptyset_p:\emptyset_{\!_S}\cim p$ does not coincide with $i_l^f$ for any compactification $l$ of $q$.
}{}

\subsubsection{Base Change}\label{Sec:Base:Change}Open pullbacks and proper pushforwards satisfy the \indx{proper-base change formula}, as in the following lemma.
\Lem{}{Assume that $f:x\to y$ is a proper morphism and $j:y'\oim y$ is an open immersion of $S$-schemes, and let 
	\begin{equation}\label{Eq:Base:Change}
	\begin{tikzpicture}[descr/.style={fill=white}]
	\node (X0) at (0,0) {$x'$};
	\node (XBar0) at (2,0) {$x$};
	\node (X1) at (0,-2) {$y'$};
	\node (XBar1) at (2,-2) {$y$};
	\node at (.25,-.4) {$\ulcorner$};
	\path[->,font=\scriptsize]
	(XBar0) edge node[right]{$f$} (XBar1)
	(X0) edge node[left]{$f'$} (X1)
	;
	\path[right hook->,font=\scriptsize]
	(X0) edge node[descr]{$\!\!\!\circ\!\!\!$}node[above]{$j'$} (XBar0)
	(X1) edge node[descr]{$\!\!\!\circ\!\!\!$}node[below]{$j$} (XBar1);
	\end{tikzpicture}
	\end{equation}
	be a Cartesian square of $S$-schemes. Then, $j^!\circ f_\shriek=f'_\shriek\circ {j'}^!$.
}{Base:Change}
\begin{proof}	
	Let $l:w\cim q$ be a compactification of $y$. Since $f$ is proper, there exists a compactification $i:z\cim p$ of $x$ in $\Comp_{\!_S}(f,l)$ and a morphism of compactifications $g:i\to l$ such that $f$ is a base change in $\sftSch/S$ of $g$ along $j_{\!_l}$. Given that $j$ and $j'$ are open immersions, let $i_l^j:w'\cim q$ (\resp. $i_i^{j'}:z'\cim p$) be the compactification of $y'$ and (\resp. $x'$) induced from $l$ (\resp. $i$), as in \S.\ref{Sec:Pull:Back}. The morphism of compactification $g:i\to l$ defines a morphism compactifications $g:i_{\!_\red}\to l_{\!_\red}$. Since the square \eqref{Eq:Base:Change} is Cartesian in $\sftSch/S$, the morphism $f'$ is a base change in $\sftSch/S$ of $g$ along $j_{i_l^j}=j_{\!_l}\circ j$. Therefore, the morphism of compactification $g:i_{\!_\red}\cim l_{\!_\red}$ defines a morphism compactifications $g:i_i^{j'}\to i_l^j$, see Remark \ref{Rem:Cartesian:Square:Morphism:Compactifications}. Then, one has
	\[
	\iota_{\!_{i_l^j}}\circ j^! \circ f_\shriek=j_{\!_l}^\ast \circ \iota_{\!_{l_{\!_\red}}}\circ f_\shriek=j_{\!_l}^\ast \circ g_\ast\circ \iota_{\!_{i_{\!_\red}}}
	\andd
	\iota_{\!_{i_l^j}}\circ f'_\shriek \circ {j'}^!=g_\ast \circ \iota_{\!_{i_i^{j'}}} \circ {j'}^!=g_\ast \circ {j'_{\!_i}}^\ast\circ \iota_{\!_{i_{\!_\red}}}.
	\]
	Consider the commutative diagram \eqref{Base:Change:Diagram:Closed} of proper $S$-schemes, where the commutativity of the left square is a result of the commutativity of the other squares and having ${i_{l}^{j}}$ a monomorphism of proper $S$-schemes. 
	\begin{equation}\label{Base:Change:Diagram:Closed}
	\begin{tikzpicture}[descr/.style={fill=white}]
	\node (z) at (-2,0) {$z_{\!_\red}$};
	\node (X0) at (0,0) {$z'$};
	\node (XBar0) at (2,0) {$p$};
	\node (X1) at (0,-2) {$w'$};
	\node (w) at (-2,-2) {$w_{\!_\red}$};
	\node (XBar1) at (2,-2) {$q$};
	\path[->,font=\scriptsize]
	(XBar0) edge node[right]{$g$} (XBar1)
	(z) edge (w)
	;
	\path[->,font=\scriptsize]
	(X0) edge (X1)
	;
	\path[right hook ->,font=\scriptsize]
	(w) edge node{$\diagup$}(X1)
	(z) edge node{$\diagup$} (X0)
	(X1) edge node{$\diagup$}node[above]{$i_l^j$} (XBar1)
	(X0) edge node{$\diagup$}node[below]{$i_i^{j'}$} (XBar0)
	;
	\path[bend left,right hook ->,font=\scriptsize]
	(z) edge node{$\diagup$}node[above]{$i_{\!_\red}$} (XBar0)
	;
	\path[bend right,right hook ->,font=\scriptsize]
	(w) edge node{$\diagup$}node[below]{$l_{\!_\red}$} (XBar1)
	;
	\end{tikzpicture}
	\end{equation}
	The diagram \eqref{Base:Change:Diagram:Closed} induces the solid commutative diagram
	\[
	\begin{tikzpicture}[descr/.style={fill=white},scale=3]
	\node (Z0) at (0,0) {$F(z')$};
	\node (XBar0) at (1.5,0) {$F(p)$};
	\node (Ci0) at (3,0) {$C_{\!\!_F}(i_i^{j'})$};
	\node (W0) at (0,-1.2) {$F({z}_{\!_\red})$};
	\node (YBar0) at (1.5,-1.2) {$F(p)$};
	\node (Cj0) at (3,-1.2) {$C_{\!\!_F}(i_{\!_\red})$};
	\node (Z1) at (1.1,.4) {$F(w')$};
	\node (XBar1) at (2.6,.4) {$F(q)$};
	\node (Ci1) at (4.1,.4) {$C_{\!\!_F}(i_l^j)$};
	\node (W1) at (1.1,-.8) {$F({w}_{\!_\red})$};
	\node (YBar1) at (2.6,-.8) {$F(q)$};
	\node (Cj1) at (4.1,-.8) {$C_{\!\!_F}(l_{\!_\red})$};
	\path[>->,font=\scriptsize]
	(W1) edge node[descr,pos=.7] {}(Z1)
	(W0) edge (Z0)
	(Z0) edge node[descr]{$\!\!\! {i_i^{j'}}_\ast \!\!\!$} (XBar0)
	(W0) edge node[descr]{$\!\!\! {i_{\!_\red}}_\ast \!\!\!$} (YBar0)
	(Z1) edge node[descr]{$\!\!\!{i_l^j}_\ast\!\!\!$} (XBar1)
	(W1) edge node[descr]{$\!\!\! {l_{\!_\red}}_\ast \!\!\!$}  node[descr,pos=.145] {} (YBar1)
	;
	\path[->,font=\scriptsize]
	(YBar0) edge node[descr]{$\!\!\!\id\!\!\!$} (XBar0)
	(YBar1) edge node[descr]{$\!\!\!\id\!\!\!$}  node[descr,pos=.7,sloped] {} (XBar1)
	;
	\path[->>,dashed, font=\scriptsize]
	(YBar1) edge node[descr]{$\!\!\! { \epsilon_{\!_{l_{\!_\red}}}}\!\!\!$} node[descr,pos=.195,sloped] {}(Cj1)
	(XBar0) edge node[descr]{$\!\!\!{ \epsilon_{\!_{i_i^{j'}}}} \!\!\!$} (Ci0)
	(YBar0) edge node[descr]{$\!\!\! { \epsilon_{\!_{i_{\!_\red}}}}\!\!\!$}    (Cj0)
	(XBar1) edge node[descr]{$\!\!\!{ \epsilon_{\!_{i_l^j}}} \!\!\!$} (Ci1)
	;

	\path[->,dotted, font=\scriptsize]
	(Cj0) edge node[descr]{$\!\!\!{j'_{\!_{i}}}^\ast\!\!\!$} (Ci0)
	(Cj1) edge node[descr]{$\!\!\!j_{\!_{l}}^\ast\!\!\!$} (Ci1)
	(Ci0) edge node[descr]{$\!\!\!{g}_\ast\!\!\!$} (Ci1)
	(Cj0) edge node[descr]{$\!\!\!{g}_\ast\!\!\!$} (Cj1)
	;
	\path[->,font=\scriptsize]
	(Z0) edge (Z1)
	(W0) edge   (W1)
	(XBar0) edge node[descr]{$\!\!\!{g}_\ast\!\!\!$} (XBar1)
	(YBar0) edge node[descr]{$\!\!\!{g}_\ast\!\!\!$} (YBar1)
	;
	\end{tikzpicture}
	\]
	in $\scC$. Since ${ \epsilon_{\!_{{i}_{\!_\red}}}}\!\!\!$ is an epimorphism in $\scC$, the universal propriety of cokernels implies that $j_{\!_l}^\ast \circ g_\ast=g_\ast \circ {j'_{\!_i}}^\ast$. The morphism $\iota_{\!_{i_l^j}}$ is an isomorphism, as $F$ is a \cdp-functor. Hence, one has $j^! \circ f_\shriek=f'_\shriek \circ {j'}^!$.
\end{proof}

\Eg{}{
	Assume that $i:v\cim x$ is a closed immersion of $S$-schemes with complementary open immersion $j:u\oim x$. Then, one has a Cartesian square
	\[
	\begin{tikzpicture}[descr/.style={fill=white}]
	\node (X0) at (0,0) {$\emptyset_{\!_S}$};
	\node (XBar0) at (2,0) {$v$};
	\node (X1) at (0,-2) {$u$};
	\node (XBar1) at (2,-2) {$x$};
	\node at (.25,-.4) {$\ulcorner$};
	\path[right hook->,font=\scriptsize]
	(XBar0) edge node{$\diagup$} node[right]{$i$} (XBar1)
	(X0) edge node{$\diagup$} (X1)
	(X0) edge node[descr]{$\!\!\!\circ\!\!\!$} (XBar0)
	(X1) edge node[descr]{$\!\!\!\circ\!\!\!$}node[below]{$j$} (XBar1);
	\end{tikzpicture}
	\]
	in $\sftSch/S$. Thus, one has $	j^!\circ i_\shriek=0$ as $F^\ps(\emptyset_{\!_S})\cong F(\emptyset_{\!_S})\cong\0$. {In fact, the sequence $F^\ps(v)\stackrel{i_\shriek}{\rightarrowtail} F^\ps(x) \stackrel{j^!}{\twoheadrightarrow} F^\ps (u)$ is a cofibre sequence in $\scC$, as in Proposition \ref{Prop:Excision}.
	}
}{Closed:Open:Sequence:0}

\Eg{}{{Recall Example \ref{Ex:Coproduct:Preservation}, and} assume that $j_x:x\incl x\bigsqcup y$ is a closed open immersion of $S$-schemes with complementary closed open immersion $j_y:y\incl x\bigsqcup y$. Then, there exists a Cartesian square
	\[
	\begin{tikzpicture}[descr/.style={fill=white}]
	\node (X0) at (0,0) {$x$};
	\node (XBar0) at (2,0) {$x$};
	\node (X1) at (0,-2) {$x$};
	\node (XBar1) at (2,-2) {$x\bigsqcup y$};
	\node at (.25,-.4) {$\ulcorner$};
	\path[->,font=\scriptsize]
	(X0) edge node[left]{$\id_x$} (X1)
	(X0) edge node[above]{$\id_x$} (XBar0);
	\path[right hook->,font=\scriptsize]
	(XBar0) edge node[right]{$j_x$} (XBar1)
	
	(X1) edge node[below]{$j_x$} (XBar1);
	\end{tikzpicture}
	\]
	in $\sftSch/S$, and hence ${j_x}^!\circ {j_x}_\shriek=(\id_{x})_\shriek\circ (\id_{x})^!=\id_{F^\ps(x)}$. Similarly ${j_y}^!\circ {j_y}_\shriek=\id_{F^\ps(y)}$. By Example \ref{Ex:Closed:Open:Sequence:0}, one also has ${j_x}^!\circ {j_y}_\shriek=0$ and ${j_y}^!\circ {j_x}_\shriek=0$. {Therefore, when $\scC$ is additive, $F^\ps (x\bigsqcup y)$ is a direct sum in $\scC$ of $F^\ps (x)$ and $F^\ps(y)$}.
}{Direct:Sum:Preservation}

\subsubsection{Excision}\label{Sec:Excision}
A functor $G:{\sftSch/S}^\prop_\open\to \scC$, to a {\waldhausen category}, induces a group homomorphism $\Kg_0(\sftSch/S)\to \Kg_0(\scC)$ only when the evident composite map $\Ob({\sftSch/S}^\prop_\open)\to \Ob(\scC)\to \Kg_0(\scC)$ respects the scissors relations. That holds when $G$ satisfies the {excision property}. We will see below that the functor $F^\ps$, given in Corollary \ref{Cor:Compactly:Supported:Counterpart}, satisfies the {excision property}, and hence it induces an \euler-\poincare characteristic $\mu_{\!_F}:\Kg_0(\sftSch/S)\to \Kg_0(\scC)$.

\Prop{}{The functor $F^\ps$, given in Corollary \ref{Cor:Compactly:Supported:Counterpart}, satisfies the {excision property}, \ie for every closed immersion $i:v\cim x$ of $S$-schemes with complementary open immersion $j:u\oim x$, the sequence $
	F^\ps (v)\stackrel{i_\shriek}{\rightarrowtail} F^\ps(x) \stackrel{j^!}{\twoheadrightarrow} F^\ps (u)$ is a cofibre sequence in $\scC$. In particular, $j^!$ is an epimorphism in $\scC$.
}{Excision}
\begin{proof}Assume that $i:v\cim x$ is a closed immersion of $S$-schemes with complementary open immersion $j:u\oim x$, and let $l:w\cim q$ be a compactification of $x$ in which the underlying scheme of $w$ is reduced. The open immersion $j_{\!_l}\circ j:u\oim q$ induces a compactification $l':w'\cim q$ of $u$ in which the underlying scheme of $w'$ is reduced, \ie $l'=i_l^j$ using the notation of \S.\ref{Sec:Pull:Back}. Hence, there exists a unique closed immersion $c:w\cim w'$ for which $l=l'\circ c$. 
	
	\medskip
	
	On the other hand, let $i':p\cim q$ be the closed immersion of the scheme-thoracic image of $j_{\!_l}\circ i$. Then, Lemma \ref{Lem:Open:Closed:Closed:Open} implies the existence of an open immersion $j'_{\!_l}:v\oim p$ for which the solid square containing it in the commutative diagram
	\begin{equation}\label{Excision:Diag:1}
	\begin{tikzpicture}[descr/.style={fill=white},>=angle 90,scale=2,text height=1.5ex, text depth=.25ex,row sep=4em, column sep=4em]
	
	\node (F11) at (0,0) {$u$};
	\node (F12) at (1,0) {$x$};
	\node (F13) at (2,0) {$q$};
	\node (F14) at (3,0) {$w'$};
	\node (F15) at (4,0) {$w$};
	
	\node (F21) at (0,-1) {$\emptyset_{\!_S}$};
	\node at (0.1,-.8) {$\llcorner$};
	\node (F22) at (1,-1) {$v$};
	\node at (1.1,-.8) {$\llcorner$};
	\node (F23) at (2,-1) {$p$};
	\node (F24) at (3,-1) {${p}'$};
	\node at (2.9,-.8) {$\lrcorner$};
	\node (F25) at (4,-1) {$z$};
	\node at (3.9,-.8) {$\lrcorner$};
	
	\draw[right hook->,font=\scriptsize]
	(F21) edge node{$\diagup$}node[left] {$\emptyset_u$}(F11)
	(F22) edge node{$\diagup$}node[left] {$i$}(F12)
	(F23) edge node{$\diagup$}node[left] {$i'$}(F13)
	(F11) edge node[descr]{$\!\!\!\circ\!\!\!$}node[above] {$j$} (F12)
	(F21) edge node[descr]{$\!\!\!\circ\!\!\!$} (F22)
	(F12) edge node[descr]{$\!\!\!\circ\!\!\!$}node[above] {$j_{\!_l}$} (F13)
	(F22) edge node[descr]{$\!\!\!\circ\!\!\!$}node[below] {$j'_{\!_l}$} (F23)
	;
	\draw[left hook->,font=\scriptsize]	
	(F14) edge node{$\diagup$}node[above] {$l'$}(F13)	
	(F15) edge node{$\diagup$}node[above] {$c$}(F14)
	;
	\draw[dashed, left hook->,font=\scriptsize]	
	(F24) edge node{$\diagup$}node[below] {$\pi'_{p}$}(F23)	
	(F25) edge node{$\diagup$}node[below] {$e$}(F24)
	
	(F24) edge node{$\diagup$}node[left] {$\pi_{w'}$}(F14)	
	(F25) edge node{$\diagup$}node[left] {$\pi_{w}$}(F15)
	;
	\end{tikzpicture}
	\end{equation}	
	of $S$-schemes is Cartesian. The leftmost square in \eqref{Excision:Diag:1} is also Cartesian in $\sftSch/S$, by the definition of $u$. Let $(\pi_{w'},p',\pi'_p)$ (\resp. $(\pi_{w},z,\pi_p)$) be the pullback of the cospan of $i'$ and $l'$ (\resp.  $i'$ and $l=l'\circ c$). Then, there exists a unique morphism $e:z\to p'$ of proper $S$-schemes that makes the diagram commute. In particular, one has $\pi_p=\pi'_p\circ e$ and $e$ is a closed immersion.
	
	\medskip
	
	{Since open complements are closed under pullbacks} and $j_{\!_l}$ is a complementary open immersion to $l$, one finds that $j'_{\!_l}$ is a complementary open immersion to $\pi_p$. Hence, $\pi_p:z\cim p$ is a compactification in $\Comp_{\!_S}(i,l)$ and $(\pi_w,i'):\pi_p\to l$ is a morphism of compactifications such that $i$ is a base change in $\sftSch/S$ of $i'$ along $j_{\!_l}$. Also, since $j_{\!_l}\circ j$ is a complementary open immersion to $l'$, the projection $\pi'_p$ is a surjective closed immersion. Hence, $\pi'_p:p'\cim p$ is a compactification in $\Comp_{\!_S}(\emptyset_u,l')$ and $(\pi_{w'},i'):\pi'_p\to l'$ is a morphism of compactifications such that $\emptyset_u$ is a base change in $\sftSch/S$ of $i'$ along $j_{\!_{l'}}$.

	\medskip
	The {Pullback Lemma} implies that the square of proper $S$-schemes
	\begin{equation}\label{Excision:Abstract:Blow:Up}
	\begin{tikzpicture}[descr/.style={fill=white},>=angle 90,scale=2,text height=1.5ex, text depth=.25ex,row sep=4em, column sep=4em]
	
	\node (F14) at (3,0) {$w'$};
	\node (F15) at (3,1) {$w$};
	
	\node (F24) at (2,0) {${p}'$};
	
	\node (F25) at (2,1) {$z$};
	\node at (2.1,.80) {$\ulcorner$};
	
	\draw[right hook->,font=\scriptsize]	
	
	(F15) edge node{$\diagup$}node[right] {$c$}(F14)
	
	(F25) edge node{$\diagup$}node[left] {$e$}(F24)
	
	(F24) edge node{$\diagup$}node[below] {$\pi_{w'}$}(F14)	
	(F25) edge node{$\diagup$}node[above] {$\pi_{w}$}(F15)
	;
	\end{tikzpicture}
	\end{equation}
	is Cartesian. Since $\pi'_p$ is a {surjective} closed immersion, the morphism $c$ induces an isomorphism
	\[
	w'{\setminus} \pi_{w'}\cong
	\left(q{\setminus} (j_{\!_l}\circ j )\right)_{\!_\red} \!\!{\setminus} \pi_{w'}=
	\left(\left(q{\setminus} (j_{\!_l}\circ j) \right) {\setminus} i' \right)_{\!_\red}=
	\left(\left(q{\setminus} j_{\!_l}\right){\setminus} i'\right)_{\!_\red}=\left(q{\setminus} j_{\!_l}\right)_{\!_\red}\!\!{\setminus} \pi_w\cong w {\setminus} \pi_w
	\]
	of $S$-schemes. Therefore, the square \eqref{Excision:Abstract:Blow:Up} is a \cdp-square of proper $S$-schemes, and hence the functor $F$ sends it to a pushout square in $\scC$.
	
	\medskip
	
	Since $(\pi_w,i'):\pi_p\to l$ (\resp. $(\pi_{w'},i'):\pi'_p\to l'$) is a morphism of compactifications and $i$ (\resp. $\emptyset_u$) is a base change of $i'$ along $j_{\!_l}$ (\resp. $j_{\!_{l'}}=j_{\!_l}\circ j$), there exists a commutative diagram 
	\begin{equation}\label{Excision:Pushouts}
	\begin{tikzpicture}[descr/.style={fill=white},>=angle 90,scale=2,text height=1.5ex, text depth=.25ex,row sep=4em, column sep=4em]
	\node (F14) at (1.5,0) {$F(w')$};
	\node (F15) at (1.5,1) {$F(w)$};
	\node (F24) at (0,0) {$F({p}')$};
	\node (F25) at (0,1) {$F(z)$};
	\node at (1.25,.25) {$\lrcorner$};
	\node (G14) at (3.75,-.5) {$F(q)$};
	\node (G15) at (3.75,.5) {$F(q)$};
	\node (G24) at (2.25,-.5) {$F(p)$};
	\node (G25) at (2.25,.5) {$F(p)$};
	\node (H14) at (6,-1) {$F^\ps(u)$};
	\node (H15) at (6,0) {$F^\ps(x)$};
	\node (H24) at (4.5,-1) {$F^\ps(\emptyset_{\!_S})$};
	\node (H25) at (4.5,0) {$F^\ps(v)$};
	\draw[>->,font=\scriptsize]	
	(F15) edge node[descr,pos=.6]{$\!\!\!c_\ast\!\!\!$}node[descr,pos=.265,sloped] {} (F14)
	(F25) edge node[left,pos=.6] {$e_\ast$}(F24)
	(F24) edge node[above] {${\pi_{w'}}_\ast$}(F14)	
	(F25) edge node[above] {${\pi_{w}}_\ast$}(F15)
	;
	\draw[>->,font=\scriptsize]	
	(F15) edge node[above] {$l_\ast$}(G15)
	(F25) edge node[descr]{$\!\!\!{\pi_{p}}_\ast\!\!\!$}(G25)
	(F24) edge node[descr]{$\!\!\!{\pi'_{p}}_\ast\!\!\!$}(G24)	
	(F14) edge node[descr]{$\!\!\!l'_\ast\!\!\!$}node[descr,pos=.265,sloped] {} (G14)
	;
	\draw[font=\scriptsize]
	(G15.265) edge node[descr,pos=.27,sloped] {} (G14.95)
	(G15.275) edge node[descr,pos=.27,sloped] {} (G14.85)
	;
	\draw[->>,font=\scriptsize]	
	(G15) edge node[descr]{$\!\!\!\epsilon'_{\!_{l}}\!\!\!$} (H15)
	(G25) edge node[descr]{$\!\!\!\epsilon'_{\!_{\pi_p}}\!\!\!$} (H25)
	(G24) edge node[descr]{$\!\!\!\epsilon'_{\!_{\pi'_p}}\!\!\!$} (H24)	
	(G14) edge node[descr]{$\!\!\!\epsilon'_{\!_{l'}}\!\!\!$}   node[descr,pos=.29,sloped] {}          (H14)
	;
	\draw[->,font=\scriptsize]	
	(H15) edge node[descr,pos=.4]{$\!\!\!j^{!}\!\!\!$}(H14)
	(H25) edge node[descr,pos=.4]{$\!\!\!0\!\!\!$} (H24)
	(H24) edge node[descr]{$\!\!\!0\!\!\!$}(H14)	
	(H25) edge node[descr]{$\!\!\!i_\shriek\!\!\!$}(H15)
	;
	\draw[>->,font=\scriptsize]	
	(G24) edge node[descr]{$\!\!\!i'_\ast\!\!\!$}(G14)	
	(G25) edge node[descr]{$\!\!\!i'_\ast\!\!\!$}(G15)
	;
	\draw[font=\scriptsize]
	(G25.275) edge (G24.85)
	(G25.265) edge (G24.95)
	;
	\end{tikzpicture}
	\end{equation}
	in $\scC$, where $\epsilon'_{\!_l}=\iota_{\!_l}^\inv\circ \epsilon_{\!_l}$, $\epsilon'_{\!_{l'}}=\iota_{\!_{l'}}^\inv\circ \epsilon_{\!_{l'}}$, $\epsilon'_{\!_{\pi_p}}=\iota_{\!_{\pi_p}}^\inv\circ \epsilon_{\!_{\pi_p}}$, and $\epsilon'_{\!_{\pi'_p}}=\iota_{\!_{\pi'_p}}^\inv\circ \epsilon_{\!_{\pi'_p}}$. Since $\pi'_p$ is a {surjective} closed immersion of proper $S$-schemes, the morphism ${\pi'_p}_\ast$ is an isomorphism. Also, $F^\ps(\emptyset_{\!_S})\cong \0$. 
	
	\medskip

	Since the functor $F$ send the square \eqref{Excision:Abstract:Blow:Up} to a pushout square in $\scC$ and the morphisms $\epsilon'_{\!_{l}}$ and $\epsilon'_{\!_{l'}}$ are epimorphisms in $\scC$, a diagram chase of \eqref{Excision:Pushouts} shows that the sequence $F^\ps (v)\stackrel{i_\shriek}{\rightarrowtail} F^\ps(x) \stackrel{j^!}{\twoheadrightarrow} F^\ps (u)$ is a cokernel sequence in $\scC$. Hence, it is a cofibre sequence as $i_\shriek$ is a cofibration in $\scC$, by Proposition \ref{Prop:PS:Properties}.\ref{itm:PS:C}.
\end{proof}

\subsubsection{Weak Monoidal}\label{MSPS:Weak:Monodial}
When the \waldhausen category $\scC$ is {symmetric monoidal} and the \cdp-functor $F:\PropSch/S\to \scC$ is weak monoidal with respect to the Cartesian product of proper $S$-schemes, we show that $F^\ps$ is also weak monoidal. The weak monoidality of {the} functor {$F^\ps$} is a formal consequence of Proposition \ref{Prop:Excision} and the weak monoidality of {$F$}, which is based on the {following} lemma.

\Lem{}{Let the diagram \eqref{Mon:Dia} be a commutative diagram in a \waldhausen category $\scC$, in which the vertical and horizontal sequences are {{cokernel}} sequences, $(\varsigma,F,\sigma)$ is a pushout of the span $F_{2,1}\stackrel{\sigma_1}{\longleftarrow}F_{1,1}\stackrel{\varsigma_1}{\longrightarrow}F_{1,2}$, and the morphism $\lambda$ is the unique morphism $F\to F_{2,2}$ induced by the universal property of pushouts making the diagram commute. Then, the sequence $F\stackrel{\lambda}{\longrightarrow} F_{2,2} \stackrel{\pi}{\longrightarrow}F_{3,3}$ is a {{cokernel}} sequence in $\scC$. 
	\begin{equation}\label{Mon:Dia}
	\begin{tikzpicture}[descr/.style={fill=white},>=angle 90,scale=1.5,text height=1.5ex, text depth=.25ex,row sep=4em, column sep=4em]
	\node (F11) at (0,0) {$F_{1,1}$};
	\node (F12) at (1.5,0) {$F_{1,2}$};
	\node (F13) at (3,0) {$F_{1,3}$};
	\node (F21) at (0,-1.5) {$F_{2,1}$};
	\node (F22) at (1.5,-1.5) {$F_{2,2}$};
	\node (F23) at (3,-1.5) {$F_{2,3}$};
	\node (F33) at (3,-3) {$F_{3,3}$};
	\node (F) at (1,-1) {$F$};
	\node at (.85,-.85) {$\lrcorner$};
	\draw[>->,font=\scriptsize]
	(F11) edge node[above] {$\varsigma_{1}$}(F12)
	(F21) edge node[below] {$\varsigma_{2}$} (F22)
	(F13) edge node[right] {$\sigma_{3}$}(F23)
	(F21) edge node[above] {$\varsigma$}(F)
	;
	\draw[->,font=\scriptsize]	
	
	(F11) edge node[left] {$\sigma_{1}$} (F21)
	(F12) edge node[right] {$\sigma_{2}$}(F22)
	
	(F12) edge node[left] {$\sigma$}(F)
	(F22) edge node[descr]{$\!\!\!\pi\!\!\!$}(F33)
	;
	\draw[->,dashed,font=\scriptsize]
	(F) edge node[above] {$\lambda$}(F22)
	;
	\draw[->>,font=\scriptsize]
	(F12) edge node[above] {$\pi_{1}$} (F13)
	(F22) edge node[above]{$\pi_{2}$} (F23)
	(F23) edge node[right] {$\pi_3$} (F33)
	;
	\end{tikzpicture}
	\end{equation}
}{Cofibrations:Cofibre:Square}
\begin{proof}
	The statement of the lemma follows from a diagram chase of \eqref{Mon:Dia}.
\end{proof}

\Prop{}{Let $(\scC,\wedge,\1)$ be a symmetric monoidal \waldhausen category, and suppose that the \cdp-functor $F:(\PropSch/S,\times,\id_{\!_S})\to (\scC,\wedge,\1)$ is weak monoidal. Then, the functor $F^\ps:{\sftSch/S}^\prop_\open\to \scC$ is weak monoidal. Moreover, $F^\ps$ is strong monoidal when $F$ is.
}{Monoidal}
\begin{proof}Since $F$ is weak monoidal, its coherence morphisms $\phi_{p,q}:F(p)\wedge F(q)\to F(p\times q)$ and $\phi_{\!_S}:\1\to F(\id_{\!_S})$ are weak equivalences in $\scC$, for every pair of proper $S$-schemes $p$ and $q$. Let $\phi^\ps_{\!_S}:\1\to F^\ps(\id_{\!_S})$ be the composite {weak equivalence} $\phi^\ps_{\!_S}\coloneqq \varphi_{\!_S}\circ \phi_{\!_S}$ in $\scC$, where $\varphi$ is the {natural isomorphism} {asserted} by Corollary \ref{Cor:Comparison:Morphism}.
	
	\medskip

	For $k=0,1$, assume that $x_k$ is an $S$-scheme, let $i_k:z_k\cim p_k$ be a compactification of $x_k$ with complementary open immersion $j_{\!_{i_k}}:x_k\oim p_k$. {The functors $-\wedge C$ and $C\wedge-$ preserve {{cofibre} sequences} for every object $C\in \scC$, see Definition \ref{Def:Waldhausen:Model:Structure}.} Thus, Proposition \ref{Prop:Excision} induces a {cofibre} sequence 
	\[
	\begin{tikzpicture}[descr/.style={fill=white},>=angle 90,scale=2,text height=1.5ex, text depth=.25ex,row sep=4em, column sep=4em]
	\node (F11) at (0,0) {$F^\ps(z_k)\wedge F^\ps(y)$};
	\node (F12) at (2.5,0) {$F^\ps(p_k)\wedge F^\ps(y)$};
	\node (F13) at (5,0) {$F^\ps(x_k)\wedge F^\ps( y)$};
	\draw[>->,font=\scriptsize]
	(F11) edge node[above] {${i_k}_\shriek$}(F12);
	\draw[->>,font=\scriptsize]
	(F12) edge node[above] {$j_{\!_{i_k}}^!$}(F13);
	\end{tikzpicture}
	\]
	in $\scC$, for every $S$-scheme $y$. Since both $p_k$ and $z_k$ are proper $S$-schemes, and $\varphi$ is a natural isomorphism, the sequence
	\[
	\begin{tikzpicture}[descr/.style={fill=white},>=angle 90,scale=2,text height=1.5ex, text depth=.25ex,row sep=4em, column sep=4em]
	\node (F11) at (0,0) {$F(z_k)\wedge F(y)$};
	\node (F12) at (2.5,0) {$F(p_k)\wedge F( y)$};
	\node (F13) at (5,0) {$F^\ps(x_k)\wedge F( y)$};
	\draw[>->,font=\scriptsize]
	(F11) edge node[above] {${i_k}_\ast$}(F12);
	\draw[->>,font=\scriptsize]
	(F12) edge node[above] {$j_{\!_{i_k}}^!$} (F13);
	\end{tikzpicture}
	\]
	is a {cofibre} sequence in $\scC$ for a proper $S$-scheme $y$. 
	
	\medskip
	
	On the one hand, the monoidal product bifunctor $\wedge:\scC\times\scC\to \scC$ induces the solid commutative diagram \eqref{Monoidal:Diagram} in $\scC$, in which horizontal and vertical sequences are {cofibre} sequences. Let $\lambda:\underline{F}\to F(p_0)\wedge F( p_1)$ be the pushout-product of ${i_0}_\ast$ and ${i_1}_\ast$, and let $\pi\coloneqq j_{\!_{i_0}}^!\circ j_{\!_{i_1}}^!$. Since $\scC$ is a {symmetric monoidal \waldhausen category} and ${i_0}_\ast$ and ${i_1}_\ast$ are cofibrations in $\scC$, the morphism $\lambda$ is also a cofibration in $\scC$. Thus, {Lemma} \ref{Lem:Cofibrations:Cofibre:Square} implies that the sequence
	\[
	\underline{F}\stackrel{\lambda}{\longrightarrow}F(p_0)\wedge F( p_1) \stackrel{\pi}{\longrightarrow}F^\ps(x_0)\wedge F^\ps( x_1)
	\]
	is a {cofibre} sequence in $\scC$.
	
	\medskip
	
		On the other hand, there exits a pushout in $\PropSch/S$ of the span of the closed immersions $\id_{z_0}\times i_1$ and $i_0\times \id_{z_1}$. In fact, there exists a bicartesian square
		
	\begin{equation}\label{Monodial:Bicartesian}
	\begin{tikzpicture}[descr/.style={fill=white},ampersand replacement=\&]
	\node (m-1-1) at (0,0) {$z_0\times z_1$};
	\node (m-1-2) at (2,0) {$z_0\times p_1$};
	\node (m-2-1) at (0,-2) {$p_0\times z_1$};
	\node (m-2-2) at (2,-2) {$z,$};
	\node at (.25,-.4) {$\ulcorner$};
	\node at (1.75,-1.6) {$\lrcorner$};
	\path[right hook->,font=\scriptsize]
	(m-1-1) edge node{$\diagup$}node	[left]	{$i_0\times\id_{z_1}$} (m-2-1)
	(m-1-2) edge node{$\diagup$}node	[right]	{$i'_0$} (m-2-2)
	(m-2-1) edge node{$\diagup$}node	[below]	{$i'_1$} (m-2-2)
	(m-1-1) edge node{$\diagup$}node	[above]	{$\id_{z_0}\times i_1 $} (m-1-2);	
	\end{tikzpicture}
	\end{equation}
	of closed immersions of proper $S$-schemes, see \cite[Th.3.11]{Sch:05:GSSWCP} and \cite[\S.2]{Cam:17:KTSV}. In particular, the square \eqref{Monodial:Bicartesian} is a \cdp-square of proper $S$-schemes, and hence it is mapped by $F$ to a pushout square of cofibrations in $\scC$. 
	
	\begin{equation}\label{Monoidal:Diagram}
	\begin{tikzpicture}[descr/.style={fill=white},>=angle 90,scale=2,text height=1.5ex, text depth=.25ex,row sep=4em, column sep=4em]
	\node (F11) at (0,0) {$F(z_0)\wedge F( z_1)$};
	\node (F12) at (2,0) {$F(z_0)\wedge F(p_1)$};
	\node (F13) at (4,0) {$F(z_0)\wedge F^\ps( x_1)$};
	\node (F21) at (0,-2) {$F(p_0)\wedge F( z_1)$};
	\node (F22) at (2,-2) {$F(p_0)\wedge F( p_1)$};
	\node (F23) at (4,-2) {$F(p_0)\wedge F^\ps( x_1)$};
	\node (F31) at (0,-4) {$F^\ps(x_0)\wedge F( z_1)$};
	\node (F32) at (2,-4) {$F^\ps(x_0)\wedge F( p_1)$};
	\node (F33) at (4,-4) {$F^\ps(x_0)\wedge F^\ps( x_1)$};
	\node (F) at (1.4,-1.4) {$\underline{F}$};
	\node at (1.3,-1.3) {$\lrcorner$};
	\draw[>->,font=\scriptsize]
	(F11) edge node[above] {${i_{1}}_\ast$}(F12)
	(F21) edge node[below] {${i_{1}}_\ast$} (F22)
	(F31) edge node[below] {${i_{1}}_\ast$} (F32)
	(F11) edge node[left] {${i_{0}}_\ast$} (F21)
	(F12) edge node[right] {${i_{0}}_\ast$}(F22)
	(F13) edge node[right] {${i_{0}}_\ast$}(F23)
	;
	\draw[->,font=\scriptsize]
	(F22) edge node[descr]{$\!\!\!\pi\!\!\!$}(F33)
	;
	\draw[->,dashed,font=\scriptsize]
	(F21) edge node[above] {$\lambda_1$}(F)
	(F12) edge node[left] {$\lambda_0$}(F)
	;
	\draw[->,dotted,font=\scriptsize]
	(F) edge node[above] {$\lambda$}(F22)
	;
	\draw[->>,font=\scriptsize]
	(F12) edge node[above] {$j_{\!_{i_1}}^!$} (F13)
	(F22) edge node[above]{$j_{\!_{i_1}}^!$} (F23)
	(F32) edge node[below] {$j_{\!_{i_1}}^!$} (F33)
	(F21) edge node[left] {$j_{\!_{i_0}}^!$} (F31)
	(F22) edge node[left]{$j_{\!_{i_0}}^!$} (F32)
	(F23) edge node[right] {$j_{\!_{i_0}}^!$} (F33)
	;
	\end{tikzpicture}
	\end{equation}

	Since $\scC$ is a Waldhausen category and $F$ is weak monoidal, the unique morphism $\lambda':\underline{F}\to F(z)$, that makes the diagram \eqref{Monodial:Bicartesian:Image} commute, is a weak equivalence in $\scC$.	
	\begin{equation}\label{Monodial:Bicartesian:Image}
	\begin{tikzpicture}[descr/.style={fill=white},>=angle 90,scale=2,text height=1.5ex, text depth=.25ex,row sep=4em, column sep=4em]
	\node (F110) at (-1,.5) {$F(z_0)\wedge F( z_1)$};
	\node (F120) at (1,.5) {$F(z_0)\wedge F( p_1)$};
	\node (F210) at (-1,-1.5) {$F(p_0)\wedge F( z_1)$};
	\node (F220) at (1,-1.5) {$\underline{F}$};
	\node at (.75,-1.35) {$\lrcorner$};
	\node (F11) at (0,-0.25) {$F(z_0\times z_1)$};
	\node (F12) at (2,-0.25) {$F(z_0\times p_1)$};
	\node (F21) at (0,-2.25) {$F(p_0\times z_1)$};
	\node (F22) at (2,-2.25) {$F(z)$};
	\node at (1.75,-2.1) {$\lrcorner$};
	\draw[>->,font=\scriptsize]
	(F210) edge node[descr]{$\!\!\!{i_{1}}_\ast\!\!\!$}node[descr,pos=.31] {} (F220)
	(F120) edge node[descr,pos=.6]{$\!\!\!{i_{0}}_\ast\!\!\!$}node[descr,pos=.355,sloped] {} (F220)
	(F11) edge (F12)
	(F21) edge node[descr]{$\!\!\!{i'_{1}}_\ast\!\!\!$}(F22)
	(F11) edge (F21)
	(F12) edge node[descr,pos=.6]{$\!\!\!{i'_{0}}_\ast\!\!\!$}(F22)
	(F110) edge node[descr]{$\!\!\!{i_{1}}_\ast\!\!\!$}(F120)
	(F110) edge node[descr,pos=.6]{$\!\!\!{i_{0}}_\ast\!\!\!$} (F210)
	;
	\draw[->,font=\scriptsize]
	(F110) edge node[descr]{$\!\!\!\phi_{z_0,z_1}\!\!\!$}(F11)
	(F210) edge node[descr]{$\!\!\!\phi_{p_0,z_1}\!\!\!$} (F21)
	(F120) edge node[descr]{$\!\!\!\phi_{z_0,p_1}\!\!\!$}(F12)
	;
	\draw[>->,dotted,font=\scriptsize]
	(F220) edge node[descr]{$\!\!\!\lambda'\!\!\!$} (F22)
	;
	\end{tikzpicture}
	\end{equation}
	Let $p$ denote the proper $S$-scheme $p_0\times p_1$, denote the closed immersion $z\cim p$ by $i$, and denote its complementary open immersion $x_0\times x_1\oim p$ by $j$. The morphism $i_\ast$ is a cofibration in $\scC$, and one has $\phi_{p_0,p_1}\circ \lambda=i_\ast\circ \lambda'$. Thus, there exists the solid commutative diagram
	\begin{equation}\label{Monodial:Last:Graph}
	\begin{tikzpicture}[descr/.style={fill=white},>=angle 90,scale=2,text height=1.5ex, text depth=.25ex,row sep=4em, column sep=4em]
	\node (F11) at (0,0) {$\underline{F}$};
	\node (F12) at (2.5,0) {$F(p_0)\wedge F(p_1)$};
	\node (F13) at (5,0) {$F^\ps({x}_0)\wedge F^\ps({x}_1)$};
	\node (F21) at (0,-1) {$F(z)$};
	\node (F22) at (2.5,-1) {$F(p)$};
	\node (F31) at (0,-2) {$F^\ps(z)$};
	\node (F32) at (2.5,-2) {$F^\ps(p)$};
	\node (F33) at (5,-2) {$F^\ps(x_0\times x_1)$};
	\draw[>->,font=\scriptsize]
	(F11) edge node[above] {$\lambda$}(F12)
	(F21) edge node[below] {$i_\ast$} (F22)
	(F31) edge node[below] {$i_\shriek$} (F32)
	;
	\draw[->>,font=\scriptsize]
	(F12) edge node[above] {$\pi$}(F13)
	(F32) edge node[below] {$j^!$}(F33)
	;
	\draw[->,font=\scriptsize]
	(F12) edge node[right]{$\phi_{p_0,p_1}$} (F22)
	(F11) edge node[left] {$\lambda'$} (F21)
	(F21) edge node[left] {$\varphi_z$} (F31)
	(F22) edge node[right] {$\varphi_{p}$} (F32)
	;
	\draw[->,dashed,font=\scriptsize]
	(F13) edge node[right] {$\phi^\ps_{{x}_0,{x}_1}$}(F33)
	;
	\end{tikzpicture}
	\end{equation}	
	in $\scC$, in which the horizontal sequences are {cofibre} sequences, and vertical morphisms are {weak equivalences}. Thus, the universal property of cokernels implies the existence of a unique morphism $\phi^\ps_{{x}_0,{x}_1}:F^\ps(x_0)\wedge F^\ps(x_1)\to F^\ps(x_0\times x_1)$ that makes the diagram commute; which is a weak equivalence in $\scC$, by the axiom \ref{itm:W3}. The uniqueness of the morphism $\phi^\ps_{{x}_0,{x}_1}$ implies the existence of a natural transformation ${\phi^\ps:F^\ps\wedge F^\ps\Rightarrow F^\ps(\times )}$ with a component $\phi^\ps_{{x}_0,{x}_1}$ for every pair of $S$-schemes $x_0$ and $x_1$. Also, a {diagram chase} of the associativity hexagons and unitality squares shows that $F^\ps$ is weak monoidal with the coherence {natural morphism} ${\phi^\ps}$. Moreover, when $\phi$ is a natural isomorphism, so is $\phi^\ps$.
\end{proof}

\subsubsection{Motivic Measures}\label{Motivic:Measure:PS}
Theorem \ref{Th:Induced:Measure} collects the main statements in \S.\ref{Extension:PS}, which allows one to associate motivic measures to \cdp-functors from proper $S$-schemes to \waldhausen categories.

\Th{}{Assume that $F:(\PropSch/S,\times,\id_{\!_S})\to (\scC,\wedge,\1)$ is a weak monoidal \cdp-functor to a {symmetric monoidal \waldhausen category}. Then, there exists a functor
	\[
	F^\ps:({\sftSch/S}^\prop_\open,\times,\id_{\!_S}) \to (\scC,\wedge,\1),
	\]
	where ${\sftSch/S}^\prop_\open$ is the category whose objects are $S$-schemes and whose morphisms are finite compositions of proper morphisms and formal inverses of open immersions, defined on proper morphisms in \eqref{PS:Proper:Pushforward} and on open immersions in \eqref{PS:Open:Pullback}, such that
	\begin{itemize}
		\item there exists a natural isomorphism $\varphi:F\stackrel{\sim}{\Rightarrow} F^\ps_{|_{\PropSch/S}}$;
		\item $F^\ps$ satisfies the excision property, \ie for every closed immersion $i:v\cim x$ of $S$-schemes with complementary open immersion $j:u\oim x$, the sequence 
		\[
		F^\ps(v)\stackrel{i_\shriek}{\longrightarrow} F^\ps(x) \stackrel{j^!}{\longrightarrow} F^\ps {(u)}
		\]
		is a cofibre sequence in $\scC$; and
		\item $F^\ps$ is weak monoidal, \ie $F^\ps$ is lax monoidal whose coherence morphisms
		\[
		\phi^\ps_{x,y}:F^\ps(x)\wedge F^\ps(y) \to F^\ps(x \times y)\andd \phi^\ps_{\!_S}:\1\to F(\id_{\!_S})
		\]
		are weak equivalences in $\scC$, for every $S$-schemes $x$ and $y$.
	\end{itemize}
	Therefore, there exists a motivic measure
	\[
	\mu_{\!_F}:\Kg_0(\sftSch/S)\to \Kg_0(\scC),
	\]
	that sends the class of a proper $S$-scheme $p$ to the class of $F(p)$.
}{Induced:Measure}
\Eg{}{Suppose that $\fk$ is a field. Then, a closed immersion $i:\Spec \fk \cim \PP^1_\fk$ is an initial object in the category of compactifications $\Comp_{\!_\fk}(\AA^1_\fk)$, and hence $F^\ps(\AA^1_\fk)=\coker i_\ast$.
}{Affine:Line}

\Rem{}{\zakharevich introduced the notion of an assembler in \cite{Zak:17:KTA}, using which she associated to the category of $\fk$-varieties, for a field $\fk$, a spectrum whose path components group is isomorphic to the \grothendieck group $\Kg_0(\sftSch/\fk)$. Then, Campbell provided an $E_\infty$-ring spectrum $\Kg({\Var/\fk})$ whose ring of path components is isomorphic to the \grothendieck ring $\Kg_0(\sftSch/\fk)$, which is conjectured to be equivalent to \zakharevich's spectrum, see \cite{Cam:17:KTSV}. {Lemma \ref{Lem:Base:Change}} and Proposition \ref{Prop:Excision} imply that a \cdp-functor $F:\PropSch/S\to \scC$ to a Waldhausen category defines a map of spectra $\Kg(F):\Kg({\Var/\fk})\to \Kg(\scC)$ that sends a point in the class $[P]\in \Kg_0({\Var/\fk})$ to a point in the class $[F(P)]\in \Kg_0(\scC)$, for every proper $\fk$-scheme $P$, see \cite[Def.5.2 and Prop.5.3]{Cam:17:KTSV}. 
}{}

\subsection{Functors Compactification}\label{Functors:Compactification}
Given a functor $F:\PropSch/S\to \scC$ to a \waldhausen category that is not a \cdp-functor, one may would like to `universally' associate to $F$ a \cdp-functor, and hence define an associated motivic measure that is closely related to $F$. Recall that the properties \ref{itm:PS:Initial} and \ref{itm:PS:Proper} imply that \cdp-functors are \cdp-cosheaves on $\PropSch/S$, as the \cdp-topology is generated by \cdp-squares, see \S.\ref{SubSec:GTAG}. Hence, a natural choice of such an association is the \cdp-cosheafification, when it exists, which is the dual of the \cdp-sheafification, see \cite{Pra:16:C}. We will restrain ourself from discussing the general process here, and only focus on aspects relevant to \S.\ref{Compactified:Yoneda}.

\medskip

Assume that $F$ satisfies \ref{itm:PS:Closed}, \ie it sends closed immersions of proper $S$-schemes to cofibrations in $\scC$. Then, for a \cdp-square
\[
\begin{tikzpicture}[descr/.style={fill=white},ampersand replacement=\&]
\node (m-1-1) at (0,0) {$z$};
\node (m-1-2) at (2,0) {$p$};
\node (m-2-1) at (0,-2) {$w$};
\node (m-2-2) at (2,-2) {$q$};
\node at (.25,-.4) {$\ulcorner$};
\path[->,font=\scriptsize]
(m-1-1) edge node[left]{$\underline{f}$} (m-2-1)
(m-1-2) edge node[right]{$f$} (m-2-2);
\path[right hook->,font=\scriptsize]
(m-2-1) edge node{$\diagup$}node[below]{$i$} (m-2-2)
(m-1-1) edge node{$\diagup$}node[above]{$\underline{i}$} (m-1-2);	
\end{tikzpicture}
\]
of proper $S$-schemes, the morphism $\underline{i}_\ast$ is a cofibration, and the pushout of the span of $\underline{i}_\ast$ and $\underline{f}_\ast$ exists in $\scC$. We denote the canonical morphism $F(w)\coprod_{F(z)}F(p)\to F(q)$ in $\scC$, induced by the universal property of pushouts, by $\alpha_{i,f}$, and consider the set of morphisms
\[
\overline{\Lambda}\coloneqq
\bigg\lbrace \alpha_{i,f}: F(w)\coprod_{F(z)}F(p)\to F(q)\mid (i,f)\in \Lambda\bigg\rbrace \bigcup \bigg\lbrace \0\to F(\emptyset)\bigg\rbrace
\]
in $\scC$, where $\Lambda$ is the set of all \cdp-squares of proper $S$-schemes. If there exists an exact functor $\scC\to \scC'$ of \waldhausen categories that sends all morphisms in $\overline{\Lambda}$ to isomorphisms in $\scC'$, the composition $F':\PropSch/S\to \scC\to \scC'$ satisfies the properties \ref{itm:PS:Closed}-\ref{itm:PS:Proper}, \ie it is a \cdp-functor. When the functor $\scC\to \scC'$ is a localisation with respect to $\overline{\Lambda}$, we say that the induced functor $F':\PropSch/S\to \scC'$ is a \indx{\cdp-compactification} of $F$. 

\medskip

When $\scC$ is a symmetric monoidal Waldhausen category and $F$ is only lax monoidal, one may seek a localisation for which the composite functor is also weak monoidal.

\medskip

In the next section, we apply this argument to the pointed \yoneda embedding.

\section{The \texorpdfstring{\cdp-}{cdp-}Waldhausen \texorpdfstring{$\Kg$-}{K-}Theory of Noetherian Schemes}\label{Compactified:Yoneda}

We utilise the \cdp-topology to construct a monoidal proper-fibred Waldhausen category $\scC^\omega_\tau$ over Noetherian schemes of finite Krull dimensions, in \S.\ref{Fibred:Noetherian}, for a topology $\tau$ that is finer than the \cdp-topology. For every Noetherian scheme $S$ of finite Krull dimension, there exists a canonical {\cdp-functor $\underline{\yon}^\tau:\PropSch/S\to \scC^\omega_\tau(S)$}, given by the $\tau$-cosheafification of the pointed \yoneda embedding, as in \eqref{cdp:Compactifiable:tau:Omega}. This functor induces, in \eqref{K:Schemes}, a surjective motivic measure
\[
\mu_{\!_\tau}:\Kg_0(\sftSch/S)\to \Kg_0\big( \scC^\omega_\tau(S)\big).
\]
We propose the $\Kg$-theory commutative ring spectrum $\Kg\big( \scC^\omega_\cdp(S)\big)$ as an extension of (modified) Grothendieck ring of $S$-schemes, see Corollary \ref{Cor:Factorisation:Modified} and Conjecture \ref{Conj:Universal:Measure:Modified}. However, we leave investigating the universality of the \cdp-functor $\underline{\yon}^\tau$ and the validity of \ref{Conj:Universal:Measure:Modified} for a future work.

\medskip

For an essentially small category $\scC$, the \yoneda embedding into the category of presheaves on $\scC$ gives a free cocompletion of $\scC$. Whereas, for a \grothendieck topology $\tau$ on $\scC$, the $\tau$-cosheafification of the \yoneda embedding gives a cocompletion of $\scC$ in the category of $\tau$-sheaves $\Shv_\tau(\scC)$, with the relations imposed by declaring $\tau$-coving sieves to be colimit cocones, see Remark \ref{Rem:Sheafification:Cocones}. The category of pointed $\tau$-sheaves admits a symmetric monoidal \waldhausen structure, whose cofibrations are monomorphisms, weak equivalences are isomorphisms, and monoidal product is given by the smash product, recall Example \ref{Ex:Waldhausen:Sheaves:Monoidal}. {In particular, when $\scC$ is the category of proper $S$-schemes and $\tau$ is a topology on $\PropSch/S$, it is interesting to consider when the $\tau$-cosheafification of the pointed \yoneda embedding is a \cdp-functor, and to use such a \cdp-functor, if it exists, to better understand the (modified) \grothendieck ring $\Kg_0(\sftSch/S)$, and probably its higher $\Kg$-theory.}

\medskip

For the rest of this section, let $\tau$ be an additively-saturated pretopology on the category Noetherian schemes of finite Krull dimensions that is finer than the \cdp-pretopology and coarser than the {proper} pretopology, \cf Remark \ref{Rem:Topology:Restrict:Attention}. Also, let $S$ be a \noetherian scheme in $\fdNoe$.

\medskip

Recall that the category of proper $S$-schemes is essentially small, and the forgetful functor $\PSh_\bullet(\PropSch/S)\to \PSh(\PropSch/S)$, that forgets the base point, admits a faithful (but not full) left adjoint $\PSh(\PropSch/S)\to \PSh_\bullet(\PropSch/S)$, given by adjoining a disjoint base point, \ie $\scX_+=(\scX \coprod \ast, \ast)$, see \cite[p.4]{Hov:99:MC}. Let $\yon_{\!_{-,+}}$ denote the composite functor $-_{\!_+}\circ \yon:\PropSch/S\to \PSh_\bullet(\PropSch/S)$. The gluing of a pair of closed subschemes of a proper $S$-scheme, along their scheme-theoretic intersection, defines a pushout square in $\PropSch/S$, which is a \cdp-square, {see the proof of Proposition \ref{Prop:PS:Properties}}. Then, the functor $\yon_{\!_{-,+}}$ is not a \cdp-functor, as it forgets all colimits. However, for every closed immersion $i:z\cim p$ in $\PropSch/S$, the morphism $\yon_{\!_{i,+}}$ is a monomorphism. Following the argument in \S.\ref{Functors:Compactification}, we may consider a localisation of the category $\PSh_\bullet(\PropSch/S)$ with respect to the set of morphisms
\[
\overline{\Lambda}\coloneqq
\bigg\lbrace \alpha_{i,f}: \yon_{\!_{w,+}}\coprod_{\yon_{\!_{z,+}}}\yon_{\!_{p,+}}\to \yon_{p,+} \mid (i,f)\in \Lambda\bigg\rbrace \bigcup \bigg\lbrace \0\to \yon_{\!_{\emptyset,+}}\bigg\rbrace,
\]
where $\Lambda$ is the set of all \cdp-squares in $\PropSch/S$. 

\medskip
The \cdp-sheafification functor $-^{\!^{\ash_\cdp}}:\PSh_\bullet(\PropSch/S)\to \Shv_{\!_{\bullet,\cdp}}(\PropSch/S)$ provides such a localisation, see Definition \ref{Def:cd:topology}. That is,
\begin{itemize}
	\item[\ref{itm:PS:Closed}] the \cdp-sheafification functor preserve monomorphisms;
	\item[\ref{itm:PS:Initial}] the \cdp-sheafification of $\yon_{\!_{\emptyset,+}}$ is isomorphic to $\0$, as $\yon_{\!_{\emptyset,+}}^{\!^\cdp}(\emptyset)=\ast$; and
	\item[\ref{itm:PS:Proper}] the functor $\yon^{\!^{\ash_\cdp}}$, \ie the composition of the \cdp-sheafification functor with the the \yoneda embedding, sends \cdp-squares to pushout squares, by Proposition \ref{Prop:cd:Cosheafification}; also, the left adjoint functor $-{_+}^{\!\!\!\!^{\ash_\cdp}}$ preserves colimits.	
\end{itemize}
Therefore, the functor 
\begin{equation}\label{cdp:Compactifiable}
\underline{\yon}\coloneqq -^{\ash_\cdp}\circ \yon_{\!_{-,+}}:\PropSch/S \to \Shv_{\!_{\bullet,\cdp}}(\PropSch/S)
\end{equation}
is a \cdp-functor. Moreover, Remark \ref{Rem:Sheafification:Cocones} shows that the \cdp-topology is the coarsest topology $\tau$ on $\PropSch/S$ for which the composite functor $-^{\ash_\tau}\circ \yon_{\!_{-,+}}$ is a \cdp-functor. The $\tau$-sheafification functor preserves monomorphisms and colimits, and it factorises through the \cdp-sheafification functor, as $\tau$ is finer than the \cdp-pretopology. Hence, the functor
\begin{equation}\label{cdp:Compactifiable:tau}
\underline{\yon}^\tau\coloneqq -^{\ash_\tau}\circ \yon_{\!_{-,+}}:\PropSch/S \to \Shv_{\!_{\bullet,\tau}}(\PropSch/S)
\end{equation}
is a \cdp-functor. To avoid bulky notations, we let $\scC_\tau(S)\coloneqq \Shv_{\!_{\bullet,\tau}}(\PropSch/S)$. When $\tau$ is the \cdp-pretopology, denote $\scC_\tau(S)$, $\underline{\yon}^\tau$, and $\yon^{\ps_\tau}$, by $\scC(S)$, $\underline{\yon}$, and $\yon^\ps$, respectively.

\medskip
The functor $\yon_{\!_{-,+}}$ is strong monoidal, with respect to the Cartesian product of proper $S$-schemes and the smash product of pointed presheaves. Since the $\tau$-sheafification functor is a left exact reflector, it preserves the smash product, as the smash product of pointed (pre)sheaves only involves finite limits and colimits, recall Example \ref{Ex:Waldhausen:Sheaves:Monoidal}. Thus, the functor $\underline{\yon}^\tau$ is also strong monoidal.

\medskip
Since the \waldhausen category $\scC_\tau(S)$ is cocomplete, its $\Kg$-theory is connected, \ie it has a trivial path components group, recall Example \ref{Ex:Eilenberg:Swindle}. To establish a non-trivial motivic measure, we need to consider a \waldhausen subcategory in $\scC_\tau(S)$, with a non-connected $\Kg$-theory, which contains the essential image of $\underline{\yon}^\tau$. We construct this subcategory by induction.
\begin{itemize}
	\item Let $\scC^0_\tau(S)$ be the full subcategory in $\scC_\tau(S)$ in which $\scX\in \scC^0_\tau(S)$ \iff $\scX\cong \underline{\yon}^\tau_{\!_p}$ for a proper $S$-scheme $p$.
	\item For an integer $n\geq 1$, let $\scC^n_\tau(S)$ be the full subcategory in $\scC_\tau(S)$ in which $\scX\in \scC^n_\tau(S)$ \iff there exists a pushout square
	\begin{equation}\label{Inductive:C:n}
	\begin{tikzpicture}[descr/.style={fill=white},ampersand replacement=\&]
	\node (m-1-1) at (0,0) {$\scY'$};
	\node (m-1-2) at (2,0) {$\scY$};
	\node (m-2-1) at (0,-2) {$\scX'$};
	\node (m-2-2) at (2,-2) {$\scX$};
	\node at (1.75,-1.6) {$\lrcorner$};
	\path[>->,font=\scriptsize]
	(m-2-1) edge (m-2-2)
	(m-1-1) edge node	[above]	{$\iota$} (m-1-2);	
	\path[->,font=\scriptsize]	
	(m-1-1) edge (m-2-1)
	(m-1-2) edge (m-2-2);
	\end{tikzpicture}
	\end{equation}
	in $\scC_\tau(S)$, in which $\scX', \scY,$ and $\scY'$ belong to $\scC^{n-1}_\tau(S)$, and $\iota$ is a monomorphism of pointed $\tau$-sheaves.	
\end{itemize}
One has $\scC^n_\tau(S)\subset \scC^{n+1}_\tau(S)$ for every $n\in\NN$. Then, the full subcategory 
\begin{equation}\label{Waldhausen:Omega}
\scC^\omega_\tau(S)\coloneqq \bigcup_{n\in \NN} \scC^n_\tau(S)
\end{equation}
in $\scC_\tau(S)$ admits a Waldhausen structure, whose cofibrations (\resp. weak equivalences) are morphisms in $\scC^\omega_\tau(S)$ that are cofibrations (\resp. weak equivalences) in $\scC_\tau(S)$, \ie monomorphisms (\resp. isomorphisms). 

\medskip

The inclusion functor $\scC^\omega_\tau(S)\incl \scC_\tau(S)$ is an exact functor of Waldhausen categories, \ie $\scC^\omega_\tau(S)$ is a \waldhausen subcategory in $\scC_\tau(S)$. In fact, $\scC^\omega_\tau(S)$ is a the smallest full \waldhausen subcategory in $\scC_\tau(S)$ that contains the essential image of $\underline{\yon}^\tau$. We abuse notation and we use $\underline{\yon}^\tau$ to also denote the unique \cdp-functor 
\begin{equation}\label{cdp:Compactifiable:tau:Omega}
\PropSch/S \to \scC^\omega_\tau(S)
\end{equation} 
that factorises $\underline{\yon}^\tau$. {Since the category $\PropSch/S$ is essentially small, one can use induction to show that $\scC^\omega_\tau(S)$ is also essentially small.}

\Lem{}{The symmetric monoidal structure on $\scC_\tau(S)$, given by the smash product, restricts to $\scC^\omega_\tau(S)$, making $\scC^\omega_\tau(S)$ into a symmetric monoidal \waldhausen category.
}{Symmetric:Monoidal}
\begin{proof}
	Since $\scC^\omega_\tau(S)$ is a full \waldhausen subcategory in $\scC_\tau(S)$ and the unit $\1_{\wedge}=\underline{\yon}^\tau_{\!_{\id_{\!_S}}}$ belongs to $\scC^\omega_\tau(S)$, it is suffices to show that the smash product restricts to $\scC^\omega_\tau(S)$. 
	
	\medskip
	
	For every $n\in \NN$, we use induction to show that the smash product restricts to a functor
	\[
	\wedge:\scC^n_\tau(S)\times \scC^n_\tau(S)\to \scC^{2n}_\tau(S).
	\]
	\begin{itemize}
		\item Let $\scX_0$ and $\scX_1$ belong to $\scC^0_\tau(S)$, \ie there exists a proper $S$-scheme $p_k$ for which $\scX_k\cong \underline{\yon}^\tau_{\!_{p_k}}$, for $k=0,1$. Since $\underline{\yon}^\tau$ is a strong monoidal functor, one has
		\[
		\scX_0\wedge\scX_1\ \cong \ \underline{\yon}^\tau_{\!_{p_0}}\!\!\wedge\ \underline{\yon}^\tau_{\!_{p_1}}\!\!\cong\ \underline{\yon}^\tau_{\!_{p_0\times p_1}} \in \scC^0_\tau(S).
		\]
		\item For an integer $n\geq 1$, assume that the smash product restricts to a functor
		\[
		\wedge:\scC^{n-1}_\tau(S)\times \scC^{n-1}_\tau(S)\to \scC^{2n-2}_\tau(S),
		\]
		and let $\scX_0$ and $\scX_1$ belong to $\scC^n_\tau(S)$, \ie there exists a pushout square
		\[
		\begin{tikzpicture}[descr/.style={fill=white},ampersand replacement=\&]
		\node (m-1-1) at (0,0) {$\scY'_k$};
		\node (m-1-2) at (2,0) {$\scY_k$};
		\node (m-2-1) at (0,-2) {$\scX'_k$};
		\node (m-2-2) at (2,-2) {$\scX_k$};
		\node at (1.75,-1.6) {$\lrcorner$};
		\path[>->,font=\scriptsize]
		(m-2-1) edge (m-2-2)
		(m-1-1) edge node	[above]	{$\iota_k$} (m-1-2);	
		\path[->,font=\scriptsize]	
		(m-1-1) edge (m-2-1)
		(m-1-2) edge (m-2-2);
		\end{tikzpicture}
		\]
		in $\scC_\tau(S)$, in which $\scX'_k, \scY_k,$ and $\scY'_k$ belong to $\scC^{n-1}_\tau(S)$, and $\iota_k$ is a monomorphism, for $k=0,1$. Since $\scC_\tau(S)$ is a symmetric monoidal \waldhausen category, the smash product with any object in $\scC_\tau(S)$ preserves such a pushout square, and hence there exists a pushout diagram
		\[
		\begin{tikzpicture}[descr/.style={fill=white},ampersand replacement=\&]
		\node (m-1-1) at (0,0) {$\scU\wedge\scY'_k$};
		\node (m-1-2) at (3,0) {$\scU\wedge\scY_k$};
		\node (m-2-1) at (0,-2) {$\scU\wedge\scX'_k$};
		\node (m-2-2) at (3,-2) {$\scU\wedge\scX_k$};
		\node at (2.75,-1.6) {$\lrcorner$};
		\path[>->,font=\scriptsize]
		(m-2-1) edge (m-2-2)
		(m-1-1) edge node	[above]	{$\id_{\scX} \wedge \iota_k$} (m-1-2);	
		\path[->,font=\scriptsize]	
		(m-1-1) edge (m-2-1)
		(m-1-2) edge (m-2-2);
		\end{tikzpicture}
		\]
		in $\scC_\tau(S)$, for every $\scU\in \scC_\tau(S)$. In particular, for $k=0$ and $\scU\in \{\scX'_1,\scY'_1,\scY_1\}$, one finds that $ \scX_0\wedge \scU\cong \scU\wedge \scX_0$ belongs to $\scC^{2n-1}_\tau(S)$. Similarly, for $k=1$ and $\scU=\scX_0$, one finds that $\scX_0\wedge \scX_1$ belongs to $\scC^{2n}_\tau(S)$.
	\end{itemize}
	Therefore, the smash product restricts to $\scC^\omega_\tau(S)$.
\end{proof}
The functor $\underline{\yon}^\tau:\PropSch/S \to \scC^\omega_\tau(S)$ is a strong monoidal \cdp-functor. Then, by Theorem \ref{Th:Induced:Measure}, there exists a strong monoidal functor
\begin{equation}\label{K:Schemes:Functor}
\yon^{\ps_\tau}:\sftSch^\prop_\open/S\to \scC^\omega_\tau(S)
\end{equation}
that satisfies the excision property, and coincides with $\underline{\yon}^\tau$ on proper $S$-schemes. Hence, it induces a motivic measure
\begin{equation}\label{K:Schemes}
\mu_{\!_\tau}:\Kg_0(\sftSch/S)\to \Kg_0\big( \scC^\omega_\tau(S)\big),
\end{equation}
which sends the class of a proper $S$-scheme $p$ to the class of $\underline{\yon}^\tau_{\!_{p}}$. Moreover, for a field $\fk$, the functor $\yon^{\ps_\tau}$ induces a map of spectra $\Kg(\yon^{\ps_\tau}):\Kg({\Var/\fk})\to \Kg\big(\scC^\omega_\tau(\fk)\big)$, from the spectrum $\Kg({\Var/\fk})$ defined in \cite{Cam:17:KTSV}.

\Eg{}{
Recall Example \ref{Ex:Affine:Line}, and consider a closed immersion $i:\Spec \fk \cim \PP^1_\fk$ with complementary open immersion $j:\AA^1_\fk \oim \PP^1_\fk$. Then,
\[
\yon^{\ps_\tau}_{\AA^1_\fk}=\coker \big(\underline{\yon}^\tau_{\!_{\Spec \fk}}\longrightarrow\ \underline{\yon}^\tau_{\!_{\PP^1_\fk}}\big)=({\yon}_{\!_{\PP^1_\fk}},\infty)^{\!^{\ash_\tau}},
\]
where $\infty$ denotes the unique $\fk$-rational point in the $\fk$-scheme $\big(\PP^1_\fk{\setminus} j(\AA^1_\fk)\big)_{\!_\red}$.
}{}

\Lem{}{The motivic measure $\mu_{\!_\tau}$ is surjective.}{Surjective:Measure}
\begin{proof}The statement of the lemma follows by induction.
\end{proof}
The connection between the motivic measure $\mu_{\!_\tau}$ and the modified Grothendieck ring of $S$-schemes might be illustrated through Corollary \ref{Cor:Factorisation:Modified}, which is based on the following proposition, and Conjecture \ref{Conj:Universal:Measure:Modified}.

\Prop{}{Let $f:x\to y$ be a morphism of $S$-schemes. Then, the following are equivalent
	\begin{enumerate}
		\item\label{itm:Factorisation:Modified:iso} $f$ is a proper morphism and $f_\shriek:\yon^{\ps_\tau}_{\!_x}\to \yon^{\ps_\tau}_{\!_y}$ is an isomorphism in $\scC_\tau^\omega(S)$; and
		\item\label{itm:Factorisation:Modified:uh} $f$ is a universal homeomorphism, \cf Proposition \ref{Prop:Topology:Isomorphism:Proper:If:Only:If}.
	\end{enumerate}
}{Factorisation:Modified}
\begin{proof} Recall our conventions in \S.\ref{Notations}, which require $S$-schemes to be of finite type. Hence, the morphism $f$ is a universal homeomorphism \iff it is finite, universally injective, and surjective, by \cite[Prop.2.4.5]{EGA4.2}. Let $l:w\cim q$ be a compactification of $y$. Since $f$ is proper under either of the conditions \eqref{itm:Factorisation:Modified:iso} or \eqref{itm:Factorisation:Modified:uh}, the category $\Comp_{\!_S}(f,l)$ is nonempty in either case, by Proposition \ref{Prop:Compactification:Proper:Morphism}. Suppose that $i:z\cim p$ is a compactification in $\Comp_{\!_S}(f,l)$, and let $g:i\to l$ be a morphism of compactifications such that $f$ is a base change in $\sftSch/S$ of $g$ along $j_{\!_l}$. For simplicity and without loss of generality, we may choose $i$ that fits into a diagram 
	\begin{equation}\label{Eq:Universal:Homeomorphism:cdf}
	\begin{tikzpicture}[descr/.style={fill=white}]
	\node (W0) at (0,0) {$z$};
	\node (X0) at (4,0) {$x$};
	\node (XBar0) at (2,0) {$p$};
	\node (W1) at (0,-2) {$w$};
	\node (X1) at (4,-2) {$y$};
	\node (XBar1) at (2,-2) {$q$};
	\node at (.25,-.4) {$\ulcorner$};
	\node at (3.75,-.4) {$\urcorner$};
	\path[->,font=\scriptsize]
	(XBar0) edge node[left]{$g$} (XBar1)
	(X0) edge node[right]{$f$} (X1)
	(W0) edge node[left]{$\underline{g}$} (W1)
	;
	\path[left hook->,font=\scriptsize]
	(X0) edge node[descr]{$\!\!\!\circ\!\!\!$}node[above]{$j_{\!_i}$} (XBar0)
	(X1) edge node[descr]{$\!\!\!\circ\!\!\!$}node[below]{$j_{\!_l}$} (XBar1)
	;
	\path[right hook->,font=\scriptsize]
	(W1) edge node{$\diagup$}node[below]{$l$} (XBar1)
	(W0) edge node{$\diagup$}node[above]{$i$} (XBar0)
	;
	\end{tikzpicture}
	\end{equation}
	of Cartesian squares in $\sftSch/S$. The morphism $f_\shriek:\yon^{\ps_\tau}_{\!_x}\to \yon^{\ps_\tau}_{\!_y}$ is an isomorphism in $\scC^\omega_\tau(S)$ \iff the square
	\begin{equation}\label{Pushout:Square:Modified:Ring}
	\begin{tikzpicture}[descr/.style={fill=white}]
	\node (W0) at (0,0) {$\underline{\yon}^\tau_{\!_z}$};
	\node (XBar0) at (2,0) {$\underline{\yon}^\tau_{\!_p}$};
	\node (W1) at (0,-2) {$\underline{\yon}^\tau_{\!_w}$};
	\node (XBar1) at (2,-2) {$\underline{\yon}^\tau_{\!_q}$};
	\path[->,font=\scriptsize]
	(XBar0) edge node[right]{$g_\ast$} (XBar1)
	(W0) edge node[left]{$\underline{g}_\ast$} (W1)
	;
	\path[>->,font=\scriptsize]
	(W1) edge node[below]{$l_\ast$} (XBar1)
	(W0) edge node[above]{$i_\ast$} (XBar0)
	;
	\end{tikzpicture}
	\end{equation}
	is a pushout square in $\scC^\omega_\tau(S)$. Since the $\tau$-sheafification functor preserves colimits, that is equivalent to the canonical morphism $\Theta:\yon_{\!_{w,+}}\coprod_{\yon_{\!_{z,+}}} \yon_{\!_{p,+}}\to \yon_{\!_{q,+}}$ of pointed presheaves being a $\tau$-local isomorphism.
	
	\medskip
	
	\begin{itemize}
		\item[{$(\ref{itm:Factorisation:Modified:uh}\Rightarrow \ref{itm:Factorisation:Modified:iso})$}] Assume that $f$ is a universal homeomorphism. Then, the set $\{l:w\cim q,g:p\to q\}$ is a \cdp-covering family of $q$ in $\PropSch/S$. Indeed, let $\fk$ be a field, and let $b:\Spec \fk \to Q$ be a morphism of schemes in $\fdNoe$, where $Q$ is the underlying scheme of $q$. Either $b$ lifts along $l$ or $j_{\!_l}$. When $b$ lifts along $j_{\!_l}$, there exists a morphism $a:\Spec \fk\to Y$ in $\fdNoe$ for which $b=j_{\!_l}\circ a$, where $Y$ is the underlying scheme of $y$. Consider the Cartesian square
		\[
		\begin{tikzpicture}[descr/.style={fill=white},ampersand replacement=\&]
		\node (m-1-1) at (0,0) {$T$};
		\node (m-1-2) at (2,0) {$X$};
		\node (m-2-1) at (0,-2) {$\Spec \fk$};
		\node (m-2-2) at (2,-2) {$Y$};
		\node at (.25,-.4) {$\ulcorner$};
		\path[->,font=\scriptsize]
		(m-1-1) edge node[left]{$\underline{f}$} (m-2-1)
		(m-1-2) edge node[right]{$f$} (m-2-2)
		(m-2-1) edge node[below]{$a$} (m-2-2)
		(m-1-1) edge node[above]{$\underline{a}$} (m-1-2);	
		\end{tikzpicture}
		\]
		in the category $\fdNoe$. The morphism $f$ is a finite universal homeomorphism, and hence $T$ is a one-point scheme $\Spec R$ and $\underline{f}$ is induced by a finite ring homomorphism $\psi:\fk \incl R$, to a local ring $R$ of Krull dimension zero. Let $\frm$ be the maximal ideal of $R$, and let $\rf\coloneqq \sfrac{R}{\frm}$. Then, the induced homomorphism $\fk \incl \rf$ is a finite field extension. Assuming that $[\rf:\fk]\neq 1$, there exist distinct ring homeomorphisms $\rf\to \rf$ over $\fk$, which contradicts with $\underline{f}$ being universally injective. Thus, one has $[\rf:\fk]=1$, \ie the residue field of $T$ at its unique point is isomorphic to $\fk$. Therefore, $a$ lifts along $f$, which lifts $b$ along $g$, and hence $\{l,q\}$ is a \cdp-covering family.

		\medskip 
		
		Let $t:T\to S$ be a proper $S$-scheme, and let $b\in \yon_{\!_{q,+}}(t)$. Either $b$ is the base point of $\yon_{\!_{q,+}}(t)$, in which case $b$ belongs to the image of $\Theta_{\!_t}$, or $b$ is a morphism $t\to q$ of proper $S$-schemes. Assuming the latter case, consider the Cartesian squares
		\[
		\begin{tikzpicture}[descr/.style={fill=white}]
		\node (W0) at (0,0) {$t_{\!_l}$};
		\node (XBar0) at (2,0) {$t$};
		\node (W1) at (0,-2) {$w$};
		\node (XBar1) at (2,-2) {$q$};
		\node at (.25,-.4) {$\ulcorner$};
		\path[->,font=\scriptsize]
		(XBar0) edge node[right]{$b$} (XBar1)
		(W0) edge node[left]{$b_{\!_l}$} (W1)
		;
		\path[right hook->,font=\scriptsize]
		(W1) edge node{$\diagup$}node[below]{$l$} (XBar1)
		(W0) edge node{$\diagup$}node[above]{$l_{\!_b}$} (XBar0)
		;
		\end{tikzpicture}\qquad
		\begin{tikzpicture}[descr/.style={fill=white}]
		\node (W0) at (0,0) {$t_{\!_g}$};
		\node (XBar0) at (2,0) {$t$};
		\node (W1) at (0,-2) {$p$};
		\node (XBar1) at (2,-2) {$q$};
		\node at (.25,-.4) {$\ulcorner$};
		\path[->,font=\scriptsize]
		(XBar0) edge node[right]{$b$} (XBar1)
		(W0) edge node[left]{$b_{\!_g}$} (W1)
		;
		\path[->,font=\scriptsize]
		(W1) edge node[below]{$g$} (XBar1)
		(W0) edge node[above]{$g_{\!_b}$} (XBar0)
		;
		\end{tikzpicture}
		\]
		in $\PropSch/S$. The set $\{l_{\!_b},g_{\!_b}\}$ is a \cdp-covering family of $t$ in $\PropSch/S$. Let $[b_{\!_l}]$ and $[b_{\!_g}]$ denote the classes of $b_{\!_l}\in \yon_{\!_{w,+}}(t_{\!_l})$ and $b_{\!_g}\in\yon_{\!_{p,+}}(t_{\!_g})$ in $\big(\yon_{\!_{w,+}}\coprod_{\yon_{\!_{z,+}}} \yon_{\!_{p,+}}\big)(t_{\!_l})$ and $\big(\yon_{\!_{w,+}}\coprod_{\yon_{\!_{z,+}}} \yon_{\!_{p,+}}\big)(t_{\!_g})$, respectively,. Then, $\Theta_{\!_{t_{\!_l}}}([b_{\!_l}])=l_\ast(b_{\!_l})=l_{\!_b}^\ast (b)$	and	$\Theta_{\!_{t_{\!_g}}}([b_{\!_g}])=g_\ast(b_{\!_g})=g_{\!_b}^\ast (b)$. Since the pretopology $\tau$ is finer than the \cdp-pretopology, the morphism  $\Theta$ is a $\tau$-local epimorphism, by Lemma \ref{Lem:Local_Morphism}.
		
		\medskip
		
		On the other hand, suppose that $t:T\to S$ is a proper $S$-scheme, and let $a_0$ and $a_1$ be sections in $\big(\yon_{\!_{w,+}}\coprod_{\yon_{\!_{z,+}}} \yon_{\!_{p,+}}\big)(t)$ such that $\Theta_{\!_t}(a_0)=\Theta_{\!_t}(a_1)$. When either $a_0$ or $a_1$ coincides with the base point in $\big(\yon_{\!_{w,+}}\coprod_{\yon_{\!_{z,+}}} \yon_{\!_{p,+}}\big)(t)$, so does the other. Assuming that $a_0\neq\ast$ and $a_1\neq\ast$, we distinguish the three cases.
		\begin{enumerate}
			\item\label{itm:Closed:Immersion} There exist sections $a'_0$ and $a'_1$ in $\yon_{\!_{w,+}}(t)$ such that $a_0=[a'_0]$ and $a_1=[a'_1]$. Then, $l_\ast (a'_0)=\Theta_{\!_t}(a_0)=\Theta_{\!_t}(a_1)=l_\ast (a'_1)$. Since $a_0\neq\ast$ and $a_1\neq\ast$, the sections $a'_0$ and $a'_1$ are morphisms $t\to w$ of proper $S$-schemes for which $l\circ a'_0=l\circ a'_1$, which implies that $a'_0=a'_1$, as $l$ is a monomorphism of proper $S$-schemes. Hence, $a_0=a_1$.
			\item There exists a section $a'_0$ in $\yon_{\!_{w,+}}(t)$ and a section $a'_1$ in $\yon_{\!_{p,+}}(t)$ such that $a_0=[a'_0]$ and $a_1=[a'_1]$. Then, $a'_0$ (\resp. $a'_1$) is a morphism $t\to w$ (\resp. $t\to p$) of proper $S$-schemes, as $a_0\neq\ast$ and $a_1\neq\ast$. Then, $l\circ a'_0=\Theta_{\!_t}(a_0)=\Theta_{\!_t}(a_1)=g\circ  a'_1$, and hence there exists a morphism $a':t\to z$ of proper $S$-schemes such that $a'_0=\underline{g}_\ast(a')$ and $a'_1=i_\ast(a')$, as the diagram \eqref{Eq:Universal:Homeomorphism:cdf} consists of Cartesian squares. Thus, $a_0=[a'_0]=[a'_1]=a_1$. The same argument applies when there exists a section $a'_0$ in $\yon_{\!_{p,+}}(t)$ and a section $a'_1$ in $\yon_{\!_{w,+}}(t)$ such that $a_0=[a'_0]$ and $a_1=[a'_1]$.
			\item There exist sections $a'_0$ and $a'_1$ in $\yon_{\!_{p,+}}(t)$ such that $a_0=[a'_0]$ and $a_1=[a'_1]$. As $a_0\neq\ast$ and $a_1\neq\ast$, the sections $a'_0$ and $a'_1$ are morphisms $t\to p$ of proper $S$-schemes such that $g\circ a'_0=\Theta_{\!_t}(a_0)=\Theta_{\!_t}(a_1)=g\circ  a'_1 $. Consider the diagram
			\begin{equation}\label{uh:First:Graph}
			\begin{tikzpicture}[descr/.style={fill=white}]
			\node (W2) at (0,2) {$t_{\!_{k,l}}$};
			\node (X2) at (4,2) {$t_{\!_{k,j_{\!_{l}}}}$};
			\node (XBar2) at (2,2) {$t$};
			\node (W0) at (0,0) {$z$};
			\node (X0) at (4,0) {$x$};
			\node (XBar0) at (2,0) {$p$};
			\node (W1) at (0,-2) {$w$};
			\node (X1) at (4,-2) {$y$};
			\node (XBar1) at (2,-2) {$q$};
			\node at (.25,1.6) {$\ulcorner$};
			\node at (3.75,1.6) {$\urcorner$};
			\node at (.25,-.4) {$\ulcorner$};
			\node at (3.75,-.4) {$\urcorner$};
			\path[->,font=\scriptsize]
			(XBar2) edge node[left]{$a'_k$} (XBar0)
			(X2) edge node[right]{$a'_{\!_{k,j_{\!_l}}}$} (X0)
			(W2) edge node[left]{$a'_{\!_{k,l}}$} (W0)
			(XBar0) edge node[left]{$g$} (XBar1)
			(X0) edge node[right]{$f$} (X1)
			(W0) edge node[left]{$\underline{g}$} (W1)
			;
			\path[left hook->,font=\scriptsize]
			(X2) edge node[descr]{$\!\!\!\circ\!\!\!$}node[above]{$j_{\!_{k,l}}$} (XBar2)
			(X0) edge node[descr]{$\!\!\!\circ\!\!\!$}node[above]{$j_{\!_i}$} (XBar0)
			(X1) edge node[descr]{$\!\!\!\circ\!\!\!$}node[below]{$j_{\!_l}$} (XBar1)
			;
			\path[right hook->,font=\scriptsize]
			(W2) edge node{$\diagup$}node[above]{$l_k$} (XBar2)
			(W1) edge node{$\diagup$}node[below]{$l$} (XBar1)
			(W0) edge node{$\diagup$}node[above]{$i$} (XBar0)
			;
			\end{tikzpicture}
			\end{equation}		
			of Cartesian squares of $S$-schemes, for $k=1,0$. Since $g\circ a'_0=g\circ  a'_1$, there exist such Cartesian squares with 
			\[
			t_{\!_l}\coloneqq t_{\!_{0,l}}=t_{\!_{1,l}}\quad,\quad
			\underline{l}\coloneqq l'_0=l'_1\quad,\quad
			t_{\!_{j_{\!_l}}}\coloneqq t_{\!_{0,j_{\!_{l}}}}=t_{\!_{1,j_{\!_{l}}}}\quad\text{, and }\quad
			{j}_{\!_{\underline{l}}}\coloneqq j'_{\!_{0,l}}=j'_{\!_{1,l}}.
			\]
			
			\medskip
			Then, one has $f\circ a'_{\!_{0,j_{\!_l}}}=f\circ a'_{\!_{1,j_{\!_l}}}$, which implies the existence of a \cdp-cover $\sigma_{\!_{j_{\!_l}}}:t'_{\!_{j_{\!_l}}}\to t_{\!_{j_{\!_l}}}$ such that 
			\begin{equation}\label{Eq:Equal:Dominant}
			a'_{\!_{0,j_{\!_l}}}\circ \sigma_{\!_{j_{\!_l}}}=a'_{\!_{1,j_{\!_l}}}\circ \sigma_{\!_{j_{\!_l}}},
			\end{equation}
			as $f_\ast:\yon_{\!_x}\to \yon_{\!_y}$ is a \cdp-local monomorphism, by Proposition \ref{Prop:Voevodsky:96:HS:Prop.3.2.5}.
			
			\medskip
			Since $\sigma_{\!_{j_{\!_l}}}$ is a proper morphism, the category $\Comp_{\!_S}(\sigma_{\!_{j_{\!_l}}},\underline{l})$ is nonempty, and there exists a compactification $l':t'_{\!_l}\cim t'$  of $t'_{\!_{j_{\!_l}}}$ that fits into a diagram
			\[
			\begin{tikzpicture}[descr/.style={fill=white}]
			\node (W0) at (0,0) {$t'_{\!_l}$};
			\node (X0) at (4,0) {$t'_{\!_{j_{\!_l}}}$};
			\node (XBar0) at (2,0) {$t'$};
			\node (W1) at (0,-2) {$t_{\!_l}$};
			\node (X1) at (4,-2) {$t_{\!_{j_{\!_l}}}$};
			\node (XBar1) at (2,-2) {$t$};
			\node at (.25,-.4) {$\ulcorner$};
			\node at (3.75,-.4) {$\urcorner$};
			\path[->,font=\scriptsize]
			(XBar0) edge node[left]{$\sigma$} (XBar1)
			(X0) edge node[right]{$\sigma_{\!_{j_{\!_l}}}$} (X1)
			(W0) edge node[left]{$\sigma_{\!_{l}}$} (W1)
			;
			\path[left hook->,font=\scriptsize]
			(X0) edge node[descr]{$\!\!\!\circ\!\!\!$}node[above]{$j_{\!_{l'}}$} (XBar0)
			(X1) edge node[descr]{$\!\!\!\circ\!\!\!$}node[below]{${j}_{\!_{\underline{l}}}$} (XBar1)
			;
			\path[right hook->,font=\scriptsize]
			(W1) edge node{$\diagup$}node[below]{$\underline{l}$} (XBar1)
			(W0) edge node{$\diagup$}node[above]{$l'$} (XBar0)
			;
			\end{tikzpicture}
			\]
			of Cartesian squares of $S$-schemes. Since $\sigma_{\!_{j_{\!_l}}}$ is a \cdp-cover, the set $\{ \underline{l}:t_{\!_l}\cim t, \sigma:t'\to t\}$ is a \cdp-covering family of $t$ in $\PropSch/S$. Indeed, let $\fk$ be a field, and let $b:\Spec \fk \to T$ be a morphism of schemes in $\fdNoe$, then either $b$ lifts along $\underline{l}$ or ${j}_{\!_{\underline{l}}}$. In the latter case, the lift $\Spec \fk \to T_{\!_{{j}_{\!_{\underline{l}}}}}$ lifts along the \cdp-cover $\sigma_{\!_{j_{\!_l}}}$, which lifts $b$ along $\sigma$, where $T_{\!_{{j}_{\!_{\underline{l}}}}}$ is the underlying scheme of $t_{\!_{{j}_{\!_{\underline{l}}}}}$. Since $\underline{l}$ is a monomorphism of schemes and ${j}_{\!_{\underline{l}}}$ is a complementary open immersion to $\underline{l}$, one sees that the morphism $\sigma_{\!_l}$ is a \cdp-cover.
			
			\medskip
			
			Let $\{i_{\!_\alpha}:t_{\!_\alpha}\to t'\mid \alpha\in A\}$ be the \cdp-covering family of $t'$ in $\PropSch/S$ by its integral components, and consider the Cartesian squares \eqref{Universal:Homeomorphism} of $S$-schemes, for every $\alpha\in A$. Then, the set $\{i_{\!_{\alpha,l}}:t_{\!_{\alpha,l}}\to t'_{\!_l}\mid\alpha\in A\}$ is a \cdp-covering family of $t'_{\!_l}$ in $\PropSch/S$, and hence the set
			\[
			\scU\coloneqq \big\{\underline{l}\circ \sigma_{\!_l}\circ i_{\!_{\alpha,l}}:t_{\!_{\alpha,l}}\to t\mid \alpha\in A  \big\}\bigcup \big\{\sigma\circ i_{\!_\alpha}:t_{\!_\alpha}\to t\mid \alpha\in A  \big\}
			\]
			is a \cdp-covering family in $\PropSch/S$. For every $\alpha\in A$, we distinguish two cases.
			\begin{enumerate}
				\item When $t_{\!_{\alpha,j_{\!_{l}}}}$ is nonempty, the open immersion $j_{\!_{l_{\!_\alpha}}}$ is dominant, and hence $a'_{\!_{0}}\circ \sigma \circ i_{\!_{\alpha}}=a'_{\!_{1}}\circ \sigma \circ i_{\!_{\alpha}}$, by \cite[Reduced-to-Separated Th.10.2.2]{Vak:15:FAG} and \eqref{Eq:Equal:Dominant}. Thus, one has $(\sigma\circ i_{\!_{\alpha}})^\ast(a_{\!_{0}})=
				[a'_{\!_{0}}\circ \sigma \circ i_{\!_{\alpha}}]=
				[a'_{\!_{1}}\circ \sigma \circ i_{\!_{\alpha}}]=
				(\sigma\circ i_{\!_{\alpha}})^\ast(a_{\!_{1}}).$
				\item When $t_{\!_{\alpha,j_{\!_{l}}}}\cong \emptyset_{\!_S}$, the morphism $i_{\!_\alpha}$ factorises through $l'$, and hence there exists $a'_{\!_{\alpha,k}}$ in $\yon_{\!_{w,+}}(t_{\!_{\alpha}})$ with $[a'_{\!_{\alpha,k}}]=(\sigma\circ i_{\!_{\alpha}})^\ast(a_{\!_{k}})\in\big(\yon_{\!_{w,+}}\coprod_{\yon_{\!_{z,+}}} \yon_{\!_{p,+}}\big)(t_{\!_{\alpha}})$, for $k=0,1$. Since $\Theta_{\!_t}(a_0)=\Theta_{\!_t}(a_1)$, one has $l_\ast(a'_{\!_{\alpha,0}})=l_\ast(a'_{\!_{\alpha,1}})$, and hence $a'_{\!_{\alpha,0}}=a'_{\!_{\alpha,1}}$ as $l$ is a monomorphism of proper $S$-schemes. Therefore, $(\sigma\circ i_{\!_{\alpha}})^\ast(a_{\!_{0}})=(\sigma\circ i_{\!_{\alpha}})^\ast(a_{\!_{1}})$.
			\end{enumerate}
			Thus, for every $\alpha\in A$, one has
			\[
			(\sigma\circ i_{\!_{\alpha}})^\ast(a_{\!_{0}})=
			(\sigma\circ i_{\!_{\alpha}})^\ast(a_{\!_{1}})\andd
			(\underline{l}\circ \sigma_{\!_l}\circ i_{\!_{\alpha,l}})^\ast(a_{\!_{0}})=
			(\underline{l}\circ \sigma_{\!_l}\circ i_{\!_{\alpha,l}})^\ast(a_{\!_{1}}).
			\]
		\end{enumerate}
		Therefore, there always exists a \cdp-covering family $\scV$ of $t$ in $\PropSch/S$ such that $\delta^\ast (a_0)=\delta^\ast (a_1)$, for every $\delta\in \scV$, \ie the morphism $\Theta$ is a $\tau$-local monomorphism, by Lemma \ref{Lem:Local_Morphism}. Therefore, the square \eqref{Pushout:Square:Modified:Ring} is a pushout square in $\scC^\omega_\tau(S)$, \ie the morphism $f_\shriek$ is an isomorphism in $\scC^\omega_\tau(S)$.
		
		\begin{equation}\label{Universal:Homeomorphism}
		\begin{tikzpicture}[descr/.style={fill=white}]
		\node (W1) at (0,0) {$t'_{\!_l}$};
		\node (X1) at (4,0) {$t'_{\!_{j_{\!_l}}}$};
		\node (XBar1) at (2,0) {$t'$};
		\node (W0) at (0,2) {$t_{\!_{\alpha,l}}$};
		\node (X0) at (4,2) {$t_{\!_{\alpha,j_{\!_{l}}}}$};
		\node (XBar0) at (2,2) {$t_{\!_\alpha}$};
		\node at (.25,1.6) {$\ulcorner$};
		\node at (3.75,1.6) {$\urcorner$};
		\path[->,font=\scriptsize]
		(XBar0) edge node[left]{$i_{\!_\alpha}$} (XBar1)
		(X0) edge node[right]{$i_{\!_{\alpha,j_{\!_{l}}}}$} (X1)
		(W0) edge node[left]{$i_{\!_{\alpha,l}}$} (W1)
		;
		\path[left hook->,font=\scriptsize]
		(X1) edge node[descr]{$\!\!\!\circ\!\!\!$}node[below]{$j_{\!_{l'}}$} (XBar1)
		(X0) edge node[descr]{$\!\!\!\circ\!\!\!$}node[above]{$j_{\!_{l_{\!_\alpha}}}$} (XBar0)
		;
		\path[right hook->,font=\scriptsize]
		(W0) edge node{$\diagup$}node[above]{$l_{\!_\alpha}$} (XBar0)
		(W1) edge node{$\diagup$}node[below]{$l'$} (XBar1)
		;
		\end{tikzpicture}
		\end{equation}
		\item[{$(\ref{itm:Factorisation:Modified:iso}\Rightarrow \ref{itm:Factorisation:Modified:uh})$}] Assume that $f$ is proper and that $f_\shriek$ is an isomorphism, \ie $\Theta$ is a $\tau$-local isomorphism.
		
		\medskip
		
		For surjectivity, suppose that $\fk$ is an algebraically closed field, let $b:\Spec \fk \to Y$ be a morphism in $\fdNoe$, let $l':V\cim Y$ be the scheme-theoretic image of $b$, and let $\underline{b}:\Spec \fk \to V$ be the unique morphism in $\fdNoe$ for which $b=l'\circ \underline{b}$. For $v\coloneqq y\circ l'$, there exists a Cartesian square
		\[
		\begin{tikzpicture}[descr/.style={fill=white}]
		\node (Z) at (2,0) {$t$};
		\node (Y) at (0,-2) {$y$};
		\node (Q) at (2,-2) {$q$};
		\node (X) at (0,0) {$v$};
		\node at (.25,-.4) {$\ulcorner$};
		\path[right hook->,font=\scriptsize]
		(Z) edge node{$\diagup$} node[right]{$\underline{l}$}(Q)
		(X) edge node{$\diagup$} node[left]{$l'$}(Y)
		(Y) edge node[descr]{$\!\!\!\circ\!\!\!$}node[below]{$j_{\!_l}$} (Q)
		(X) edge node[descr]{$\!\!\!\circ\!\!\!$}node[above]{$\underline{j}$} (Z)
		;
		\end{tikzpicture}
		\]
		of $S$-schemes, where $\underline{l}$ is the scheme-theoretic image of the immersion $j_{\!_l}\circ l'$, by Lemma \ref{Lem:Open:Closed:Closed:Open}. Since $\underline{l}$ is a section in $\yon_{\!_{q,+}}(t)$, there exists a $\tau$-cover $\sigma:t'\to t$ and a section $a$ in $\yon_{\!_{w,+}}(t')\coprod_{\yon_{\!_{z,+}}(t')}\yon_{\!_{p,+}}(t')$ such that $\Theta_{t'}(a)=\underline{l}\circ \sigma$. The section $a$ is not the base point in $\yon_{\!_{w,+}}(t')\coprod_{\yon_{\!_{z,+}}(t')}\yon_{\!_{p,+}}(t')$ as $\underline{l}$ is not the base point of $\yon_{\!_{q,+}}(q')$. Since $\tau$ is coarser than the {proper} pretopology, the morphism $\sigma$ is surjective, and hence the morphism $\underline{j}\circ \underline{b}$ lifts to a morphism $b':\Spec \fk \to T'$ along $\sigma$, where $T$ is the underlying scheme of $t'$. 
		
		\medskip
		
		Assume for the sake of contradiction that there exists a morphism $a':t'\to w$ for which $a=[a']$, and hence $l\circ a'=\Theta_{t'}(a)=\underline{l}\circ \sigma$, then the morphism $j_{\!_l}\circ b$ lifts along $l$, which contradict the definition of $j_{\!_l}$ as a complementary open immersion to $l$. Therefore, there exists a morphism $a':t'\to p$ for which $a=[a']$, \ie $g\circ a'=\underline{l}\circ \sigma$. 
		
		\medskip
		
		The morphism $a'\circ b'$ lifts either along $i$ or $j_{\!_i}$. Assuming that $a'\circ b'$ lifts along $i$ implies that exists a morphism $a':t'\to w$ for which $a=[a']$, which leads to a contradiction. Thus, $a'\circ b'$ lifts along $j_{\!_i}$ to a morphism $\underline{a}:\Spec \fk \to X$, \ie $a'\circ b'=j_{\!_i}\circ\underline{a}$. Then,
		\[
		j_{\!_l}\circ f\circ \underline{a}= g\circ j_{\!_i} \circ \underline{a}=g\circ a'\circ b'=\underline{l}\circ \sigma\circ b'=\underline{l}\circ \underline{j}\circ \underline{b}=j_{\!_l}\circ l'\circ \underline{b}=j_{\!_l}\circ b,
		\]
		and hence $f$ is surjective, as $j_{\!_l}$ is a monomorphism in $\fdNoe$.
		
		\medskip

		On the other hand, for universal injectivity, suppose that $\fk$ is an algebraically closed field, and let $b_0,b_1:\Spec \fk \to X$ be morphisms of schemes in $\fdNoe$ such that $f\circ b_0=f\circ b_1$, let $i_0:Z_0\cim P$ and $i_1:Z_1\cim P$ be the scheme-theoretic images of $j_{\!_i}\circ b_0$ and $j_{\!_i}\circ b_0$, respectively, and let $b'_1:\Spec \fk \to Z_0$ (\resp. $b'_0:\Spec \fk \to Z_1$) denote the unique morphism in $\fdNoe$ for which $j_{\!_i}\circ b_0=i_0\circ b'_1$ (\resp. $j_{\!_i}\circ b_1=i_1\circ b'_0$). Let $z_0\coloneqq p\circ i_0$ (\resp. $z_1\coloneqq p\circ i_1$), and denote by $g_0:z'\to z_1$ (\resp. $g_1:z'\to z_0$) the base change in $\PropSch/S$ of $g\circ i_0$ (\resp. $g\circ i_1$) along $g\circ i_1$ (\resp. $g\circ i_0$). Then, there exists a unique morphism of schemes $b:\Spec \fk \to Z'$ in $\fdNoe$ such that $b'_0=g_0\circ b$ and $b'_1=g_1\circ b$, where $Z'$ is the underlying scheme of $z'$. 	 Thus, $i_0\circ g_1, i_1\circ g_0\in \yon_{p,+}(z')$ and $g_\ast(i_0\circ g_1)=g_\ast(i_1\circ g_0)$. Since $\Theta$ is a $\tau$-local monomorphism and $g\circ j_{\!_i}\circ b_0=g\circ j_{\!_i}\circ b_1$ does not factorise through $l$, there exists a $\tau$-cover $\sigma:\underline{z}\to z'$ for which $i_0\circ g_1\circ \sigma = i_1\circ g_0\circ \sigma$. The morphism $b$ lifts along $\tau$-covers as the latter are surjective, in particular, it lifts along $\sigma$. Thus, one has $j_{\!_i}\circ b_0=j_{\!_i}\circ b_1$, and hence $b_0=b_1$, as $j_{\!_i}$ is a monomorphism in $\fdNoe$. Then, $f$ is universally injective, by \cite[\S.3.5.5]{EGA1}.
		
		\medskip
		
		Therefore, the proper universally injective morphism $f$ is finite, by \cite[Th.8.11.1]{EGA4.3}, and hence it is a universal homeomorphism, by \cite[Prop.2.4.5]{EGA4.2}.
	\end{itemize}
\end{proof}
\Cor{}{The motivic measure $\mu_{\!_\tau}$ factorises through the the motivic measure $\mu_{\!_\uh}$, which takes value in the modified Grothendieck ring of $S$-schemes.}{Factorisation:Modified}

Some partial results, based on an analogue of \cite[Th.3.2.9]{Voe:96:HS}, led us in the direction of the following conjecture. However, we do not have a full prove, yet.
\Conj{}{The motivic measure $\mu_{\!_\cdp}$  is isomorphic to the motivic measure $\mu_{\!_\uh}$.	
}{Universal:Measure:Modified}
The validity of this conjecture would extend the modified Grothendieck ring of $S$-schemes, through the Waldhausen $\Kg$-theory spectrum $\Kg\big(\scC^\omega(S)\big)$, which we investigate in a future work.

\subsection{The Commutative Ring Spectrum Structure}
The spectrum $\Kg\big(\scC^\omega_\tau(S)\big)$ admits a canonical commutative ring spectrum structure, \ie a homotopy commutative monoid structure in {the category of $S^1$-spectra of pointed topological spaces}.

\medskip

The functor $\ast\to \PropSch/S $ that sends the unique object of $\ast$ to $\id_{S}$ is continuous, with respect to the indiscrete topology on $\ast$ and any topology on $\PropSch/S$. Hence, it induces an exact functor of Waldhausen categories $u_{\!_S}^\tau:\PSh_\bullet(\ast)\cong\Set_\bullet \to \scC_\tau(S)$, which is a left adjoint to the the global section functor $\scC_\tau(S)\to \Set_\bullet$, \ie $u_{\!_S}^\tau$ is given by sending a pointed sets to its constant pointed $\tau$-sheaf. Recall that the category $\FSet_\bullet$ of pointed finite sets admits a symmetric monoidal Waldhausen structure, as seen in Example \ref{Ex:Waldhausen:Finite:Sets}. The category $\scC^\omega_\tau(S)$ contains the unit $1_{\wedge}=\underline{\yon}^\tau_{\!_{\id_{\!_S}}}$, admits all finite colimits, and the $\tau$-sheafification functor commutes with colimits. Thus, the functor $u_{\!_S}^\tau$ restricts to the exact functor of Waldhausen categories $\upsilon_{\scC^\omega_\tau(S)}:\FSet_\bullet \to \scC^\omega_\tau(S)$, given in \eqref{Unit:Waldhausen:Upsilon}, which we denote by $ \upsilon_{\!_S}^\tau$.

\medskip
On the other hand, $\scC^\omega_\tau(S)$ is a symmetric monoidal Waldhausen category, by Lemma \ref{Lem:Symmetric:Monoidal}. Hence, there exists a paring $\otimes:\Kg\big(\scC^\omega_\tau(S)\big)\wedge \Kg\big(\scC^\omega_\tau(S)\big)\to \Kg\big(\scC^\omega_\tau(S)\big)$, see \cite[p.342]{Wal:85:AKTS}. That makes $\Kg\big(\scC^\omega_\tau(S)\big)$ into a commutative ring spectrum, see \cite[Cor.2.8]{BM:11:DKDIAKTS}. 

\medskip

\subsection{The Monoidal Proper-Fibred Waldhausen Category}\label{Fibred:Noetherian}
The $\Kg$-theory commutative ring spectrum $\Kg\big(\scC^\omega_\tau(S)\big)$, for a Noetherian scheme $S$ of finite Krull dimension, arises from a fibre of a \indx{monoidal proper-fibred Waldhausen category} over Noetherian schemes of finite Krull dimensions\footnote{See \cite[\S.1]{CD:13:TCMM} for a treatment of $\scP$-fibred categories, for a set $\scP$ of morphisms of schemes.}. That is, there exists a strong monoidal pseudofunctor $\scC^\omega_\tau:\fdNoe^\op\to \Wald^{\wedge}_2$, where $\Wald_2^{\wedge}$ is the $2$-category of essentially small symmetric monoidal Waldhausen categories, weak monoidal exact functors between them, and monoidal natural transformations between the latter, such that 
\begin{itemize}
	\item for every scheme $S\in \fdNoe$ the fibre $\scC^\omega_\tau(S)$ is the symmetric monoidal Waldhausen category constructed in \eqref{Waldhausen:Omega}, by Lemma \ref{Lem:Restricted:Pullback:Exact};
	\item for every proper morphism $f:S\to T$ in $\fdNoe$, the pullback $f^\ast:\scC^\omega_\tau(T)\to \scC^\omega_\tau(S)$ admits an exact left adjoint $f_{\#}:\scC^\omega_\tau(S)\to \scC^\omega_\tau(T)$, by Lemma \ref{Lem:Restricted:Pushforward:Exact};
	\item $\scC^\omega_\tau$ satisfies the proper-base change property, by Lemma \ref{Lem:Base:Change:Noetherian}; and
	\item $\scC^\omega_\tau$ satisfies the proper-projection formula, by Lemma \ref{Lem:Noetherian:Projection:Formula}.
\end{itemize}
Then, applying the \waldhausen's $\Kg$-theory $2$-functor\footnote{See \cite[\S.1]{FM:17:WACQ} for the treatment of the $2$-categorical \waldhausen's $\Kg$-theory.} induces the monoidal proper-fibred commutative ring spectrum $\Kg(\scC^\omega_\tau):\fdNoe^\op \to {\CRing\Spt_2}$. In fact, the strong monoidal \cdp-functor $\underline{\yon}^\tau:\PropSch/S \to \scC^\omega_\tau(S)$, given in \eqref{cdp:Compactifiable:tau:Omega}, arises from the \indx{geometric section} of $\scC^\omega_\tau$, see \cite[\S.1.1.c]{CD:13:TCMM}.

\Rem{}{
	The statements in the rest of this subsection were motivated by Dan Petersen's answer in \cite{MO:Pat:14:HGR}, which recalls Ekedahl's approach to higher \grothendieck groups of varieties. The statements of Lemma \ref{Lem:Base:Change:Noetherian} and Lemma \ref{Lem:Noetherian:Projection:Formula} are essentially consequences of \cite[Ex.1.1.11 and Ex.1.1.28]{CD:13:TCMM}.
}{}

\subsubsection{Inverse Image}\label{Pullbacks:Essential:Part}Recall the canonical proper-fibred category $\PropSch/-:\fdNoe\to \CAT_2$, as in \cite[Ex.1.1.4 and Ex.1.1.11]{CD:13:TCMM}. A morphism $f:S\to T$ in $\fdNoe$ induces a functor $f^\inv:\PropSch/T \to \PropSch/S$, that sends a proper $T$-scheme to its base change along $f$. That in turn induces a direct image functor $f_\ast:\PSh_\bullet(\PropSch/S)\to \PSh_\bullet(\PropSch/T)$, given by precomposition with $(f^\inv)^\op$, \ie for a presheaf $\scX\in \PSh_\bullet(\PropSch/S)$ and for a proper $T$-scheme $q$, one has $
f_\ast(\scX)(q)=\scX \big(f^\inv(q)\big)$. The functor $f_\ast$ admits a left adjoint $f^\ast_\pre:\PSh_\bullet(\PropSch/T)\to \PSh_\bullet(\PropSch/S)$, called the \indx{inverse image functor} along $f$, and it is given by a left \kan extension along $(f^\inv)^\op$, see \S.\ref{Continuous:Maps}. The functor $-_{\!_+}$ commutes with colimits, for being a left adjoint. Then, using the coend formula \cite[\S.X.4.(1)]{ML:98:CWM}, one sees that 
\begin{align*}
f^\ast_\pre(\yon_{\!_{q,+}})(p)&=
\bigg(\int^{x\in \PropSch/T} {\PropSch/S\big(p,f^\inv(x)\big)} \times \PropSch/T(x,q)\bigg)_{\!_+}\\&\cong
\bigg(\int^{x\in \PropSch/T} \PropSch/S\big(p,f^\inv(x)\big) \times \PropSch/S\big(f^\inv(x),f^\inv(q)\big)\bigg)_{\!_+}\\&\cong
\bigg( \PropSch/S\big(p,f^\inv(q))\bigg)_{\!_+}\cong
\yon_{\!_{f^\inv(q),+}}(p),
\end{align*}
for every proper $T$-scheme $q$ and for every proper $S$-scheme $p$, \ie $f^\ast_\pre(\yon_{\!_{q,+}})\cong\yon_{\!_{f^\inv(q),+}}$. In fact, this is a defining property for $f^\ast_\pre$, as every object in $\PSh_\bullet(\PropSch/T)$ is a colimit of a diagram in the essential image of $\yon_{\!_{-,+}}$.

\medskip

Since the base change functor $f^\inv$ commutes with fibre products, it is continuous with respect to the $\tau$-pretopology, see \cite[\S.III.Prop.1.6]{SGA4}. Thus, the direct image functor $f_\ast$ preserves $\tau$-sheaves, and it restricts to a functor $f_\ast :\scC_\tau(S)\to \scC_\tau(T)$, which admits a left adjoint $f^\ast$, given by the composition of $f^\ast_\pre$ with the associated $\tau$-sheaf functor $-^{\!^{\ash_\tau}}$. Since the $\tau$-sheafification functor commutes with colimits, one has
\begin{equation}\label{Geometric:Section:Eq:1}
f^\ast(\underline{\yon}^\tau_{\!_{q}})\cong \underline{\yon}^\tau_{\!_{f^\inv(q)}},
\end{equation}
for every proper $T$-scheme $p$.

\Lem{}{Assume that $f:S\to T$ is a morphism in $\fdNoe$. Then, the functor $f^\ast:\scC_\tau(T) \to \scC_\tau(S)$ is a strong monoidal exact functor of symmetric monoidal Waldhausen categories.}{Pullback:Exact}
\begin{proof}The functor $f^\inv$ is Cartesian, as limits commute with each other, and hence the functor $f^\ast_\pre$ is left exact, see \cite[A.Ex.4.1.10]{Joh:02:SE}. Also both the $\tau$-sheafification functor and the inclusion functor, of $\tau$-sheaves into presheaves, are left exact. Thus, the left adjoint functor $f^\ast$ is left exact. In particular, the functor $f^\ast$ preserves monomorphisms, finite colimits, the unit of the monoidal structure, and smash products of pointed $\tau$-sheaves, as the latter only involves finite limits and colimits of pointed $\tau$-sheaves. 
\end{proof}

Although the functor $f_\ast:\PSh_\bullet(\PropSch/S)\to \PSh_\bullet(\PropSch/T)$ admits a right adjoint given by the right \kan extension, its resection $f_\ast :\scC_\tau(S)\to \scC_\tau(T)$ does not necessarily admit a right adjoint. In particular, it is not necessarily exact.

\Lem{}{Assume that $f:S\to T$ is a morphism in $\fdNoe$. Then, the functor $f^\ast:\scC_\tau(T) \to \scC_\tau(S)$ restricts to a strong monoidal exact functor $f^{\ast,\omega}:\scC^\omega_\tau(T) \to \scC^\omega_\tau(S)$.
}{Restricted:Pullback:Exact}
\begin{proof}The statement follows from Lemma \ref{Lem:Pullback:Exact}, provided the restriction $f^{\ast,\omega}$ exists. Since $f^\ast$ commutes with pushouts, it suffices to show that $f^\ast$ restricts to a functor $f^{\ast,0}:\scC^0_\tau(T) \to \scC^\omega_\tau(S)$ to induce a strong monoidal exact functor $f^{\ast,\omega}:\scC^\omega_\tau(T) \to \scC^\omega_\tau(S)$, which holds by \eqref{Geometric:Section:Eq:1}.
\end{proof}
When no confusion arises, we abuse notation, and refer to $f^{\ast,\omega}$ by $f^\ast$.

\Cor{}{Assume that $f:S\to T$ is a morphism in $\fdNoe$. Then, the there exists a {morphism of commutative ring spectra} $
	f^\ast:\Kg\big(\scC^\omega_\tau(T)\big) \to \Kg\big(\scC^\omega_\tau(S)\big)$, that sends a point in the component $[\underline{\yon}^\tau_{\!_q}]$ to a point in the component $[\underline{\yon}^\tau_{\!_{f^\inv(q)}}]$, for every proper $T$-scheme $q$.
}{}

\Eg{}{Assume that $f:S\to T$ is a morphism in $\fdNoe$, let $y$ be a $T$-scheme, and let $l:w\cim q$ be a compactification of $y$. Since $f^\ast$ is exact and complementary open immersions are closed under pullbacks, one has
	\[
	f^\ast(\yon^{\ps_\tau}_{\!_y})\cong
	\coker \big(f^\ast(\underline{\yon}^\tau_{\!_w})\rightarrowtail f^\ast(\underline{\yon}^\tau_{\!_q}) \big)\cong 
	\coker \big(\underline{\yon}^\tau_{\!_{f^\inv(w)}}\rightarrowtail \underline{\yon}^\tau_{\!_{f^\inv(q)}} \big)\cong
	\yon^{\ps_\tau}_{\!_{f^\inv(y)}},
	\]
	where the $S$-scheme $f^\inv(y)$ is a base change in $\sftSch/T$ of $y$ along $f$. Then, in particular,
	\begin{equation}\label{Pullback:0}
	\pi_0(f^\ast)\big([\yon^{\ps_\tau}_{\!_y}]\big)=[f^\ast(\yon^{\ps_\tau}_{\!_y})]=[\yon^{\ps_\tau}_{\!_{f^\inv(y)}}].
	\end{equation}
	Since the motivic measure $\mu_{\!_{\tau,T}}$ is surjective, by Lemma \ref{Lem:Surjective:Measure}, one sees that the ring homomorphism $
	\pi_0(f^\ast):\Kg_0\big(\scC^\omega_\tau(T)\big) \to \Kg_0\big(\scC^\omega_\tau(S)\big)$ is determined by \eqref{Pullback:0}. We may abuse notation and refer to $\pi_0(f^\ast)$ by $f^\ast$.
}{}

Let $f:S\to T$ and $g:T\to U$ be morphisms in $\fdNoe$. Since \kan extensions are determined up to canonical natural isomorphisms, one has canonical {monoidal} natural isomorphisms ${(g\circ f)}^{\ast}\stackrel{\sim}{\Rightarrow} f^{\ast} \circ g^{\ast}$ and ${(\id_{\!_S})}^{\ast}\stackrel{\sim}{\Rightarrow}\id_{\!_{\scC^\omega_\tau(S)}}$, {which satisfy the cocycle condition.}

\Cor{}{There exists a pseudofunctor $\scC^\omega_\tau:\fdNoe^\op \to \Wald^{\wedge}_2$, which sends a Noetherian scheme $S$ of finite Krull dimension to the Waldhausen category $\scC^\omega_\tau(S)$, given in \eqref{Waldhausen:Omega}, and sends $f^\op$, for a morphism $f:S\to T$ in $\fdNoe$, to the strong monoidal exact functor $f^{\ast,\omega}$, given in Lemma \ref{Lem:Restricted:Pullback:Exact}. Then, the \waldhausen's $\Kg$-theory $2$-functor induces a {pseudofunctor} $\Kg(\scC^\omega_\tau):\fdNoe^\op \to {\CRing\Spt_2}$.
}{Waldhausen:Noetherian:Open}

\subsubsection{Proper Direct Image}
Suppose that $f:S\to T$ is a proper morphism in $\fdNoe$. Then, the functor $f^\inv$, given in \S.\ref{Pullbacks:Essential:Part}, admits a left adjoint $f_{\!_\circ}:\PropSch/S \to \PropSch/T$, given by composition with $f$. Thus, the functor $f^!:\PSh_\bullet(\PropSch/T)\to \PSh_\bullet(\PropSch/S)$, given by precomposition with $f_{\!_\circ}$, is a left adjoint to $f_\ast:\PSh_\bullet(\PropSch/S)\to \PSh_\bullet(\PropSch/T)$, and hence $f^!$ is canonically isomorphic to $f^\ast_\pre$. The functor $f^!$ admits a left adjoint $f_{\#,\pre}$, given by a left \kan extension along $f_{\!_\circ}^\op$. Since the functor $f_{\circ}$ preserves $\tau$-coving families, $f^!$ restricts to a functor $f^!:\scC_\tau(T)\to \scC_\tau(S)$, which is canonically isomorphic to $f^\ast$. The functor $f^!$ admits a left adjoint $f_{\#}:\scC_\tau(S)\to \scC_\tau(T)$, called the \indx{proper direct image functor} along $f$, and it is given by the $\tau$-sheafification of $f_{\#,\pre}$. Similar to the inverse image functor, for a proper $S$-scheme $p$, one has
\begin{equation}\label{Geometric:Section:Eq:2}
f_{\#}(\underline{\yon}^\tau_{\!_{p}})\cong \underline{\yon}^\tau_{\!_{f_{\!_\circ}(p)}},
\end{equation}
which is a defining property of $f_{\#}$.

\Lem{}{
	Assume that $f:S\to T$ is a proper morphism in $\fdNoe$. Then, the functor $f_{\#}:\scC_\tau(S) \to \scC_\tau(T)$ is an exact functor of Waldhausen categories.
}{Pushforward:Exact}
\begin{proof} It is sufficient to show the functor $f_{\#}$ commutes with monomorphisms, as it commutes with colimits for being a left adjoint. 
	
	\medskip

	Assume that $\iota:\scX\to \scY$ is a monomorphism in $\scC_\tau(S)$, let $q$ be a proper $T$-scheme, and let $t_0,t_1\in f_{\#}(\scX)(q)$ such that $f_{\#}(\iota)_q(t_0)=f_{\#}(\iota)_q(t_1)$. 
	
	\medskip 
	
	By the definition of the $\tau$-sheafification functor, as in \cite[Def.2.63]{Vis:08:NGTFCDT}, there exists a $\tau$-covering family $\scU=\{\sigma_\alpha:q_\alpha\to q\mid \alpha\in A\}$, and there exists a section $t_{k,\alpha}\in f_{\#,\pre}(\scX)(q_\alpha)$ such that $\sigma_\alpha^\ast(t_k)=t_{k,\alpha}^{\ash}$ for every $\alpha\in A$, for $k=0,1$. For $\alpha\in A$, pulling back along $\sigma_\alpha$ yields 
	\[
	{\big(f_{\#}(\iota)_{q_\alpha}(t_{0,\alpha})\big)}^{\ash}=f_{\#}(\iota)_{q_\alpha}(t_{0,\alpha}^{\ash})=f_{\#}(\iota)_{q_\alpha}(t_{1,\alpha}^{\ash})={\big(f_{\#,\pre}(\iota)_{q_\alpha}(t_{1,\alpha})\big)}^{\ash}.
	\]
	Thus, there exits a $\tau$-covering family $\scU_\alpha=\{\sigma_{\alpha,\beta}:q_{\alpha,\beta}\to q_\alpha\mid \beta\in B_\alpha\}$ for which
	\begin{align*}
	f_{\#,\pre}(\iota)_{q_{\alpha,\beta}} \big(\sigma_{\alpha,\beta}^\ast(t_{0,\alpha})\big)&=
	\sigma_{\alpha,\beta}^\ast \big(f_{\#,\pre}(\iota)_{q_\alpha}(t_{0,\alpha})\big)=
	\sigma_{\alpha,\beta}^\ast \big(f_{\#,\pre}(\iota)_{q_\alpha}(t_{1,\alpha})\big)\\&=
	f_{\#,\pre}(\iota)_{q_{\alpha,\beta}} \big(\sigma_{\alpha,\beta}^\ast(t_{1,\alpha})\big).
	\end{align*}
	The functor $f_{\#,\pre}$ is a left \kan extension along $f_{\!_\circ}^\op$. For a proper $T$-scheme $q'$, the coend formula \cite[\S.X.4.(1)]{ML:98:CWM} implies that the underlying set of $f_{\#,\pre}(\scX)(q')$ can be given by
	\[
	\bigg(\bigsqcup_{p\in \PropSch/S} \PropSch/T(q',f\circ p)\times \scX(p)\bigg)\big/\sim,
	\]
	where $\sim$ is the smallest equivalence relation that identifies $(g,s)\in \PropSch/T(q',f\circ p)\times \scX(p)$ and $(g',s')\in \PropSch/T(q',f\circ p')\times \scX(p')$ whenever there exists a morphism $h:p\rightarrow p'$ of proper $S$-schemes for which $g'=f_{\!_\circ}(h)\circ g $ and $s=h^\ast (s')$; whereas the point of $f_{\#,\pre}(\scX)(q')$ is given by the unique class $[(g,\ast)]$, which is independent of the choice of the proper $S$-scheme $p$ and the morphism $g: q'\to f\circ p$ of proper $T$-schemes. Also, one has $f_{\#,\pre}(\iota)_{q'} \left([(g,s)]\right)=\big[\big(g,\iota_{p}(s)\big)\big]$, for every proper $S$-scheme $p$, every morphism $g:q'\to f\circ p$ of proper $T$-schemes, and for every section $s\in \scX(p)$.
	
	\medskip
	
	For $k=0,1$, for $\alpha\in A$, and for $\beta\in B_\alpha$, let $p_{k,\alpha,\beta}$ be a proper $S$-scheme, let $g_{k,\alpha,\beta}:q_{\alpha,\beta}\to f\circ p_{k,\alpha,\beta}$ be a morphism of proper $T$-schemes, and let $s_{k,\alpha,\beta}\in \scX(p_{k,\alpha,\beta})$ for which $\sigma_{\alpha,\beta}^\ast(t_{k,\alpha})=[(g_{k,\alpha,\beta},s_{k,\alpha,\beta})]$. Since 
	\begin{align*}
	\big[\big(g_{0,\alpha,\beta},\iota_{p_{0,\alpha,\beta}}(s_{0,\alpha,\beta})\big)\big]&=
	f_{\#,\pre}(\iota)_{q_{\alpha,\beta}}\big([(g_{0,\alpha,\beta},s_{0,\alpha,\beta})]\big)=
	f_{\#,\pre}(\iota)_{q_{\alpha,\beta}}\big([(g_{1,\alpha,\beta},s_{1,\alpha,\beta})]\big)\\&=
	\big[\big(g_{1,\alpha,\beta},\iota_{p_{1,\alpha,\beta}}(s_{1,\alpha,\beta})\big)\big],
	\end{align*}
	there exist a proper $S$-scheme $p_{\alpha,\beta}$, a morphism $g_{\alpha,\beta}:q_{\alpha,\beta}\to f\circ p_{\alpha,\beta}$ of proper $T$-schemes, and a morphism
	$h_{k,\alpha,\beta}:p_{\alpha,\beta}\to p_{k,\alpha,\beta}$ of proper $S$-schemes, such that $g_{k,\alpha,\beta}=f_{\!_\circ}(h_{k,\alpha,\beta})\circ g_{\alpha,\beta}$ and
	\[
	\iota_{p_{\alpha,\beta}}\big(h_{0,\alpha,\beta}^\ast(s_{0,\alpha,\beta})\big)=
	h_{0,\alpha,\beta}^\ast \big(\iota_{p_{0,\alpha,\beta}}(s_{0,\alpha,\beta})\big)=h_{1,\alpha,\beta}^\ast \big(\iota_{p_{1,\alpha,\beta}}(s_{1,\alpha,\beta})\big)=\iota_{p_{\alpha,\beta}}\big(h_{1,\alpha,\beta}^\ast(s_{1,\alpha,\beta})\big).
	\]
	Since $\iota$ is a monomorphisms, one has $h_{0,\alpha,\beta}^\ast(s_{0,\alpha,\beta})=h_{1,\alpha,\beta}^\ast(s_{1,\alpha,\beta})$, and hence 
	\[
	\sigma_{\alpha,\beta}^\ast(t_{0,\alpha})=[(g_{0,\alpha,\beta},s_{0,\alpha,\beta})]=[(g_{1,\alpha,\beta},s_{1,\alpha,\beta})]=\sigma_{\alpha,\beta}^\ast(t_{1,\alpha}).
	\]
	Thus, $
	(\sigma_\alpha\circ \sigma_{\alpha,\beta})^\ast (t_0)=\sigma_{\alpha,\beta}^\ast(t_{0,\alpha})^{\ash}=\sigma_{\alpha,\beta}^\ast(t_{1,\alpha})^{\ash}=(\sigma_\alpha\circ \sigma_{\alpha,\beta})^\ast (t_1)$. Since $\{\sigma_\alpha\circ \sigma_{\alpha,\beta}:q_{\alpha,\beta}\to q\mid \alpha\in A, \beta\in B_\alpha\}$ is a $\tau$-covering family in $\PropSch/T$ and $f_{\#}(\scX)$ is a $\tau$-sheaf, one has $t_0=t_1$. 
\end{proof}

\Lem{}{Let $f:S\to T$ be a proper morphism in $\fdNoe$. Then, the functor $f_{\#}:\scC_\tau(S) \to \scC_\tau(T)$ restricts to an exact functor of Waldhausen categories $f_\#^{\omega}:\scC^\omega_\tau(S) \to \scC^\omega_\tau(T)$.}{Restricted:Pushforward:Exact}
\begin{proof}
	The proof is essentially the same as of the proof of Lemma \ref{Lem:Restricted:Pullback:Exact}, utilising Lemma \ref{Lem:Pushforward:Exact} and \eqref{Geometric:Section:Eq:2} instead of Lemma \ref{Lem:Pullback:Exact} and \eqref{Geometric:Section:Eq:1}.
\end{proof}
When no confusion arises, we abuse notation, and refer to $f_\#^{\omega}$ by $f_\#$.
\Cor{}{Assume that $f:S\to T$ is a proper morphism in $\fdNoe$. Then, the functor $f_\#:\scC^\omega_\tau(S) \to \scC^\omega_\tau(T)$ induces a {morphism of spectra} $f_\#:\Kg\big(\scC^\omega_\tau(S)\big) \to \Kg\big(\scC^\omega_\tau(T)\big)$, that sends a point in the component $[\underline{\yon}^\tau_{\!_p}]$ to a point in the component $[\underline{\yon}^\tau_{\!_{f_{\!_\circ}(p)}}]$, for a proper $S$-scheme $p$.
}{}
Also, we may abuse notation if no confusion arises and refer to $\pi_0(f_\#)$ by $f_\#$.

\Eg{}{Assume that $f:S\to T$ is a proper morphism in $\fdNoe$, let $x$ be an $S$-scheme, and let $i:z\cim p$ be a compactification of $x$. Since $f_\#$ is exact and $f_{\!_\circ}$ preserves complementary open immersions, one has
	\[
	f_\#(\yon^{\ps_\tau}_{\!_x})\cong
	\coker \big(f_\#(\underline{\yon}^\tau_{\!_z})\rightarrowtail f_\#(\underline{\yon}^\tau_{\!_p}) \big)\cong 
	\coker \big(\underline{\yon}^\tau_{\!_{f_{\!_\circ}(z)}}\rightarrowtail \underline{\yon}^\tau_{\!_{f_{\!_\circ}(p)}} \big)\cong
	\yon^{\ps_\tau}_{\!_{f\circ x}}.
	\]
	Thus, in particular, $f_\#\big([\yon^{\ps_\tau}_{\!_x}]\big)=[f_\#(\yon^{\ps_\tau}_{\!_x})]=[\yon^{\ps_\tau}_{\!_{f\circ x}}]$.
}{}
Suppose that $f:S\to T$ and $g:T\to U$ are proper morphisms in $\fdNoe$. Then, there exist canonical natural isomorphisms $(g\circ f)_\#\stackrel{\sim}{\Rightarrow} g_\# \circ f_\#$ and $(\id_{\!_S})_\#\stackrel{\sim}{\Rightarrow}\id_{\!_{\scC^\omega_\tau(S)}}$, which satisfy the cocycle condition.
\Cor{}{The fibred Waldhausen category $\scC^\omega_\tau$, given in Corollary \ref{Cor:Waldhausen:Noetherian:Open}, is in fact a pre-proper-fibred Waldhausen category, \ie there exists a pseudofunctor $\scC^\omega_\tau:\fdNoe^\prop \to \Wald_2$, where $\Wald_2$ is the $2$-category of essentially small Waldhausen categories, exact functors between them, and natural transformations between the latter, which sends a proper morphism $f:S\to T$ in $\fdNoe$ to the exact functor $f_\#^{\omega}$, given in Lemma \ref{Lem:Restricted:Pushforward:Exact}, that is a left adjoint to $f^\ast$. Then, the \waldhausen's $\Kg$-theory $2$-functor induces a {pseudofunctor} $\Kg(\scC^\omega_\tau):\fdNoe^\prop \to \Spt_2$.
}{Noetherian:Proper:Direct:Image}
In contract to the inverse image, the proper direct image is not necessarily strong monoidal. That is, for a proper morphism $f:S\to T$ in $\fdNoe$, one has $f_\# (\1_{\!_S})=
f_\#( \underline{\yon}^\tau_{\!_{\id_{\!_S}}}) \cong
\underline{\yon}^\tau_{\!_{f_{\!_\circ}(\id_{\!_S})}}\cong
\underline{\yon}^\tau_{\!_{f}}$, which is not necessarily isomorphic to $\1_{\!_T}$ for a proper morphism $f$. However, since $f_\#$ is a left adjoint to the strong monoidal functor $f_\ast$, it is oplax monoidal, see \cite[\S.1.1.24]{CD:13:TCMM}.

\Eg{}{Suppose that $p$ is a proper $S$-scheme. Then, $\big[\underline{\yon}^\tau_{\!_p}\big]=\big[(p_\# \circ p^\ast ) (\1_{\!_S})\big]\in \Kg_0\big(\scC^\omega_\tau(S)\big)$.
}{}

\subsubsection{Open Direct Image}\label{Open:Direct:Image}For a morphism $f:S\to T$ in $\fdNoe$, the right adjoint direct image functor $f_\ast :\scC_\tau(S)\to \scC_\tau(T)$ is not necessarily an exact functor of Waldhausen categories, as it may not commute with pushouts. However, when $j$ is an open immersions, the functor $j_\ast$ is a strong monoidal exact functor, as in Corollary \ref{Cor:Open:Direct:Image}. Then, it suffices to show that $j_\ast$ restricts to a functor $j_\ast^0:\scC^0_\tau(S) \to \scC^\omega_\tau(T)$, in order to induce a strong monoidal exact functor $j_\ast^\omega:\scC^\omega_\tau(S) \to \scC^\omega_\tau(T)$. However, it is not clear to us that such a restriction exists. Also, we intended to use the open direct image functor to extend the proper direct image to all morphisms in $\fdNoe$, but it does not seem to provide an extension independent from the choice of compactifications. 

\medskip

For an open immersion $j:S\oim T$ in $\fdNoe$, we first show that the functor $j^\inv$ is almost cocontinuous, as in Definition \ref{Def:Almost:Cocontinuous}, then we apply Lemma \ref{Lem:Almost:Cocontinuous} to deduce that $j_\ast$ is a strong monoidal exact functor. However, we do not pursue studying $j^\inv$ further here, and we leave it for a further work.

\Lem{}{Let $j:S\oim T$ be an open immersion in $\fdNoe$. When $\tau$ is either the \cdp-pretopology or the proper pretopology, the functor $j^\inv$ is continuous and almost cocontinuous with respect to the $\tau$-pretopology.}{Open:Direct:Image}
\begin{proof}
	The functor $j^\inv$ is continuous with respect to the $\tau$-pretopology for preserving $\tau$-covering families, see \cite[\S.III.Prop.1.6]{SGA4}.
	
	\medskip
	
	Assume that $q:Q\to T$ is a proper $T$-scheme, let $q':Q'\to S$ be a base change of $q$ along $j$, and let $\scU=\{\sigma_\alpha:p_\alpha\to q'\mid \alpha\in A\}$ be a $\tau$-covering family in $\PropSch/S$. Recall that schemes of finite type over Noetherian schemes are Noetherian, by \cite[Tag 01T6]{stacks-project}. Since $Q'$ is Noetherian, the open immersion $ j'\coloneqq q^\inv(j)$ is quasi-compact, and hence of finite type, see \cite[Tags 01P0, 01TU, and 01TW]{stacks-project}. Thus, for every $\alpha\in A$, the category of compactifications $\Comp_Q( j'\circ \sigma_\alpha)$ is nonempty, by \nagata's Compactification Theorem. Hence, there exists a proper $Q$-scheme $z_\alpha:Z_\alpha\to Q$ which admits an open immersion $j_\alpha: j'\circ \sigma_\alpha \oim z_\alpha$ in $\sftSch/Q$. Consider the commutative diagram
	\[
	\begin{tikzpicture}[descr/.style={fill=white}]
	\node (NewZ) at (4,0) {$S$};
	\node (NewQ) at (4,-2) {$T$};
	
	\node (Z) at (2,0) {$Q'$};
	\node (fp) at (0,0) {$Z_\alpha\times_{\!_Q} Q'$};
	\node (Y) at (0,-2) {$Z_\alpha$};
	\node (Q) at (2,-2) {$Q$};
	\node (X) at (-1,1) {$P_\alpha$};
	\node at (.25,-.4) {$\ulcorner$};
	\node at (2.25,-.4) {$\ulcorner$};
	\path[right hook->,font=\scriptsize]
	(NewZ) edge node[descr]{$\!\!\!\circ\!\!\!$} node[right]{$j$}(NewQ)
	(Z) edge node[descr]{$\!\!\!\circ\!\!\!$} node[right]{$j'$}(Q)
	(fp) edge node[descr]{$\!\!\!\circ\!\!\!$} node[right]{$\underline{j}_\alpha$}(Y)
	;
	\path[->,font=\scriptsize]
	(Y) edge node[below]{$z_\alpha$} (Q)
	(fp) edge node[below]{$\underline{z}_\alpha$} (Z)
	(Q) edge node[below]{$q$} (NewQ)
	(Z) edge node[above]{$q'$} (NewZ)
	;
	\path[->,dashed,font=\scriptsize]
	(X) edge node[descr]{$\!\!\!k_\alpha\!\!\!$}(fp)
	;
	\path[bend right,right hook->,font=\scriptsize]
	(X) edge node[descr]{$\!\!\!\circ\!\!\!$} node[left]{$j_{\!_\alpha}$}(Y);
	\path[bend left,->,font=\scriptsize]
	(X) edge node [above]{$\sigma_\alpha$} (Z)
	;
	\end{tikzpicture}
	\]
	in $\fdNoe$, where $P_\alpha$ is the underlying scheme of $p_\alpha$, and $k_\alpha$ is the unique morphism $P_\alpha\to Z_\alpha\times_{\!_Q} Q'$ in $\fdNoe$, induced by the universal property of fibre products that makes the diagram commute.
	
	\medskip

	Since $j_{\!_\alpha}$ is an open immersion, so is $k_\alpha$. Also, $k_\alpha$ is proper, as $\sigma_\alpha$ is proper. Then, $k_\alpha$ is a closed open immersion, by \cite[Cor.18.12.6]{EGA4.4}. Let $l_\alpha:W_\alpha\cim Z_\alpha$ be the scheme-theoretic image of the immersion $j_\alpha$, and let $j'_\alpha:P_\alpha\to W_\alpha$ be the unique morphism of $Q$-schemes for which $\underline{j}_\alpha\circ k_\alpha=l_\alpha\circ j'_\alpha$. Then, $j'_\alpha$ is an open immersion, and the square
	\[
	\begin{tikzpicture}[descr/.style={fill=white}]
	\node (Z) at (0,-2) {$W_\alpha$};
	\node (Y) at (2,0) {$Z_\alpha\times_{\!_Q}Q'$};
	\node (Q) at (2,-2) {$Z_\alpha$};
	\node (X) at (0,0) {$P_\alpha$};
	\node at (.25,-.4) {$\ulcorner$};
	\path[right hook->,font=\scriptsize]
	(Z) edge node{$\diagup$} node[below]{$l_\alpha$}(Q)
	(X) edge node{$\diagup$} node[above]{$k_\alpha$}(Y)
	(Y) edge node[descr]{$\!\!\!\circ\!\!\!$}node[right]{$\underline{j}_\alpha$} (Q)
	(X) edge node[descr]{$\!\!\!\circ\!\!\!$}node[left]{$j'_\alpha$} (Z)
	;
	\end{tikzpicture}
	\]
	is Cartesian, by Lemma \ref{Lem:Open:Closed:Closed:Open}.  
	
	\medskip
	
	Let $z:Z\cim Q$ be a closed immersion complementary to $j':Q'\oim Q$, one may choose $Z$ to have the reduced induced structure, but such a choice does not affect the argument. Then, we will see that the set of proper morphisms
	\[
	\scV \coloneqq \{z_\alpha\circ l_\alpha:W_\alpha\to Q\mid \alpha\in A\}\bigsqcup \{z:Z\cim Q\}
	\]
	is a $\tau$-covering family of $Q$. 
	\begin{itemize}
		\item When $\tau$ is the proper pretopology, it is evident that $\scV$ is a proper covering family.
		\item On the other hand, when $\tau$ is the \cdp-pretopology, for every field $\fk$, every morphism $x:\Spec \fk \to Q$ in $\fdNoe$ lifts either through $z$ or $j'$. When $x$ lifts through $j'$ to a morphism $x':\Spec \fk \to Q'$, since $\scU$ is a \cdp-covering, there exists $\alpha_x\in A$, such that $x'$ lifts through $\sigma_{\alpha_x}$, \ie the $\fk$-point $x$ lifts through $z_{\alpha_x}\circ l_{\alpha_x}$. Thus, for every field $\fk$, every $\fk$-point in $Q$ lifts through a morphism in $\scV$, and hence $\scV$ is a \cdp-covering family.
	\end{itemize}
	For every $\alpha\in A$, the morphism $j^\inv(z_\alpha\circ l_\alpha)$ is isomorphic to $p_\alpha$, and hence factorises through it. On the other hand, the empty sieve is a $\tau$-covering sieve on the empty $S$-scheme $j^\inv(q\circ z)\cong\emptyset_{\!_{S}}$. Therefore, the functor $j^\inv$ is almost cocontinuous, as in Definition.\ref{Def:Almost:Cocontinuous}.
\end{proof}
This proof shows, in particular, that the functor $j^\inv:\PropSch/T\to \PropSch/S$ is essentially surjective, when $j:S\oim T$ is an open immersion in $\fdNoe$.

\Cor{}{Let $j:S\oim T$ be an open immersion in $\fdNoe$. When $\tau$ is either the \cdp-pretopology or the proper pretopology, the direct image functor $j_\ast:\scC_\tau(S)\to \scC_\tau(T)$ is a strong monoidal exact functor, \cf \cite[Prop.4.5]{GK:15:PAG}.}{Open:Direct:Image}
\begin{proof}
	A direct result of Lemma \ref{Lem:Open:Direct:Image} and Lemma \ref{Lem:Almost:Cocontinuous}.
\end{proof}

\subsubsection{Proper Base Change}
The inverse image and proper direct image functors satisfy the \indx{proper-base change property}, as in \cite[\S.1.1.9]{CD:13:TCMM}.

\Lem{}{Assume that $f$ is a proper morphism in $\fdNoe$, and let
	\[
	\begin{tikzpicture}[descr/.style={fill=white}]
	\node (X0) at (0,0) {$S'$};
	\node (XBar0) at (2,0) {$S$};
	\node (X1) at (0,-2) {$T'$};
	\node (XBar1) at (2,-2) {$T$};
	\node at (.25,-.4) {$\ulcorner$};
	\path[->,font=\scriptsize]
	(XBar0) edge node[right]{$f$} (XBar1)
	(X0) edge node[left]{$f'$} (X1)
	(X0) edge node[above]{$g'$} (XBar0)
	(X1) edge node[below]{$g$} (XBar1);
	\end{tikzpicture}
	\]
	be a Cartesian square in $\fdNoe$. Then, the exchange natural transformation $\varTheta:f'_\#\circ {g'}^{\ast}{\Rightarrow} g^{\ast}\circ f_\#:\scC^\omega_\tau(S)\to \scC^\omega_\tau(T')$, induced by the adjunctions $f'_\#\dashv {f'}^\ast$ and $f_\#\dashv f^\ast$ is an isomorphism, which induces a canonical homotopy $\varTheta:f'_\#\circ {g'}^{\ast}{\Rightarrow} g^{\ast}\circ f_\#:\Kg\big(\scC^\omega_\tau(S)\big)\to \Kg\big(\scC^\omega_\tau(T')\big)$.
}{Base:Change:Noetherian}
\begin{proof} The functor $f'_\#\circ {g'}^{\ast}$ (\resp. $g^{\ast}\circ f_\#$) is given by the $\tau$-sheafification of a left \kan extension along $(f'_{\!_\circ}\circ {g'}^\inv)^\op$ (\resp. $(g^\inv\circ f_{\!_\circ})^\op $), and the natural transformation $\varTheta$ is induced from the canonical natural transformation $f'_{\!_\circ}\circ {g'}^\inv \Rightarrow g^\inv\circ f_{\!_\circ}:\PropSch/S\to \PropSch/T'$, by the universal property of Kan extensions.
	
	\medskip
	
	For a proper $S$-scheme $p:P\to S$, considering the diagram
	\[
	\begin{tikzpicture}[descr/.style={fill=white}]
	\node (P0) at (-2,0) {$P'$};
	\node (PBar0) at (-2,-2) {$P$};
	\node (X0) at (0,0) {$S'$};
	\node (XBar0) at (0,-2) {$S$};
	\node (X1) at (2,0) {$T'$};
	\node (XBar1) at (2,-2) {$T,$};
	\node at (.25,-.4) {$\ulcorner$};
	\node at (-1.75,-.4) {$\ulcorner$};
	\path[->,font=\scriptsize]
	(XBar0) edge node[below]{$f$} (XBar1)
	(X0) edge node[above]{$f'$} (X1)
	(X0) edge node[left]{$g'$} (XBar0)
	(X1) edge node[right]{$g$} (XBar1)
	(PBar0) edge node[below]{$p$} (XBar0)
	(P0) edge node[above]{${g'}^\inv(p)$} (X0)	
	(P0) edge (PBar0)
	;
	\end{tikzpicture}
	\]
	of Cartesian squares in $\fdNoe$, one sees that the canonical morphism $f'_{\!_\circ}\big({g'}^\inv(p)\big)\to g^\inv\big(f_{\!_\circ} (p)\big)$ in $\PropSch/T'$ is an isomorphism, and hence the induced natural transformation $\varTheta:f'_\#\circ {g'}^{\ast}{\Rightarrow}g^{\ast}\circ f_\#$ is a natural isomorphism, \cf \cite[Ex.1.1.11]{CD:13:TCMM}.
\end{proof}
\Cor{}{Let $i:S\cim T$ be a closed immersion in $\fdNoe$. Then, the functor $i_\#:\scC^\omega_\tau(S)\to \scC^\omega_\tau(T)$ is fully faithful.
}{}
\begin{proof}The statement of the corollary follows from \cite[Cor.1.1.20]{CD:13:TCMM} as $\scC^\omega_\tau$ is a proper-fibred category over $\fdNoe$, by Corollary \ref{Cor:Waldhausen:Noetherian:Open}, Lemma \ref{Lem:Restricted:Pushforward:Exact}, Lemma \ref{Lem:Base:Change:Noetherian}, and \cite[\S.1]{CD:13:TCMM}.
\end{proof}

\Eg{}{
	Assume that $i:V\cim S$ is a closed immersion in $\fdNoe$ with complementary open immersion $j:U\oim S$. Then, one has the Cartesian squares
	\[
	\begin{tikzpicture}[descr/.style={fill=white}]
	\node (X0) at (0,0) {$\emptyset$};
	\node (XBar0) at (2,0) {$V$};
	\node (X1) at (0,-2) {$U$};
	\node (XBar1) at (2,-2) {$S$};
	\node at (.25,-.4) {$\ulcorner$};
	\path[right hook->,font=\scriptsize]
	(XBar0) edge node{$\diagup$} node[right]{$i$} (XBar1)
	(X0) edge node{$\diagup$}node[left]{$\emptyset_U$} (X1)
	(X0) edge node[descr]{$\!\!\!\circ\!\!\!$}node[above]{$\emptyset_V$} (XBar0)
	(X1) edge node[descr]{$\!\!\!\circ\!\!\!$}node[below]{$j$} (XBar1);
	\end{tikzpicture}\qquad
	\begin{tikzpicture}[descr/.style={fill=white}]
	\node (X0) at (0,0) {$V$};
	\node (XBar0) at (2,0) {$V$};
	\node (X1) at (0,-2) {$V$};
	\node (XBar1) at (2,-2) {$S$};
	\node at (.25,-.4) {$\ulcorner$};
	\path[right hook->,font=\scriptsize]
	(XBar0) edge node{$\diagup$} node[right]{$i$} (XBar1)
	(X0) edge node{$\diagup$}node[left]{$\id_V$} (X1)
	(X0) edge node{$\diagup$}node[above]{$\id_V$} (XBar0)
	(X1) edge node{$\diagup$}node[below]{$i$} (XBar1);
	\end{tikzpicture}
	\]
	in $\fdNoe$. Since $\scC_\tau(\emptyset)$ is isomorphic to the terminal Waldhausen category with one object and one morphism, one has a natural isomorphism $(j^{\ast}\circ i_\#)(\scX) \cong \underline{\yon}^\tau_{\emptyset_U}=\0$ for every $\scX\in \scC^\omega_\tau(V)$
	\ie $j^{\ast}\circ i_\# $ is naturally isomorphic to the zero functor $\scC^\omega_\tau(V)\to \scC^\omega_\tau(U)$. Also, the adjoint unit $\id_{\scC^\omega_\tau(V)}\Rightarrow i^{\ast}\circ i_\#$ is a natural isomorphism, and hence $i^\ast$ is essentially surjective.
}{Closed:Open:Sequence:0}

\subsubsection{Proper Projection Formula}
The inverse image and proper direct image functors also satisfy the \indx{proper-projection formula property}, as in \cite[\S.1.1.26]{CD:13:TCMM}.
\Lem{}{Assume that $f:S\to T$ is a proper morphism in $\fdNoe$. Then, $f$ satisfies the projection formula, \ie for pointed $\tau$-sheaves $\scX$ in $\scC^\omega_\tau(S)$ and $\scY$ in $\scC^\omega_\tau(T)$, the projection natural transformation $f_\#\big(\scX\wedge_{\!_S} f^{\ast}(\scY) \big) \Rightarrow f_\# (\scX)\wedge_{\!_T} \scY$, induced by the adjunction $f_\#\dashv f^\ast$, is an isomorphism in $\scC^\omega_\tau(T)$, \cf \cite[Ex.1.1.28]{CD:13:TCMM}. 
}{Noetherian:Projection:Formula}
\begin{proof} Suppose that $p:P\to S$ (\resp. $q:Q\to T$) is a proper $S$-scheme (\resp. $T$-scheme). Similar to the proof of Lemma \ref{Lem:Base:Change:Noetherian}, considering the diagram 
	\[
	\begin{tikzpicture}[descr/.style={fill=white}]
	\node (P0) at (-2,0) {$P'$};
	\node (PBar0) at (-2,-2) {$P$};
	\node (X0) at (0,0) {$Q'$};
	\node (XBar0) at (2,0) {$Q$};
	\node (X1) at (0,-2) {$S$};
	\node (XBar1) at (2,-2) {$T,$};
	\node at (.25,-.4) {$\ulcorner$};
	\node at (-1.75,-.4) {$\ulcorner$};
	\path[->,font=\scriptsize]
	(X1) edge node[below]{$f$} (XBar1)
	(X0) edge (XBar0)
	(X0) edge node[left]{$f^\inv(q)$} (X1)
	(XBar0) edge node[right]{$q$} (XBar1)
	(PBar0) edge node[below]{$p$} (X1)
	(P0) edge (X0)	
	(P0) edge (PBar0)
	;
	\end{tikzpicture}
	\]
	of Cartesian squares in $\fdNoe$, one sees that the projection morphism $f_{\!_\circ} \big(p\times_{\id_{\!_S}} f^\inv(q)\big)\to f_{\!_\circ}(p)\times_{\id_{\!_T}} q$, induced by the adjunction $f_{\!_\circ}\dashv f^\inv$, is an isomorphism in $\PropSch/T$, and hence the induced projection natural transformation $f_\#\big(\underline{\yon}^\tau_{\!_p}\wedge_{\!_S} f^{\ast}(\underline{\yon}^\tau_{\!_q}) \big)\Rightarrow
	f_\# (\underline{\yon}^\tau_{\!_p})\wedge_{\!_T} \underline{\yon}^\tau_{\!_q}$ is an isomorphism in $\scC^\omega_\tau(T)$. {Then, the statement of the proposition follows from the symmetric product $\wedge_{\!_S}$ and $\wedge_{\!_T}$ being biexact, and the functors $f_\#$ and $f^{\ast}$ being exact.}
\end{proof}
\Cor{}{Let $f:S\to T$ be a proper morphism in $\fdNoe$. Then, there exists a canonical path $f_\#\big(\underline{x} \cdot f^{\ast}(\underline{y}) \big)\to f_\# (\underline{x})\cdot \underline{y}$ in $\Kg\big(\scC^\omega_\tau(T)\big)$, for every $\underline{x}\in \Kg\big(\scC^\omega_\tau(S)\big)$ and $\underline{y}\in \Kg\big(\scC^\omega_\tau(T)\big)$. Thus, one has $f_\#\big([\yon^{\ps_\tau}_{\!_x}] \cdot f^{\ast}([\yon^{\ps_\tau}_{\!_y}]) \big)=f_\#([\yon^{\ps_\tau}_{\!_x}])\cdot [\yon^{\ps_\tau}_{\!_y}]$, for every $S$-scheme $x$ and $T$-scheme $y$.}{Noetherian:Projection:Formula}

\subsection{Counting Rational Points}
The motivic measure of counting rational points over a finite field is a shadow of a point on the \cdp-site of proper schemes over the ground field, as in Corollary \ref{Cor:Counting:Points:cdp}.

\medskip

Fix a finite field $\FF_q$ with $q$ elements. Recall that a point on the \cdp-site $(\PropSch/\FF_q,\cdp)$ is an adjunction $u^\ast:\Shv_\cdp(\PropSch/\FF_q)\rightleftarrows \Set:u_\ast$, in which $u^\ast$ is left exact, see \S.\ref{Toposes:Points}. Since both $u_\ast$ and $u^\ast$ are left exact, they both preserve the {final} objects; and hence they induce an adjunction $u^\ast_\bullet:\scC({\FF_q})\rightleftarrows \Set_\bullet:u_{\ast,\bullet}$, for having $\scC({\FF_q})\cong \ast\downarrow \Shv_\cdp(\PropSch/\FF_q)$ and $\Set_\bullet\cong \ast \downarrow\Set$. Moreover, $u^\ast_\bullet$ is also left exact, and hence it is a strong monoidal exact functor.

\medskip

Recall that if a functor $u:\PropSch/\FF_q\to \Set$ is flat and continuous with respect to the \cdp-pretopology, it defines a point $(u^\ast,u_\ast)$ in the \cdp-site, where 
\[
u^\ast=-\otimes_{\PropSch/\FF_q} u\andd u_\ast=\Hom^{\PropSch/\FF_p}(u,-)
\]
are the stalks and skyscraper functors associated to $u$, respectively, see \S.\ref{Toposes:Points}.

\Lem{}{The functor $\Gamma_\bullet:\scC({\FF_q}) \to \Set_\bullet$, induced by the global section functor $\Gamma:\Shv_\cdp(\PropSch/\FF_q)\to \Set$ is a strong monoidal exact functor. Moreover, for every $\FF_q$-scheme $X$, one has an isomorphism of pointed sets $\Gamma_\bullet(\yon^\ps_{\!_X})\cong X(\FF_q)_+$.
}{}
\begin{proof}
	The corepresentable functor $u\coloneqq \yon^{\Spec \FF_q}:\PropSch/\FF_q\to \Set$ is flat and continuous with respect to the \cdp-pretopology, as seen below.
	\begin{itemize}
		\item Since $\PropSch/\FF_q$ is Cartesian, every corepresentable functor is flat, as its category of elements is cofiltered, see \cite[\S.VII.6.Def.2 and Th.3]{MLM:92:SGL}. In particular, $u$ is flat.
		\item In the light of \S.\ref{Toposes:Points}, to show that $u$ is continuous with respect to the \cdp-pretopology, it suffices to show that $\Hom^{\PropSch/\FF_p}(u,S)$ is a \cdp-sheaf for every set $S\in \Set$ and that the sheafification morphism $\eta_\scP:\scP\to \scP^{\!^{\ash_\cdp}}$ is mapped to a bijection by the functor $-\otimes_{\PropSch/\FF_q} u$, for every presheaf $\scP\in \PSh(\PropSch/\FF_q)$. 
		\begin{itemize}
			\item For a \cdp-square \eqref{Diag:cd:Square} in $\PropSch/\FF_q$, a rational point $x:\Spec \FF_q\to X$ factorises uniquely though either $A$ or $Y$, or both (in which case it factorises uniquely through $B$). Thus, the functor $u$ maps every \cdp-square in $\PropSch/\FF_q$ to a pushout square in $\Set$, and hence {the functor} $u^\op$ maps \cdp-squares to Cartesian squares in $\Set^\op$. Since limits commute with each other and representable functors preserve limits, the presheaf $\Hom^{\PropSch/\FF_p}(u,S)$ maps \cdp-squares to pullback squares. Also, it maps the empty $\FF_q$-scheme to a terminal set. Hence, $\Hom^{\PropSch/\FF_p}(u,S)$ is a \cdp-sheaf, for every set $S\in \Set$, by Lemma \ref{Lem:cd:Sheaf:Square}.

			\item For every presheaf $\scP\in \PSh(\PropSch/\FF_q)$, one has
			\[
			\qquad\qquad\scP \otimes_{\PropSch/\FF_q} u=\int^{P\in \PropSch/\FF_q} \PSh(\PropSch/\FF_q)(\yon_{\!_P},\scP)\times u(P)
			\cong \scP(\Spec \FF_q).
			\]
			Since the \cdp-pretopology is completely decomposed, every \cdp-covering family of $\Spec \FF_{q}$ splits. That is, every \cdp-covering family $\scU=\{\sigma_\alpha:P_\alpha\to \Spec \FF_q \mid \alpha\in A\}$ in $\PropSch/S$ admits a refinement $\scV=\{\delta_\beta:Q_\beta\to \Spec \FF_q \mid \beta\in B\}$ such that $Q_{\beta_0}=\Spec \FF_q$ and $\delta_{\beta_0}=\id_{\Spec \FF_q}$ for some $\beta_0\in B$. Therefore, $\eta_{\scP,\Spec \FF_q}:\scP(\Spec \FF_q)\to \scP^{\!^{\ash_\cdp}}(\Spec \FF_q)$ is a bijection, and so is $\eta_\scP \otimes_{\PropSch/\FF_q} u$.
		\end{itemize}
	\end{itemize} 
	Then, there exists a \cdp-point $u^\ast:\Shv_\cdp(\PropSch/\FF_q)\rightleftarrows \Set:u_\ast$, with the stalks and skyscrapers functors $u^\ast=-\otimes_{\PropSch/\FF_q} u$ and $u_\ast=\Hom^{\PropSch/\FF_p}(u,-)$. In particular, for a \cdp-sheaf $\scX$ on $\PropSch/\FF_q$, one has
	\[
	u^\ast (\scX)\cong\scX(\Spec \FF_q)\cong \Shv_\cdp(\PropSch/\FF_q)(\yon_{\!_{\Spec \FF_q}}^{\!^{\ash_\cdp}},\scX)\cong \Shv_\cdp(\ast,\scX)=\Gamma(\scX).
	\]
	Therefore, the induced functor $\Gamma_\bullet:\scC({\FF_q})\to \Set_\bullet$
	is a strong monoidal exact functor of Waldhausen categories, for being a left exact left adjoint functor.
	
	\medskip
	
	For an $\FF_q$-scheme $X$, let $i:Z\cim P$ be a compactification of $X$ in $\sftSch/\FF_q$. Since every \cdp-covering family of $\Spec \FF_{q}$ splits, one has
	\begin{equation}\label{Gamma:cdp:Point:Counting}
	\begin{split}
	\Gamma_\bullet(\yon^\ps_{\!_X})&
	\cong \coker\big(\Gamma_\bullet(\underline{\yon}_{\!_Z})\rightarrowtail \Gamma_\bullet(\underline{\yon}_{\!_P})\big)
	\cong \coker \big(\yon_{\!_{Z,+}}^{\!^{\ash_\cdp}}(\FF_q) \rightarrowtail \yon_{\!_{P,+}}^{\!^{\ash_\cdp}}(\FF_q) \big)\\&
	\cong \coker \big(\yon_{\!_{Z,+}}(\FF_q) \rightarrowtail \yon_{\!_{P,+}}(\FF_q) \big)
	\cong \coker \big(Z(\FF_q)_+ \rightarrowtail P(\FF_q)_+ \big)
	\cong X(\FF_q)_+.
	\end{split}
	\end{equation}
\end{proof}

\Lem{}{The strong monoidal exact functor $\Gamma_\bullet:\scC({\FF_q})\to \Set_\bullet$ restricts to a strong monoidal exact functor $\Gamma^\omega_\bullet:\scC^\omega({\FF_q})\to \FSet_\bullet$.
}{}
\begin{proof}Since $\scC^\omega({\FF_q})$ is a full symmetric monoidal \waldhausen subcategory in $\scC({\FF_q})$, the statement of the lemma follows from the existence of the restriction $\Gamma^\omega_\bullet$.
	
	\medskip
	Recall Example \ref{Ex:Counting:Rational:Points}, since $\Gamma_\bullet$ commutes with pushouts, it suffices to show that $\Gamma_\bullet$ restricts to a functor $\Gamma^0_\bullet:\scC^0({\FF_q})\to \FSet_\bullet$, in order to induce a strong monoidal exact functor $\Gamma^\omega_\bullet:\scC^\omega({\FF_q})\to \FSet_\bullet$, which is a result of \eqref{Gamma:cdp:Point:Counting}. 
\end{proof}
\Cor{}{The functor $\Gamma^\omega_\bullet$ induces a morphism of commutative ring spectra
	\[
	\Gamma_\bullet:\Kg\big(\scC^\omega({\FF_q})\big)\to \SS,
	\]
	that sends a point in the component $[\yon^\ps_{\!_X}]$ to a point in the component $[X(\FF_q)_+]$, for every $\FF_q$-scheme $X$. Hence, the motivic measure \eqref{K:Schemes} factorises the classical motivic measure of counting rational points over $\FF_q$ as $\mu_{\#}=\pi_0(\Gamma_\bullet)\circ \mu_{\!_\cdp}$.
}{Counting:Points:cdp}
The homomorphism $\pi_n(\Gamma_\bullet)$, for $n\geq 1$, might be thought of as higher point counting measures.
\Rem{}{This argument applies for any base scheme $S$ in $\fdNoe$ and for any $\tau$-point $(u^\ast,u_\ast)$ on $\PropSch/S$, whenever the stalks functor $u^\ast_\bullet$ restricts to a functor $\scC^\omega_\tau(S)\to \FSet_\bullet$.
}{}
\Cor{{\cite[Prop.5.21]{Cam:17:KTSV}}}{The composition $\Gamma^\omega_\bullet\circ \upsilon_{\FF_q}:\FSet_\bullet\to \FSet_\bullet$ is an exact equivalence of \waldhausen categories, where $\upsilon_{\FF_q}$ is the exact functor $\upsilon_{\scC^\omega({\FF_q})}$ in \eqref{Unit:Waldhausen:Upsilon}. Thus, the map of spectra $\Kg(\Gamma^\omega_\bullet\circ \upsilon_{\FF_q}):\SS\to \SS$ is a homotopy equivalence, and hence the spectrum $\Kg\big(\scC^\omega({\FF_q})\big)$ splits through $\SS$. That is, $\Kg\big(\scC^\omega({\FF_q})\big)\cong \SS \vee \tilde{\Kg}\big(\scC^\omega({\FF_q})\big)$, where $\tilde{\Kg}\big(\scC^\omega({\FF_q})\big)$ is a cofibre of $\Kg(\upsilon_{\FF_q})$.
}{}

\appendix

\section{Grothendieck Sites}\label{Cat:GroSite}A \grothendieck topology is a generalisation of {topological coverings} to abstract categories, which enables the development of cohomology theories on abstract categories. 

\medskip

In this section, we recall some results on Grothendieck sites used in the paper. Interested readers may consult \cite{Joh:02:SE}, \cite{MLM:92:SGL}, and \cite{KS:06:CS}.

\medskip

Throughout this section, let $\scC$ be a locally small category. A \indx{sieve} $S$ on an object $U\in \scC$ is an inclusion $S\subset \yon_{\!_U}:\scC^\op\to \Set$. For a set $\scU$ of morphisms with a common codomain in $\scC$. For a sieve $S$ on $U\in \scC$ and for a morphism $\varphi:V\rightarrow U$ in $\scC$, there exists a sieve $\varphi^{\ast}S$ on $V$, given on an object $W\in \scC$ by $\varphi^{\ast}S(W)\coloneqq\{\psi:W\to V \text{ in }\scC\mid \varphi\circ \psi \in S(W)\}$, and it is called the restriction of $S$ along $\varphi$.

\medskip

A ({\em \grothendieck}) {\index{topology@{Grothendieck topology}}\em topology} $\tau$ on the category $\scC$ is a set $\tau=\{\Cov_\tau(U)\mid U\in \scC\}$, in which $\Cov_\tau(U)$ is a set of sieves on $U$, for every object $U\in \scC$, subject to the maximal sieve, stability, and transitivity\footnote{The transitivity axiom is also referred to by the \indx{local character} axiom, as in \cite[\S.C.2.1.p.541]{Joh:02:SE}.} axioms, as in \cite[\S.III.2.Def.1]{MLM:92:SGL}. Then, the pair $(\scC,\tau)$ is called a {\index{topology@{Grothendieck topology}!{Grothendieck site}}({\em\grothendieck}) \em site}, and a sieve $S$ on $U\in \scC$ is called a $\tau$-\indx{covering sieves} if it belongs to $\Cov_\tau(U)$. Moreover, when $\scC$ is essentially small, the site  $(\scC,\tau)$ is said to be \indx{essentially small}.

\Eg{}{
	A sieve is said to be \indx{effective epimorphic} if it forms a colimit cocone, and it is said to be \indx{universally effective epimorphic} if all its restrictions are \indx{effective epimorphic}. Every category $\scC$ admits a topology whose covering sieves are universally effective epimorphic sieves, called the {\index{topology@{Grothendieck topology}!canonical --}\em canonical topology} on $\scC$, see \cite[p.542-543]{Joh:02:SE}. A topology that is contained in the canonical topology is said to be \indx{subcanonical}.}{}

The intersection of topologies on $\scC$ is a topology on $\scC$. Hence, given a set $\scS$ of sets of sieves on $\scC$, the intersection of all topologies on $\scC$ that contain $\scS$ is a topology on $\scC$, called the \indx{topology generated by $\scS$}. In some occasions, it may be simpler to specify topologies in terms of sets of morphisms that generate them, as in \S.\ref{SubSec:GTAG}. 

\medskip

In particular, when $\scC$ has fibre products, it is usually convenient to specify a topology by a \grothendieck pretopology that generates it. A {\index{topology@{Grothendieck topology}!Grothendieck pretopology}\em \grothendieck pretopology} on $\scC$ is a set $\tau=\{\Cov_\tau(U)\mid U\in \scC\}$, in which $\Cov_\tau(U)$ is a set of families of morphisms with a common codomain that satisfies closure conditions similar to those of a topology, see \cite[\S.III.2.Def.2]{MLM:92:SGL}. A family $\scU$ in $\Cov_\tau(U)$ is called a $\tau$-\indx{covering family} of $U$. Also, the unique element of a singleton $\tau$-covering family of $U\in \scC$ is called a $\tau$-\indx{cover} of $U$.

\medskip

A \indx{refinement} of a family of morphisms $\scU=\{\sigma_\alpha:U_\alpha\to U\mid \alpha\in A\}$ is a pair $(f,\scU')$ of a map $f:A'\to A$ and a family of morphisms $\scU'=\{\sigma'_{\alpha'}:U'_{\alpha'}\to U\mid \alpha'\in A'\}$ such that $\sigma'_{\alpha'}$ factorises through $\sigma_{f(\alpha)}$, for every $\alpha'\in A'$.

\Def{}{Assume that $\scC$ has fibre products. A pretopology $\tau$ on $\scC$ is said to be \indx{saturated} if every family of morphisms in $\scC$ with a common codomain that admits a refinement by a $\tau$-covering family is a $\tau$-covering family. Assume that $\scC$ admits finite coproducts, the pretopology $\tau$ is said to be \indx{additively-saturated} if for every $\tau$-covering family $\scU=\{\sigma_\alpha:U_\alpha\to U\mid \alpha\in A\}$, the set $A$ is finite and the singleton 
$
	\{\coprod_{\alpha\in A}\sigma_\alpha:\coprod_{\alpha\in A}U_\alpha\to U\}
$
	is a $\tau$-covering family.
}{}

Every pretopology admits a \indx{saturation}, that is a pretopology in which covering families are precisely families that admit refinements in the given pretopology, see \cite[Def.2.52 and Prop.2.53]{Vis:08:NGTFCDT}. Also, \indx{additive-saturations} are defined similarly.

\medskip

In practice, the pretopologies one considers are almost never saturated, and their saturations allow redundant covering families that does not reflect the intended properties of the topology. For example, for a pretopology $\tau$ on $\scC$ and a $\tau$-covering family $\scU=\{\sigma_\alpha:U_\alpha\to U\mid \alpha\in A\}$, the set $\scU\bigsqcup \{f\}$ is a covering family in the saturation of $\tau$, for any morphism $f:V\to U$ in $\scC$. It is often more practical to consider additively-saturated pretopologies, as in Example \ref{Ex:Cover:Saturated:Not:Cover}.

\medskip

Every pretopology generates a topology, and different pretopologies may define the same topology. In particular, a pretopology and its (additive-)saturation define the same topology. For a pretopology $\tau$, we may abuse notation and use $\tau$ to also denote the topology it defines.

\subsection{The Category of Sheaves}\label{Sheaves:Cat}
A functor $P:\scC^\op\to \Set$ is called a \indx{presheaf} {\em of }({\em small}){ sets} on $\scC$, the functor category $\Funct(\scC^\op,\Set)$ is called the \indx{category of presheaves} on $\scC$, and it is usually denoted by $\PSh(\scC)$. For a topology $\tau$ on $\scC$, a presheaf $P:\scC^{op}\rightarrow \Set$ is said to be a \indx{$\tau$-separated presheaf}, a \indx{$\tau$-weak sheaf}, or a \indx{$\tau$-sheaf} if for every object $U\in \scC$ and for every $\tau$-covering sieve $S$ on $U$, the canonical map
\begin{equation}\label{Eq:Sheaf:Descent}
\iota_S^\ast: \PSh(\scC)\big(\yon_{\!_U},P\big)\to \PSh(\scC)\big(S,P\big),
\end{equation}
induced by the inclusion $\iota_S:S\incl \yon_{\!_U}$, is injective, surjective, or bijective, respectively. The full subcategory in $\PSh(\scC)$ of $\tau$-sheaves on $\scC$ is usually denoted by $\Shv_\tau(\scC)$.

\medskip

When $\scC$ is essentially small, the subcategory of $\tau$-sheaves on $\scC$ is a Cartesian reflective subcategory in the category of presheaves on $\scC$, \ie the inclusion $\Shv_\tau(\scC)\incl \PSh(\scC)$ admits a left adjoint\footnote{A left adjoint to an inclusion functor is called a \indx{reflector}.} which is left exact, called the $\tau$-\indx{sheafification functor} or the \indx{associated $\tau$-sheaf functor}, and {it is} denoted by $-^{\ash_\tau}$, see \cite[Th.3.3.12]{Bor:94:HCA3}. For a presheaf $P$ on $\scC$, there exists a $\tau$-separated presheaf $P^{+_\tau}$, given for an object $U\in \scC$ by the filtered colimit
\[
P^{+_\tau}(U)\coloneqq \colim\ \PSh(\scC)\big(i_U(-),P\big),
\]
where $i_U:\Cov_\tau(U)\incl \PSh(\scC)$ is the canonical such inclusion, and $\Cov_\tau(U)$ is consider as a preordered category. Then, the associated $\tau$-sheaf $P^{\ash_\tau}$ can be given by $P^{\ash_\tau}\coloneqq P^{+_\tau +_\tau}$, see \cite[\S.3.3]{Bor:94:HCA3}. Alternatively, to a presheaf $P$ one defines a $\tau$-separated presheaf $P^{\asp_\tau}$, given on an object $U\in \scC$ by the quotient $P^{\asp_\tau}(U)\coloneqq P(U)/\sim$, where $\sim$ is the relation on $P(U)$, with $p\sim p' $ for $p,p'\in P(U)$ \iff there exists a $\tau$-covering sieve on which the restrictions of $p$ and $p'$ coincide, and given on morphisms by the universal property of quotients. Then, the associated $\tau$-sheaf $P^{\ash_\tau}$ can be given by $P^{\ash_\tau}\coloneqq {(P^{\asp_\tau})}^{+_\tau}$, see the proof of \cite[Th.2.64.(ii)]{Vis:08:NGTFCDT}. For a section $p\in P(U)$, we denote its image in $P^{\asp_\tau}(U)$ (\resp. $P^{\ash_\tau}(U)$) by $p^{\asp}$ (\resp. $p^{\ash}$).

\medskip
The subcategory of $\tau$-sheaves on $\scC$ is closed under limits in the category $\PSh(\scC)$, and hence it is a complete category, whose limits are given object-wise. In fact, when $\scC$ is essentially small, the category of $\tau$-sheaves on $\scC$ is a bicomplete category, whose colimits are given by the $\tau$-sheafification of colimits in the category of presheaves, that the $\tau$-sheafification functor preserves colimits.

\medskip

\Rem{}{
	Giving a topology on an essentially small category $\scC$ is a way of formally declaring specific cocones to be colimit cocones in the resulting category of sheaves. Recall that the category $\PSh(\scC)$ is the free cocompletion for $\scC$, and hence it formally adds all small colimits, forgetting the colimits that already exist in $\scC$. The sheaf condition \eqref{Eq:Sheaf:Descent} shows that the cocone of a covering sieve is mapped into a colimit cocone in the category of sheaves. 
	{That is}, for a topology $\tau$ on $\scC$, for a $\tau$-sheaf $Q$ on $\scC$, and for a $\tau$-covering sieve $S$ on $U\in \scC$, one has
	\begin{align*}
	\Shv_\tau(\scC)\big(\colim (\yon_-^{\ash_\tau}\pi_S),Q\big)&
	\cong \PSh(\scC)\big(\colim ({\yon_-\pi_S)},Q\big)
	\cong \PSh(\scC)\big(S,Q\big)\\&
	\cong \PSh(\scC)\big(\yon_{\!_U},Q\big)
	\cong \Shv_\tau(\scC)\big(\yon_{\!_U}^{\ash_\tau},Q\big),
	\end{align*}
	where $\pi_S:\El(S)\to \scC$ is the canonical projection functor from the category of elements $\El(S)$, and hence $\colim (\yon_-^{\ash_\tau}\pi_S)$ is isomorphic to $\yon_{\!_U}^{\ash_\tau}$ in $\Shv_\tau(\scC)$. In particular, the canonical topology is the coarsest topology that retrieves universal colimits cocones that exist in $\scC$.}{Sheafification:Cocones}

\subsubsection{Local Epimorphisms, Monomorphisms, and Isomorphisms}\label{Sub:Local:Epi:Mono:Iso}
For a topology $\tau$ on an essentially small category $\scC$, a morphism of presheaves is said to be a {\index{Local!-- epimorphisms}\em$\tau$-local epimorphism} (\resp. {\index{Local!-- monomorphisms}\em$\tau$-local monomorphism}, \resp. {\index{Local!-- isomorphisms}\em$\tau$-local isomorphism}) if its $\tau$-sheafification is an epimorphism (\resp. a monomorphism, \resp. an isomorphism). In particular, epimorphisms (\resp. monomorphisms) are $\tau$-local epimorphisms (\resp. $\tau$-local monomorphisms), that the $\tau$-sheafification functor preserves epimorphisms for being a reflector and preserves monomorphisms for being left exact. Moreover, the category $\Shv_\tau(\scC)$ is a reflective localisation of $\PSh(\scC)$ with respect to $\tau$-local isomorphisms. 
\Eg{}{The component of the $\tau$-sheafification adjunction unit is given for a presheaf $P$ by a morphism $\eta^\tau_P: P\to P^{\ash_\tau}$ in $\PSh(\scC)$ for which there exists a bijection
	\[
	{\eta_P^\tau}^\ast:\PSh(\scC)(P^{\ash_\tau},Q)\to \PSh(\scC)(P,Q),
	\]
	for every $\tau$-sheaf $Q\in \Shv_\tau(\scC)\subset \PSh(\scC)$, and hence $\eta_P^\tau$ is a $\tau$-local isomorphism.
}{}

\Eg{}{The set of $\tau$-local epimorphisms on $\scC$ determines the topology $\tau$ on $\scC$. The {\index{topology@{Grothendieck topology}!discrete --}\em discrete} (\resp. {\index{topology@{Grothendieck topology}!indiscrete --}\em indiscrete}) {\em topology} on an essentially small category $\scC$ is the topology $\tau$ on $\scC$ whose $\tau$-local epimorphisms are all morphisms (\resp. all epimorphisms), see \cite[Ex.16.1.9]{KS:06:CS}.
}{}

While morphisms of $\tau$-sheaves that are object-wise surjective are epimorphisms of $\tau$-sheaves, the inverse does not hold as one may deduce from the following lemma.

\Lem{}{Let $\scC$ be an essentially small category with pullbacks, let $\tau$ be a pretopology on $\scC$, and let $f:P\to Q$ be a morphism of presheaves on $\scC$. Then, $f$ is
	\begin{itemize}
		\item a $\tau$-local epimorphism \iff for every object $U\in \scC$ and for every section $q\in Q(U)$, there exists a ${\tau}$-covering family $\scU=\{\sigma_\alpha: U_\alpha\to U\mid \alpha\in A\}$ and a section $p_\alpha\in P(U_\alpha)$ such that $\sigma_\alpha^\ast (q)=f_{U_\alpha}(p_\alpha)$, for every $\alpha\in A$; and
		\item a $\tau$-local monomorphism \iff for every object $U\in \scC$ and for every pair of sections $p,p'\in P(U)$ for which $f_U(p)=f_U(p')\in Q(U)$ there exists a ${\tau}$-covering family $\scU=\{\sigma_\alpha: U_\alpha\to U\mid \alpha\in A\}$ such that $\sigma_\alpha^\ast (p)=\sigma_\alpha^\ast ( p')$, for every $\alpha\in A$. 	
	\end{itemize}
}{Local_Morphism}
\begin{proof}
	See \cite[Lem.3.16]{Jar:15:LHT}.
\end{proof}

\subsection{Continuous Functors}\label{Continuous:Maps}
Let $(\scC,\tau)$ and $(\scD,\varsigma)$ be essentially small sites. A functor $f^\inv :\scD\to \scC$ induces a functor $f_\ast:\PSh(\scC)\to \PSh(\scD)$, given by precomposition with $(f^\inv)^\op$, which is called the \indx{direct image} functor along $f^\inv$. The functor $f_\ast$ admits a left adjoint $f^\ast_\pre:\PSh(\scD)\to \PSh(\scD)$, given by a left \kan extension $\Lan_{(f^\inv)^\op}$. The functor $f^\inv:\scD\to \scC$ is said to be {\index{Functor!continuous}\em continuous}, with respect to the topologies $\tau$ and $\varsigma$, if $f_\ast$ sends $\tau$-sheaves to $\varsigma$-sheaves. For a continuous functor $f^\inv$, there exists an adjunction $f^\ast:\Shv_\varsigma(\scD)\rightleftarrows\Shv_\tau(\scC):f_\ast$, where $f^\ast$ is given by the composition of $f^\ast_\pre$ with the associated $\tau$-sheaf functor $-^{\ash_\tau}$, and it is called the \indx{inverse image} functor along $f^\inv$. 

\medskip

Since the category of sheaves on an essentially small site is a reflective localisation of the category of presheaves with a left exact reflector, the functor $f^\inv$ is continuous \iff $f^\ast_\pre$ preserves local isomorphisms. Also, the functor $f^\inv$ is continuous \iff the sieve generated by $f^\inv(S)$ is a $\tau$-sieve in $\scC$, for every $\varsigma$-sieve $S$ in $\scD$, see \cite[\S.C.2.3]{Joh:02:SE}. In particular, when $f^\inv$ is a Cartesian functor between Cartesian categories and the topologies $\tau$ and $\varsigma$ are defined by pretopologies, the functor $f^\inv$ is continuous if it preserves covering families, see \cite[\S.III.Prop.1.6]{SGA4}.

\medskip

In addition to the notion of continuous functors, we need to recall the notion of almost cocontinuous functors, which admits a well-behaved direct image.

\Def{{\cite[Tag 04B7]{stacks-project}}}{Let $\scC$ and $\scD$ be essentially small categories with pullbacks, and let $\tau$ and $\varsigma$ be pretopologies on $\scC$ and $\scD$, respectively. A functor $f^\inv :\scD\to \scC$ is said to be \indx{almost cocontinuous} if for every object $V\in \scD$ and for every $\tau$-covering family $\scU=\{\sigma_\alpha:U_\alpha\to f^\inv(V)\mid \alpha\in A\}$ there exists a $\varsigma$-covering family $\scV=\{\delta_\beta:V_\beta\to V\mid \beta\in B\}$ such that for every $\beta\in B$ either 
	\begin{enumerate}
		\item the morphisms $f^\inv(\delta_\beta)$ factorises through $\sigma_\alpha$, for some $\alpha\in A$; or
		\item the empty sieve is a $\tau$-covering sieve on $f^\inv (V_\beta)$.
	\end{enumerate}
}{Almost:Cocontinuous}
\Lem{}{Assume that $\scC$ and $\scD$ are essentially small categories, let $\tau$ and $\varsigma$ be pretopologies on $\scC$ and $\scD$, respectively, and let $f^\inv :\scD\to \scC$ be a continuous and an almost cocontinuous functor,with respect to the pretopologies $\tau$ and $\varsigma$. Then, the direct image functor $	f_\ast:\Shv_\tau(\scC)\to \Shv_\varsigma(\scD)$	commutes with pushouts.
}{Almost:Cocontinuous}
\begin{proof}
	See \cite[Tag 04B9]{stacks-project}.
\end{proof}

\subsubsection{Points of Sites}\label{Toposes:Points}A \indx{point $p$ of a site} $(\scC,\tau)$ is an adjunction $p^\ast:\Shv_\tau(\scC)\rightleftarrows  \Set:p_\ast$, in which $p^\ast$ is left exact. The inverse image functor $p^\ast$ is called the \indx{stalks functor} at $p$, whereas the direct image functor $p_\ast$ is called the \indx{$\tau$-skyscraper sheaf functor} at $p$.

\medskip

Points of an essentially small site $(\scC,\tau)$ correspond to flat functors $\scC\to \Set$ that are continuous, with respect to $\tau$ and the canonical topology on $\Set$, see \cite[\S{VII.5.Cor.4}]{MLM:92:SGL}. Recall that a functor $u:\scC\to \Set$ is flat \iff its category of elements $\El(u)$ is cofiltered, see \cite[\S.VII.6.Th.3]{MLM:92:SGL}. On the other hand, since the category $\Set$ is cocomplete, every functor $u:\scC\to \Set$ induces the $u$-tensor-$\iHom$ adjunction
\[
-\otimes_\scC u: \PSh(\scC) \rightleftarrows\Set: \Hom^\scC(u,-),
\]
given by left Kan extensions of $u$ and the \yoneda embedding along each other, see \cite[\S.2]{Kan:58:FICSSC}. Since there exists a canonical equivalence of categories $\Set\cong \Shv_\can(\Set)$, a functor $u:\scC\to \Set$ is continuous, with respect to $\tau$ and the canonical topology on $\Set$, if the geometric morphism $-\otimes_\scC u:\PSh(\scC)\rightleftarrows \Set:\Hom^\scC(u,-)$ factorises through the $\tau$-sheafification geometric morphism $\PSh(\scC)\rightleftarrows\Shv_\tau(\scC)$, \ie if $\Hom^\scC(u,S)$ is a $\tau$-sheaf for every set $S\in \Set$ and the $\tau$-sheafification morphism $\eta_P:P\to P^{\ash_\tau}$ is mapped to an isomorphism by $-\otimes_\scC u$, for every presheaf $P\in \PSh(\scC)$, see \cite[\S.VII.5.Lem.3]{MLM:92:SGL} and \cite[\S.C.Lem.2.3.8]{Joh:02:SE}.

\subsection{Grothendieck Topologies in Algebraic Geometry}\label{SubSec:GTAG}
We conclude this section by recalling the topologies used in this paper and some results on their representative sheaves. For an elaborative treatment of topologies used in algebraic geometry, see \cite{GK:15:PAG}.

\medskip

A finite family of morphisms $\scU=\{\sigma_\alpha:U_\alpha\to U\mid \alpha\in A\}$ in $\fdNoe$ is said to be a {\index{proper covering}\em proper covering family} of $U$, if
\begin{itemize}
	\item the morphism $\sigma_\alpha$ is proper for every $\alpha\in A$; and
	\item $\scU$ is jointly surjective, \ie the underlying map of the coproduct morphism $\coprod_{\alpha\in A}\sigma_\alpha:\coprod_{\alpha\in A}U_\alpha\to U$ is a surjection of sets.
\end{itemize}
A proper covering family $\{\sigma_\alpha:U_\alpha\to U\mid \alpha\in A\}$ is said to be a {\index{\cdp-covering}\em \cdp-covering family} if it is \indx{completely decomposed}, \ie for every $u\in U$ there exist $\alpha\in A$ and $u_{\alpha}\in U_{\alpha}$ such that $\sigma_{\alpha}(u_{\alpha})=u$ and the induced morphism of residue fields $\rf(u)\to \rf(u_{\alpha})$ is an isomorphism\footnote{This condition is also referred to by the \indx{\nisnevich condition}, as it first appeared in \cite{Nis:89:CDTSADSSAKT}.}. The {\index{proper topology}\em proper} (\resp. {\index{\cdp-topology}\em \cdp\footnote{The notation \cdp\ appears in \cite{GK:15:PAG}, but some authors use $pro\ cdh$ instead; others use $abs$ in reference to abstract blow up squares.}}) {\em pretopology} on the category $\fdNoe$ is the pretopology whose covering families are
\begin{itemize}
	\item the proper (\resp. \cdp) covering families in $\fdNoe$; and
	\item the empty covering family of the empty scheme.
\end{itemize}
In particular, the  \cdp-pretopology is coarser than the  proper pretopology, and finer than the closed pretopology in which nonempty covering families consist of  closed immersions.

\Rem{}{Let $\scU\coloneqq\{\sigma_\alpha:U_\alpha\to U\mid \alpha\in A\}$ be a proper covering family in $\fdNoe$. Then, $\scU$ is a \cdp-covering family, \iff, for every field $\fk$, the map
	\[
	\big(\coprod_{\alpha\in A}\sigma_\alpha\big)_\ast:\fdNoe \big(\Spec \fk, \coprod_{\alpha\in A}U_\alpha\big)\to \fdNoe\big(\Spec \fk, U\big)
	\]
	is surjective. Thus, the \cdp-topology coincides with the \indx{envelop topology}, used in \cite{GS:09:MWCAV}.
}{cdp:Lifting}

For every scheme $S\in \fdNoe$, the proper, the closed, and the \cdp-pretopologies restrict to the categories $\sftSch/S$ and $\PropSch/S$. The closed pretopology is not subcanonical on $\sftSch/S$, and hence the \cdp-pretopology, and the proper pretopology are not subcanonical. That is, a surjective closed immersion $i:z\cim x$ is a closed cover of $x\in \sftSch/S$. However, $i^\ast:\yon_{\!_z}(x)\to \yon_{\!_z}(z)$ is not always a bijection. For example, let $S=\Spec \fk$ for a field $\fk$, and let $i$ be the surjective closed immersion $\Spec \fk[t]/(t^2)\cim \Spec \fk[t]/(t^3)$ in $\sftSch/S$.

\subsubsection{Completely Decomposed Structures}

\Def{}{Suppose that $x$ is an $S$-scheme. A Cartesian square 
	\begin{equation}\label{Diag:cd:Square}
	\begin{tikzpicture}[descr/.style={fill=white},ampersand replacement=\&]
	\node (m-1-1) at (0,0) {$b$};
	\node (m-1-2) at (2,0) {$y$};
	\node (m-2-1) at (0,-2) {$a$};
	\node (m-2-2) at (2,-2) {$x$};
	\node at (.25,-.4) {$\ulcorner$};
	\path[->,font=\scriptsize]
	(m-1-1) edge node[left]{$f'$} (m-2-1)
	(m-1-2) edge node[right]{$f$} (m-2-2);
	\path[right hook->,font=\scriptsize]
	(m-2-1) edge node{$\diagup$}node[below]{$i$} (m-2-2)
	(m-1-1) edge node{$\diagup$}node[above]{$i'$} (m-1-2);	
	\end{tikzpicture}
	\end{equation}
	in $\sftSch/S$ is called a \indx{\cdp-square}\footnote{Some sources refer to the square \eqref{Diag:cd:Square} by an \indx{abstract blow up} square, as in \cite[Def.12.21]{MVW:06:LNMC}; others reserve the term for a square that satisfies additional properties, as in \cite[Def.2.2.4]{SV:00:RCCS}.} over $x$ if $f$ is a proper morphism and $i$ is a closed immersion, such that the base change $\big(y\setminus i'\big) \to \big(x\setminus i\big)$ is an isomorphism. 
}{cd:topology}

\Lem{}{The \cdp-topology on the category of $S$-schemes coincides with the \grothendieck topology generated by the covering families
	\begin{itemize}
		\item $\{f:y\to x,i:a\cim x\}$, for every pair of morphisms $f:y\to x$ and $i:a\cim x$ that fit into a \cdp-square in $\sftSch/S$; and 
		\item the empty covering family of the empty $S$-scheme.
	\end{itemize}
	Equivalently, a presheaf of sets $P\in \PSh(\sftSch/S)$ is a \cdp-sheaf \iff
	\begin{itemize}
		\item $P$ sends every \cdp-square to a Cartesian square; and
		\item $P$ sends the empty $S$-scheme to a terminal set.
	\end{itemize}
}{cd:Sheaf:Square}
\begin{proof}
	See \cite[Cor.2.17]{Voe:10:HTSSCDT} and \cite[Th.2.2]{Voe:10:UMHCNcdhT}.
\end{proof}

\Prop{}{The \cdp-sheafification of the Yoneda embedding takes every \cdp-square of $S$-schemes to a cocartesian square of \cdp-sheaves on $\sftSch/S$.
}{cd:Cosheafification}
\begin{proof}
	See \cite[Lem.2.11 and Cor.2.16]{Voe:10:HTSSCDT} and \cite[Th.2.2]{Voe:10:UMHCNcdhT}.
\end{proof}

The analogue of Lemma \ref{Lem:cd:Sheaf:Square} and Proposition \ref{Prop:cd:Cosheafification} holds on $\PropSch/S$, see \cite[Lem.2.3]{Voe:10:UMHCNcdhT}.

\subsubsection{Representable Sheaves}
For canonical sites, representable presheaves are sheaves, and morphisms between representable sheaves correspond to morphisms in the original category. Since the pretopologies we are interested in here are not subcanonical, we devote this section to understanding morphisms between their representable sheaves on the categories $\sftSch/S$ and $\PropSch/S$. The argument below essentially follows \cite[\S.3.2]{Voe:96:HS}.

\Rem{}{For every (proper) $S$-scheme $p$, the representable presheaf $\yon_{\!_p}$ is additive, by the very definition of colimits, \ie for (proper) $S$-schemes $z$ and $w$, the canonical morphism $\yon_{\!_p}(z\coprod w)\to \yon_{\!_p}(z)\times\yon_{\!_p}(w)$ is an isomorphism. Also, the $\tau$-sheaf $\yon_{\!_p}^{\!^{\ash_\tau}}$ is additive for every pretopology $\tau$ on (proper) $S$-schemes that is finer than the closed pretopology. For an additively-saturated pretopology $\tau$ on the category of (proper) $S$-schemes that is finer than the closed pretopology, and for a $\tau$-covering family $\scU=\{\sigma_\alpha:z_\alpha\to z\mid \alpha\in A\}$ of (proper) $S$-schemes, one has a $\tau$-covering family 
	\[
	\scU'\coloneqq\bigg\{\coprod_{\alpha\in A} \sigma_\alpha:\coprod_{\alpha\in A} z_\alpha\to z \bigg\}
	\]
	of (proper) $S$-schemes. The additivity of $\yon_{\!_p}$ and $\yon_{\!_p}^{\!^{\ash_\tau}}$ implies that sections of $\yon_{\!_p}$ and $\yon_{\!_p}^{\!^{\ash_\tau}}$ on $\scU$ correspond to their sections on $\scU'$. Thus, without loss of generality, when considering the presheaves $\yon_{\!_p}$ and $\yon_{\!_p}^{\!^{\ash_\tau}}$, one may assume the involved $\tau$-covering families are singletons.
}{}

\Eg{}{Additively-saturated pretopologies on the category of (proper) $S$-schemes, that are finer than the closed pretopology, include:
	\begin{enumerate}
		\item the {proper} pretopology, see \cite[Tags 01T1, 01KH, and 0BX5]{stacks-project};
		\item the \cdp-pretopology, see Remark \ref{Rem:cdp:Lifting}; and
		\item the \indx{finite pretopology} (\resp. \indx{\cdf-pretopology}), which is coarser than the proper pretopology (\resp. \cdp-pretopology), whose nonempty covering families consist of finite morphisms, see \cite[Tag 0CYI]{stacks-project}.
	\end{enumerate}
	While the {proper} pretopology (\resp. \cdp-pretopology) on the category $\PropSch/S$ is saturated, as morphisms between proper $S$-schemes are {proper}, its counterpart on the category $\sftSch/S$ is not saturated. Also, the finite pretopology, and the \cdf-pretopology are not saturated on the category of (proper) $S$-schemes.
}{Cover:Saturated:Not:Cover}

\Lem{{\cite[Lem.3.2.2]{Voe:96:HS}}}{
	Assume that $\tau$ is an additively-saturated pretopology on (proper) $S$-schemes that is finer than the closed pretopology, such that $\tau$-covers are surjective, and let $p$ and $q$ be $S$-schemes, such that $p$ is reduced. Then,
	\begin{itemize}
		\item the canonical map $\sftSch/S(p,q)\to \Shv_{\tau}(\sftSch/S)(\yon_{\!_p}^{\!^{\ash_\tau}},\yon_{\!_q}^{\!^{\ash_\tau}})$ is an injection; and
		\item when $p$ and $q$ are proper $S$-schemes, the canonical map 
		\[
		\PropSch/S(p,q)	\to \Shv_{\tau}(\PropSch/S)(\yon_{\!_p}^{\!^{\ash_\tau}},\yon_{\!_q}^{\!^{\ash_\tau}})
		\]
		is an injection.
	\end{itemize}
}{Voevodsky:96:HS:Lem.3.2.2}
\begin{proof}
	Assume that $f_0,f_1:p\to q$ are morphisms of (proper) $S$-schemes such that $f_{0,\ast}=f_{1,\ast}:\yon_{\!_p}^{\!^{\ash_\tau}}\to \yon_{\!_q}^{\!^{\ash_\tau}}$. Then, in particular, for the section $\id_{_p}^{\ash}\in \yon_{\!_p}^{\!^{\ash_\tau}}(p)$, one has $\big(f_{0,\ast}(\id_{_p})\big)^{\ash}=\big(f_{1,\ast}(\id_{_p})\big)^{\ash}$, and hence there exists a $\tau$-cover $\sigma:z\to p$ such that $f_0\circ \sigma=\sigma^\ast \big(f_{0,\ast}(\id_{_p})\big)=\sigma^\ast \big(f_{1,\ast}(\id_{_p})\big)=f_1\circ \sigma\in \yon_{\!_q}(z)$. Since $p$ is reduced and $\tau$-covers are surjective, a diagram chase shows that $\sigma$ is an epimorphism in the category of (proper) $S$-schemes, and hence $f_0=f_1$. 
\end{proof}

\Prop{{\cite[Prop.3.2.5]{Voe:96:HS}}}{
	Assume that $\tau$ is an additively-saturated pretopology on (proper) $S$-schemes, and let $f:p\to q$ be a morphism of (proper) $S$-schemes. Then, 
	\begin{enumerate}
		\item\label{itm:Classifiaction:Proposition:e} the morphism $f_\ast:\yon_{\!_p}^{\!^{\ash_\tau}}\to \yon_{\!_q}^{\!^{\ash_\tau}}$ is an epimorphism \iff $f$ is a cover in the saturation of $\tau$; and
		\item\label{itm:Classifiaction:Proposition:m} assuming that $\tau$ is finer than the closed pretopology and that $\tau$-covers are surjective, the morphism $f_\ast:\yon_{\!_p}^{\!^{\ash_\tau}}\to \yon_{\!_q}^{\!^{\ash_\tau}}$ is a monomorphism \iff $f$ is universally injective.
	\end{enumerate}
}{Voevodsky:96:HS:Prop.3.2.5}
\begin{proof}\ 
	\begin{enumerate}
		\item Assume that $f$ is a cover in the saturation of $\tau$, \ie there exists a morphism $\sigma':p'\to p$ of (proper) $S$-schemes such that $\sigma\coloneqq f\circ \sigma'$ is a $\tau$-cover. The morphism $\sigma_\ast: \yon_{\!_{p'}}\to \yon_{\!q}$ is a $\tau$-local epimorphism because it factorises as an epimorphism $\yon_{\!_{p'}}\to \im \sigma_\ast$ followed by the inclusion $\im \sigma_\ast \subset \yon_{\!q}$ of the $\tau$-covering sieve generated by $\sigma$. Thus, the morphism $f_\ast:\yon_{\!_p}\to \yon_{\!_q}$ is a $\tau$-local epimorphism, and hence the morphism $f_\ast:\yon_{\!_p}^{\!^{\ash_\tau}}\to \yon_{\!_q}^{\!^{\ash_\tau}}$ is an epimorphism of $\tau$-sheaves.
		
		\medskip
		
		On the other hand, assume that $f_\ast:\yon_{\!_p}^{\!^{\ash_\tau}}\to \yon_{\!_q}^{\!^{\ash_\tau}}$ is an epimorphism of $\tau$-sheaves, \ie $f_\ast:\yon_{\!_p}\to \yon_{\!_q}$ is a $\tau$-local epimorphism. For $\id_{\!_q}\in \yon_{\!_q}(q)$, there exists a $\tau$-cover $\sigma:w\to q$ and a section $a\in \yon_{\!_p}(w)$ such that $		f\circ a=f_\ast({a})=\sigma^\ast(\id_{\!_q})=\sigma$, by Lemma \ref{Lem:Local_Morphism}. Thus, the morphism $f$ is a cover in the saturation of $\tau$.
		
		\item The proof of the \textit{if} implication essentially follows \cite{And:17:ERHS}, which corrects a mistake in the proof of \cite[Prop.3.2.5.(i)]{Voe:96:HS}.
		
		\medskip
		
		Assume that $f:p\to q$ is universally injective, let $z$ be a (proper) $S$-scheme, and let ${a}_0$ and ${a}_1$ be sections in $\yon_{\!_p}(z)$ such that $f\circ {a}_0=f_\ast({a}_0)=f_\ast ({a}_1)=f\circ {a}_1$. Consider the commutative solid diagram
		\[
		\begin{tikzpicture}[descr/.style={fill=white},ampersand replacement=\&]
		\node (m-2-2) at (0,0) {$p\times_{q}p$};
		\node (m-1-1) at (-1,1) {$z$};
		\node (m-2-3) at (2,0) {$p$};
		\node (m-3-2) at (0,-2) {$p$};
		\node (m-3-3) at (2,-2) {$q$};
		\node at (.25,-.4) {$\ulcorner$};
		\path[bend left, ->,font=\scriptsize]
		(m-1-1) edge node [above] {$a_0$} (m-2-3)
		;
		\path[bend right, ->,font=\scriptsize]
		(m-1-1) edge node [left] {$a_1$} (m-3-2)
		;
		\path[->,font=\scriptsize]
		(m-3-2) edge node [below] {$f$} (m-3-3)
		(m-2-2) edge (m-2-3)
		(m-2-2) edge (m-3-2)
		(m-2-3) edge node [right] {$f$} (m-3-3);
		\path[dotted, ->,font=\scriptsize]
		(m-1-1) edge node {$\lambda$} (m-2-2)
		;
		\end{tikzpicture}
		\]
		of (proper) $S$-schemes, and let $\lambda:z\to p\times_q p$ be the unique such morphism that makes the whole diagram commute. Since $f:p\to q$ is universally injective, the diagonal morphism $\Delta_{_f}:p \to p\times_q p$ is a surjective closed immersion, by \cite[Tag 01S4]{stacks-project}. Let $i:z_{\!_\red}\cim z$ be the close immersion of the maximal reduced closed subscheme in $z$. The morphism $i$ is a $\tau$-cover, as the $\tau$-pretopology is finer than the closed pretopology. The morphism $\lambda\circ i$ factories through every surjective closed immersion to $p\times_q p$, in particular, it factorises through the diagonal morphism $\Delta_{_f}$, which implies that $a_0\circ i=a_1\circ i$. Thus, the morphism $f_\ast:\yon_{\!_p}\to \yon_{\!_q}$ is a $\tau$-local monomorphism, by Lemma \ref{Lem:Local_Morphism}, and hence the morphism $f_\ast:\yon_{\!_p}^{\!^{\ash_\tau}}\to \yon_{\!_q}^{\!^{\ash_\tau}}$ is a monomorphism of $\tau$-sheaves.

		\medskip

		On the other hand, assume that $f_\ast:\yon_{\!_p}^{\!^{\ash_\tau}}\to \yon_{\!_q}^{\!^{\ash_\tau}}$ is a monomorphism of $\tau$-sheaves, and consider the commutative diagram 
		\begin{equation}\label{Monomorphism:Universally:Injective}
		\begin{tikzpicture}[descr/.style={fill=white},ampersand replacement=\&]
		\node (m-2-2) at (0,0) {$p\times_{q}p$};
		\node (m-1-1) at (-1,1) {$p$};
		\node (m-2-3) at (2,0) {$p$};
		\node (m-3-2) at (0,-2) {$p$};
		\node (m-3-3) at (2,-2) {$q$};
		\node at (.25,-.4) {$\ulcorner$};
		\path[bend left, ->,font=\scriptsize]
		(m-1-1) edge node [above] {$\id_{_p}$} (m-2-3)
		;
		\path[bend right, ->,font=\scriptsize]
		(m-1-1) edge node [left] {$\id_{_p}$} (m-3-2)
		;
		\path[->,font=\scriptsize]
		(m-3-2) edge node [below] {$f$} (m-3-3)
		(m-2-2) edge node [below] {$\pi_1$}(m-2-3)
		(m-2-2) edge node [right] {$\pi_0$} (m-3-2)
		(m-2-3) edge node [right] {$f$} (m-3-3)
		(m-1-1) edge node[descr]{$\!\!\!\Delta_{_f}\!\!\!$} (m-2-2)
		;
		\end{tikzpicture}
		\end{equation}
		of (proper) $S$-schemes. Recall that both the \yoneda embedding and the $\tau$-sheafification functor preserve finite limits. In particular, the morphisms $\pi_{0,\ast}$ and $\pi_{1,\ast}$ are base changes of $f_\ast$ along itself, and hence they are monomorphisms of $\tau$-sheaves. In fact, the morphisms $\pi_{0,\ast}$ and $\pi_{1,\ast}$ are isomorphisms of $\tau$-sheaves, as $\id_{_p,\ast}$ is an epimorphism of $\tau$-sheaves and the category of $\tau$-sheaves of sets is a balanced category. Thus, $\Delta_{f,\ast}$ is an epimorphism of $\tau$-sheaves, and hence $\Delta_{f}$ is a cover in the saturation of $\tau$, by \eqref{itm:Classifiaction:Proposition:e}. In particular, $\Delta_{f}$ is surjective, and hence $f$ is universally injective, by \cite[Tag 01S4]{stacks-project}.
	\end{enumerate}
\end{proof}

\Cor{}{Let $\tau$ be an additively-saturated pretopology on the category of (proper) $S$-schemes that is finer than the closed pretopology and coarser than the proper pretopology, and let $f:p\to q$ be a morphism of (proper) $S$-schemes. Then, the morphism $f_\ast:\yon_{\!_p}^{\!^{\ash_\tau}}\to \yon_{\!_q}^{\!^{\ash_\tau}}$ is an isomorphism only if the morphism $f$ is a universal homeomorphism.
}{Topology:Isomorphism:Only:If}
\begin{proof}Assume that $f_\ast:\yon_{\!_p}^{\!^{\ash_\tau}}\to \yon_{\!_q}^{\!^{\ash_\tau}}$ is an isomorphism, then $f$ is a universally injective cover in the saturation of $\tau$, by Proposition \ref{Prop:Voevodsky:96:HS:Prop.3.2.5}. In particular, there exists a morphism $\sigma':p'\to p$ of (proper) $S$-schemes such that $\sigma\coloneqq f\circ \sigma'$ is a $\tau$-cover. Since $\tau$ is coarser than the proper pretopology, the morphism $\sigma$ is surjective and universally closed, and hence a universal topological epimorphism\footnote{A \indx{\textup{(}universal\textup{)} topological epimorphism} $f$ is a morphism of schemes for which the underlying continuous map (of every base change) of $f$ is a quotient map, see \cite[\S.3.1]{Voe:96:HS}.}. This implies that $f$ is also a universal topological epimorphism, and hence every base change of $f$ is both an injection and a topological epimorphism. That is the underlying continuous map of every base change of $f$ is a monomorphism and an extremal epimorphism in the category of topological spaces, and hence a homeomorphism, see \cite[\S.2.6-\S.2.9]{Nak:89:C14CT}. Therefore, $f$ is a universal homeomorphism.
\end{proof}

\Rem{}{In the sequel, we restrict our attention to additively-saturated pretopologies on the category of proper $S$-schemes that are finer than the \cdf-pretopology and coarser than the {proper} pretopology.
}{Topology:Restrict:Attention}
\Eg{}{Pretopologies on the category of proper $S$-schemes that satisfy the assumptions of Remark \ref{Rem:Topology:Restrict:Attention} include the finite pretopology, the \cdf-pretopology, the {proper} pretopology, and the \cdp-pretopology.
}{}

\Prop{}{Let $\tau$ be a pretopology on $\PropSch/S$ as in Remark \ref{Rem:Topology:Restrict:Attention}, and let $f:p\to q$ be a morphism of proper $S$-schemes. Then, the morphism $f_\ast:\yon_{\!_p}^{\!^{\ash_\tau}}\to \yon_{\!_q}^{\!^{\ash_\tau}}$ is an isomorphism \iff $f$ is a universal homeomorphism.
}{Topology:Isomorphism:Proper:If:Only:If}

\begin{proof}Since the \cdf-pretopology is finer than the closed pretopology, the \textit{only if} implication is the statement of Corollary \ref{Cor:Topology:Isomorphism:Only:If}.
	
	\medskip
	
	Assume that $f$ is a universal homeomorphism. Then, $f$ is a surjective universally injective finite morphism, by \cite[Prop.2.4.5]{EGA4.2}. In particular, for every field $\fk$, the induced map
	\[
	f_\ast:\fdNoe \big(\Spec \fk, P\big)\to \fdNoe\big(\Spec \fk, Q\big)
	\]
	is an injection, by \cite[Tag 01S4]{stacks-project}, where $P$ and $Q$ are the underlying schemes for $p$ and $q$, respectively. 
	
	\medskip
	
	For a field $\fk$, let $y:\Spec \fk \to Q$ be a morphism of schemes, and consider the Cartesian square
	\[
	\begin{tikzpicture}[descr/.style={fill=white},ampersand replacement=\&]
	\node (m-1-1) at (0,0) {$Z$};
	\node (m-1-2) at (2,0) {$P$};
	\node (m-2-1) at (0,-2) {$\Spec \fk$};
	\node (m-2-2) at (2,-2) {$Q$};
	\node at (.25,-.4) {$\ulcorner$};
	\path[->,font=\scriptsize]
	(m-1-1) edge node[left]{$\underline{f}$} (m-2-1)
	(m-1-2) edge node[right]{$f$} (m-2-2)
	(m-2-1) edge node[below]{$y$} (m-2-2)
	(m-1-1) edge node[above]{$\underline{y}$} (m-1-2);	
	\end{tikzpicture}
	\]
	in the category $\fdNoe$. The morphism $f$ is a finite universal homeomorphism, and hence $Z$ is a one-point scheme $\Spec R$ and $\underline{f}$ is induced by a finite ring homomorphism $\psi:\fk \incl R$, to a local ring $R$ of Krull dimension zero. Let $\frm$ be the maximal ideal of $R$, and let $\rf\coloneqq \sfrac{R}{\frm}$. Then, the induced homomorphism $\fk \incl \rf$ is a finite field extension. Assuming that $[\rf:\fk]\neq 1$, there exist distinct ring homeomorphisms $\rf\to \rf$ over $\fk$, which contradicts with $\underline{f}$ being universally injective. Thus, one has $[\rf:\fk]=1$, \ie the residue field of $Z$ at its unique point is isomorphic to $\fk$. Hence, $y$ lifts along $f$, and $f_\ast$ is surjective for every field $\fk$. Therefore, $f$ is a \cdf-cover that is universally injective, and hence a universally injective $\tau$-cover. Therefore, the morphism $f_\ast:\yon_{\!_p}^{\!^{\ash_\tau}}\to \yon_{\!_q}^{\!^{\ash_\tau}}$ is an isomorphism, by Proposition \ref{Prop:Voevodsky:96:HS:Prop.3.2.5}.
\end{proof}

\Eg{{\cite[Lem.3.2.1]{Voe:96:HS}}}{Let $\tau$ be a pretopology on $\PropSch/S$ as in Remark \ref{Rem:Topology:Restrict:Attention}, and let $i:z\cim p$ be a surjective closed immersion of proper $S$-schemes. Then, the morphism $i_\ast:\yon_{\!_z}^{\!^{\ash_\tau}}\to \yon_{\!_p}^{\!^{\ash_\tau}}$ is an isomorphism. In particular, for the closed immersion of the maximal reduced closed subscheme $i:p_{_\red}\cim p$, the morphism $i_\ast:\yon_{\!_{p_{_\red}}}^{\!^{\ash_\tau}}\to \yon_{\!_p}^{\!^{\ash_\tau}}$ is an isomorphism.
}{Voevodsky:96:HS:Lem.3.2.1.Cor}

\begin{small}
	\linespread{1}

\end{small}
\end{document}